\documentclass[12pt,twoside]{amsart}
\usepackage{amssymb,mathrsfs}

\usepackage{amssymb,mathrsfs}
\usepackage{tikz-cd}
\usepackage{dynkin-diagrams}

\usepackage{fancyhdr}

\usepackage{cite}

\newcommand{\SGrunninghead}[1]{{\scriptsize #1}}
\fancypagestyle{plain}{%
  \fancyhf{}%
  \fancyhead[LE,RO]{\SGrunninghead{\thepage}}%
}
\makeatletter
\renewcommand{\sectionmark}[1]{%
  \begingroup
    \edef\SG@next{\endgroup
      \noexpand\markboth{\thesection.\space \unexpanded{#1}}%
                         {\thesection.\space \unexpanded{#1}}}%
  \SG@next}
\renewcommand{\subsectionmark}[1]{%
}
\let\SG@orig@sect\@sect
\def\@sect#1#2#3#4#5#6[#7]#8{%
  \SG@orig@sect{#1}{#2}{#3}{#4}{#5}{#6}[#7]{#8}%
  \ifnum#2<\@m
    \@ifundefined{#1mark}{}{\csname #1mark\endcsname{#7}}%
  \fi
}

\def\l@section{\@tocline{1}{0pt}{1pc}{}{\sffamily}}
\def\l@subsection{\@tocline{2}{0pt}{1pc}{5pc}{\sffamily}}
\makeatother
\begin{document}

\title{Remarks on Symplectic Geometry}

\author{Jae-Hyun Yang}

\address{Yang Institute for Advanced Study
\newline\indent
Seoul 07989, Republic of Korea
\vskip 2mm
and
\vskip 2mm
%}
%\address{
Department of Mathematics
\newline\indent
Inha University
\newline\indent
Incheon 22212, Republic of Korea}
\email{jhyang@inha.ac.kr}

\newtheorem{theorem}{Theorem}[section]
\newtheorem{lemma}[theorem]{Lemma}
\newtheorem{corollary}[theorem]{Corollary}
\newtheorem{proposition}[theorem]{Proposition}
\newtheorem{remark}[theorem]{Remark}
\newtheorem{definition}[theorem]{Definition}
\newtheorem{conjecture}[theorem]{Conjecture}
\newtheorem{example}[theorem]{Example}
\newtheorem{exercise}[theorem]{Exercise}
\newtheorem{problem}[theorem]{Problem}
\newtheorem{question}[theorem]{Question}

\renewcommand{\theequation}{\thesection.\arabic{equation}}
\renewcommand{\thetheorem}{\thesection.\arabic{theorem}}
\renewcommand{\thelemma}{\thesection.\arabic{lemma}}
\newcommand{\bbr}{\mathbb R}
\newcommand{\bbs}{\mathbb S}
\newcommand{\bn}{\bf n}
\newcommand\charf {\mbox{{\text 1}\kern-.24em {\text l}}}
\newcommand\fg{{\mathfrak g}}
\newcommand\fk{{\mathfrak k}}
\newcommand\fp{{\mathfrak p}}
\newcommand\g{\gamma}
\newcommand\G{\Gamma}
\newcommand\ka{\kappa}
\newcommand\al{\alpha}
\newcommand\be{\beta}
\newcommand\lrt{\longrightarrow}
\newcommand\s{\sigma}
\newcommand\ba{\backslash}
\newcommand\lmt{\longmapsto}
\newcommand\CP{{\mathscr P}_g}
\newcommand\CM{{\mathcal M}}
\newcommand\BC{\mathbb C}
\newcommand\BZ{\mathbb Z}
\newcommand\BR{\Bbb R}
\newcommand\BQ{\mathbb Q}
\newcommand\Rmn{{\mathbb R}^{(m,n)}}
\newcommand\PR{{\mathcal P}_n\times {\mathbb R}^{(m,n)}}
\newcommand\Gnm{GL_{n,m}}
\newcommand\Gnz{GL_{n,m}({\mathbb Z})}
\newcommand\Gjnm{Sp_{n,m}}
\newcommand\Gnml{GL(n,{\mathbb R})\ltimes {\mathbb R}^{(m,n)}}
\newcommand\Snm{SL_{n,m}}
\newcommand\Snz{SL_{n,m}({\mathbb Z})}
\newcommand\Snml{SL(n,{\mathbb R})\ltimes {\mathbb R}^{(m,n)}}
\newcommand\Snzl{SL(n,{\mathbb Z})\ltimes {\mathbb Z}^{(m,n)}}
\newcommand\la{\lambda}
\newcommand\GZ{GL(n,{\mathbb Z})\ltimes {\mathbb Z}^{(m,n)}}
\newcommand\DPR{{\mathbb D}(\PR)}
\newcommand\Rnn{{\mathbb R}^{(n,n)}}
\newcommand\Yd{{{\partial}\over {\partial Y}}}
\newcommand\Vd{{{\partial}\over {\partial V}}}
\newcommand\Ys{Y^{\ast}}
\newcommand\Vs{V^{\ast}}
\newcommand\DGR{{\mathbb D}(\Gnm)}
\newcommand\DKR{{\mathbb D}_K(\Gnm)}
\newcommand\DKS{{\mathbb D}_{K_0}(\Snm)}
\newcommand\fa{{\frak a}}
\newcommand\fac{{\frak a}_c^{\ast}}
\newcommand\SPR{S{\mathcal P}_n\times \Rmn}
\newcommand\DSPR{{\mathbb D}(\SPR)}
\newcommand\BD{{\mathbb D}}
\newcommand\SP{{\mathfrak P}_n}
\newcommand\BH{{\mathbb H}}

\thanks{
\endgraf
Keywords and phrases\,: classical mechanics, symplectic geometry, pseudoholomorphic curves,
\\ \indent symplectic topology, Arnold's conjecture, moment map, symplectic group actions, symplectic
\\ \indent embeddings, Gromov-Witten invariants, Floer homology, Fukaya category, mirror symmetry}

\begin{abstract}
The relatively-new subject, symplectic geometry, has been studied the past five decades.
Symplectic geometry is the mathematical subject studying the geometry of symplectic manifolds.
A symplectic manifold is an even-dimensional smooth manifold
equipped with a closed non-degenerate two form. In this paper,
we survey the progress on the study of symplectic geometry the past five decades.
We deal with the brief history of symplectic geometry and symplectic topology,
Arnold's conjecture, pseudoholomorphic curves,
the convexity properties of a moment map, recent results on Hamiltonian and non-Hamiltonian
symplectic group actions, the classification of symplectic group actions, recent results on the symplectic embedding problems, the theory of Gromov-Witten invariants, Lagrangian Floer homology,
the Fukaya category and mirror symmetry.
\end{abstract}

\maketitle

\vskip 8mm

\tableofcontents

\vskip 10mm

%%%%%%%%%%%%%%%%%%%%%%%%%%%%%%%%%%%%%%%%%%%%%%%%%%%%%%%%%%%%%%%%%%%%%%%%%%%%%%%%%%%%%%%%%%%%%%%%%%%%%%%%%%%%%%%%%%%%%%%%%%%%%%%%%%%%%%%%%%%%%%%%%%%%%%%%
%%%%%%%%%%%%%%%%%%%%%%%%%%%%%%%%%%%%%%%%%%%%%%%%%%%%%%%%%%%%%%%%%%%%%%%%%%%%%%%%%%%%%%%%%%%%%%%%%%%%%%%%%%%%%%%%%%%%%%%%%%%%%%%%%%%%%%%%%%%%%%%%%%%%%%%%
%%%%%%%%%%%%%%%%%%%%%%%%%%%%%%%%%%%%%%%%%%%%%%%%%%%%%%%%%%%%%%%%%%%%%%%%%%%%%%%%%%%%%%%%%%%%%%%%%%%%%%%%%%%%%%%%%%%%%%%%%%%%%%%%%%%%%%%%%%%%%%%%%%%%%%%%
%%%%%%%%%%%%%%%%%%%%%%%%%%%%%%%%%%%%%%%%%%%%%%%%%%%%%%%%%%%%%%%%%%%%%%%%%%%%%%%%%%%%%%%%%%%%%%%%%%%%%%%%%%%%%%%%%%%%%%%%%%%%%%%%%%%%%%%%%%%%%%%%%%%%%%%%
%%%%%%%%
%%%%%%%%
%%%%%%%%                   Section 1  Introduction
%%%%%%%%
%%%%%%%%
%%%%%%%%%%%%%%%%%%%%%%%%%%%%%%%%%%%%%%%%%%%%%%%%%%%%%%%%%%%%%%%%%%%%%%%%%%%%%%%%%%%%%%%%%%%%%%%%%%%%%%%%%%%%%%%%%%%%%%%%%%%%%%%%%%%%%%%%%%%%%%%%%%%%%%%%
%%%%%%%%%%%%%%%%%%%%%%%%%%%%%%%%%%%%%%%%%%%%%%%%%%%%%%%%%%%%%%%%%%%%%%%%%%%%%%%%%%%%%%%%%%%%%%%%%%%%%%%%%%%%%%%%%%%%%%%%%%%%%%%%%%%%%%%%%%%%%%%%%%%%%%%%
%%%%%%%%%%%%%%%%%%%%%%%%%%%%%%%%%%%%%%%%%%%%%%%%%%%%%%%%%%%%%%%%%%%%%%%%%%%%%%%%%%%%%%%%%%%%%%%%%%%%%%%%%%%%%%%%%%%%%%%%%%%%%%%%%%%%%%%%%%%%%%%%%%%%%%%%
%%%%%%%%%%%%%%%%%%%%%%%%%%%%%%%%%%%%%%%%%%%%%%%%%%%%%%%%%%%%%%%%%%%%%%%%%%%%%%%%%%%%%%%%%%%%%%%%%%%%%%%%%%%%%%%%%%%%%%%%%%%%%%%%%%%%%%%%%%%%%%%%%%%%%%%%

\clearpage
\begin{section}{{\bf Introduction}}
\setcounter{equation}{0}
\vskip 3mm
Symplectic geometry is a branch of mathematics that arose from classical mechanics and
studies the geometry of a symplectic manifold.
A symplectic manifold $(M,\omega)$ is an even-dimensional smooth manifold
equipped with a {\it closed non-degenerate} two form $\omega$. A symplectic form
or a symplectic structure $\omega$ is quite a mysterious, subtle and intricate thing.
Darboux's theorem states that every symplectic structure on a symplectic manifold $(M,\omega)$ of
dimension $2\,n$ is locally diffeomorphic to the standard symplectic structure on
$(\BR^{2n},\omega_0)$ given by
\begin{equation*}
 \omega_0=\sum_{j=1}^n dx_j\wedge dy_j\qquad {\rm on}\ \BR^{2n}.
\end{equation*}
According to Darboux's theorem, there are no local invariant results.
Symplectic geometry is concerned with the study of a notion of signed {\it area}
rather than length, distance or volume. It can be less intuitive than Riemannian geometry.
The diffeomorphisms that preserve the symplectic structure are called the symplectomorphism.
We are familiar with Riemannian geometry, in particular Euclidean geometry where the
curvature tensor gives a whole field of invariants.
So because of the existence of a
symplectic structure, symplectic geometry is in sharp contrast to Riemannian geometry.
Indeed, symplectic structures give rise to an abundance of interesting and deep problems
which are important geometrically and topologically.
In the world of symplectic geometry, it took many years for mathematicians to understand
the importance of a symplectic structure in the classical mechanics (Hamiltonian system) since
the Lagrangian-Jacobian-Hamiltonian ages.
\vskip 2mm
In the late 1960s, Arnold\,(1937--2010), Marsden and Weinstein were stimulated by the
innovative and ground-breaking works of Souriau\,(1922--2012) and began to study symplectic
manifolds linking to classical mechanics\,(cf.\,\cite{A-M,A1,A2,A3, M-W, W1,W-1,W-2,W2,W3}).
In the early 1980s, the convexity properties of
a moment map were investigated by Atiyah, Guillemin, Sternberg, Mumford and
Kirwan\,(cf.\,\cite{A, A-B, Gu-St1, Gu-St2, Ki, Mc1, M-F-K}).
In 1985 Gromov introduced the pseudoholomorphic (or J-holomorphic) curve technique and the symplectic capacity into symplectic geometry
to prove the famous {\sf Nonsqueezing\ Theorem} in his seminal paper,
{\sf ``Pseudoholomorphic curves in symplectic manifolds"}\,\cite{Gr1}. Thereafter
Arnold's conjecture and the concept of
pseudoholomorphic curves motivated the occurrence of the new mathematical subject, {\sf symplectic\ topology} which has been developed by the McDuff school \cite{Mc2, Mc5, Mc6} et al. Pseudoholomorphic curves are used as a tool in the four-dimensional symplectic topology, the study of the symplectic embeddings,
the theory of ECH (embedded contact homology),
and in the theory of the Gromov-Witten invariants. Taubes\,\cite{Ta1, Ta2} proved that in four dimensional symplectic manifolds,
the Gromov-Witten invariants coincide with the Seiberg-Witten invariants. The theory of Floer homology (cf.\,\cite{F1, F1+, Z}) is based on pseudoholomorphic curves with boundary lying on a Lagrangian submanifold. This led on to the notion of the Fukaya category \cite{Fu, FO1, FO2}
which is closely related to the homological mirror symmetry formulated by Kontsevich (cf.\,\cite{Ko}).
Later the symplectic capacity which is a notion of monotonic symplectic invariant was pioneered by Ekeland and Hofer\,\cite{E-H}, and then has been developed by Hofer and his collaborators (cf.\,\cite{Ho, H-Z}) from the angle of dynamical systems and Hamiltonian dynamics.
After the middle 1980s, due to the
important works of a convexity of a moment map \cite{A, A-B,Gu,Gu-St1,Gu-St2,Ki}, the study of
Hamiltonian and non-Hamiltonian symplectic group actions has been actively and intensively
developed by many mathematicians. After the early 1990s, due to Gromov's seminal works on nonsqueezing theorem and pseudoholomorphic curves, the study of symplectic embeddings has been implemented
intensively by McDuff, Hutchings et al.
The past five decades
there have been many major developments on symplectic geometry, for example, classical mechanics, statistical mechanics, quantum mechanics, classification of symplectic and Hamiltonian group actions, symplectic embedding problems, the theory of Gromov-Witten invariants, quantum cohomolgy \cite{K-M}, symplectic capacities, the ECH theory \cite{Hu1,Hu2,Hu3,Hu4}, the theory of Floer homology,
Fukaya category and homological mirror symmetry etc.
Nowadays symplectic geometry and symplectic topology became important subjects in the mathematical community. In particular, recently
Floer homology, the Fukaya category and homological mirror symmetry are major research topics
in the areas of symplectic geometry and symplectic topology.
Covering these three topics would be beyond the scope of this paper. However we give a brief
description of these three topics in the subsection 6.4.
\vskip 2mm
McDuff and Salamon dealt with symplectic topology and pseudoholomorphic curves in details
in their books \cite{Mc5, Mc6}. Cannas da Silva \cite{C}, Pelayo \cite{P2} and Schlenk \cite{Sc}
surveyed the convexity of a moment map, symplectic group actions and symplectic embedding
problems respectively in some details. In this paper, we give a brief overview of symplectic
topology, pseudoholomorphic curves, the convexity of a moment map, symplectic group actions
and symplectic embedding problems following the references \cite{Mc5, Mc6, C, P2, Sc} and
update some recent results concerning symplectic group actions \cite{J-T}
and symplectic embedding \cite{Cr-1, Cr-2}.

\vskip 2mm
The purpose of this article is to survey various results related to symplectic geometry
the past five decades. This article is organized as follows.
%%%%%%%%%%%%%%%%%%%%%%%%%%%%%%%%%%%%%%%%%%%%%%%%%%%%%%%%%%%%%%%%%%%%%%%%%%%%%%%%%%%%%%%%%%%%%%%
%%%%%%%%%%%%%%%%%%%%%%%%%%%%%%%%%%%%%%%%%%%%%%%%%%%%%%%%%%%%%%%%%%%%%%%%%%%%%%%%%%%%%%%%%%%%%%%
In section 2, we give a brief description of the history of
``{\sf symplectic geometry}". We review briefly Siegel's work on the geometry of the Siegel
upper half plane that is a symplectic manifold (cf.\,\cite{Si}) and Souriau's influence on
symplectic geometry.
%%%%%%%%%%%%%%%%%%%%%%%%%%%%%%%%%%%%%%%%%%%%%%%%%%%%%%%%%%%%%%%%%%%%%%%%%%%%%%%%%%%%%%%%%%%%%%%
%%%%%%%%%%%%%%%%%%%%%%%%%%%%%%%%%%%%%%%%%%%%%%%%%%%%%%%%%%%%%%%%%%%%%%%%%%%%%%%%%%%%%%%%%%%%%%%
In section 3, we give basic notions, definitions and examples of symplectic manifolds.
We also review some results which are Moser's theorem,
Marsden-Weinstein-Meyer theorem, Duistermaat-Heckman theorem, Weinstein Lagrangian
neighborhood theorem, and Weinstein tubular neighborhood theorem.
%%%%%%%%%%%%%%%%%%%%%%%%%%%%%%%%%%%%%%%%%%%%%%%%%%%%%%%%%%%%%%%%%%%%%%%%%%%%%%%%%%%%%%%%%%%%%%%
%%%%%%%%%%%%%%%%%%%%%%%%%%%%%%%%%%%%%%%%%%%%%%%%%%%%%%%%%%%%%%%%%%%%%%%%%%%%%%%%%%%%%%%%%%%%%%%
In section 4, we discuss Arnold's conjecture and pseudoholomorphic curves. We mention
Poincar{\'e}-Birkhoff theorem, Arnold's conjecture, Gromov's nonsqueezing theorem,
symplectic rigidity and Donaldson's work on symplectic submanifolds.
%%%%%%%%%%%%%%%%%%%%%%%%%%%%%%%%%%%%%%%%%%%%%%%%%%%%%%%%%%%%%%%%%%%%%%%%%%%%%%%%%%%%%%%%%%%%%%%
%%%%%%%%%%%%%%%%%%%%%%%%%%%%%%%%%%%%%%%%%%%%%%%%%%%%%%%%%%%%%%%%%%%%%%%%%%%%%%%%%%%%%%%%%%%%%%%
In section 5, we review the convexity properties of a moment map obtained by Atiyah, Guillemin,
Sternberg and Kirwan (cf.\,\cite{A, A-B, Gu, Gu-St1, Gu-St2, Ki}). Using the convexity
properties of a moment map, Delzant \cite{De} classified all symplectic-toric manifolds in terms
of a set of very special polytopes, so-called the Delzant polytopes. We state the Delzant's work.
%%%%%%%%%%%%%%%%%%%%%%%%%%%%%%%%%%%%%%%%%%%%%%%%%%%%%%%%%%%%%%%%%%%%%%%%%%%%%%%%%%%%%%%%%%%%%%%
%%%%%%%%%%%%%%%%%%%%%%%%%%%%%%%%%%%%%%%%%%%%%%%%%%%%%%%%%%%%%%%%%%%%%%%%%%%%%%%%%%%%%%%%%%%%%%%
In the final section we survey various results related to symplectic geometry and symplectic
topology obtained during the past five decades.
%%%%%%%%%%%%%%%%%%%%%%%%%%%%%%%%%%%%%%%%%%%%%%%%%%%%%%%%%%%%%%%%%%%%%%%%%%%%%%%%%%%%%%%%%%%%%%%%%%%%%%
In the subsection 6.1, we review on the classification problem of symplectic
and Hamiltonian group actions, and in the subsection 6.2 we survey various results on the symplectic
embedding problems. In the subsection 6.3, the theory of Gromov-Witten invariants is discussed, and
there is a mention about the work of Shen and Zhou on the Landau-Ginzburg/Calabi-Yau correspondence for
elliptic orbifold projective lines. In the subsection 6.4, we give a very sketchy description
of Lagrangian Floer homologies, the Fukaya category and homological mirror symmetry.

\vskip 0.31cm \noindent {\bf Notations:} \ \ We denote by
$\BQ,\,\BR$ and $\BC$ the field of rational numbers, the field of
real numbers and the field of complex numbers respectively. We
denote by $\BZ$ and $\BZ^+$ the ring of integers and the set of
all positive integers respectively. $\BQ^+$ (resp. $\BR^+$) denotes the set of all positive
rational (resp. real) numbers. We denote by $\BZ_+$ (resp. $\BQ_+,\ \BR_+$) the set of all
non-negative integers (resp. rational numbers, real numbers).
$\BQ^{\times}$ (resp. $\BR^{\times},\ \BC^{\times}$)
denotes the group of nonzero rational (resp. real, complex) numbers.
The symbol ``:='' means that
the expression on the right is the definition of that on the left.
For two positive integers $k$ and $l$, $F^{(k,l)}$ denotes the set
of all $k\times l$ matrices with entries in a commutative ring
$F$. For a square matrix $A\in F^{(k,k)}$ of degree $k$,
${\rm Tr}(A)$ denotes the trace of $A$. For any $B\in F^{(k,l)},\
^t\!B$ denotes the transpose of $B$.
For $A\in F^{(k,l)}$ and $B\in F^{(k,k)}$, we set $B[A]=\,^tABA$ (Siegel's notation).
For a positive integer $g$, $I_g$ denotes the identity matrix of degree $g$.
For a complex matrix $A$,
${\overline A}$ denotes the complex {\it conjugate} of $A$.
${\rm diag}(a_1,\ldots,a_g)$ denotes the $g\times g$ diagonal matrix with diagonal entries
$a_1,\ldots,a_g$. For a smooth manifold $M$, we denote by $C_c (M)$ (resp. $C_c^{\infty}(M)$ the algebra of all continuous (resp. infinitely differentiable) functions on $M$ with compact support, and by
${\mathfrak X}(M)$ the Lie algebra of all smooth vector fields on $M$.
$\mathbb{H}$ denotes the Poincar{\'e} upper half plane and ${\mathbb D}$ denotes the Poincar{\'e} disk.
If $z\in\BC,$ we put $e(z):=e^{2\pi iz}$. The contraction $\imath_X$ of a $k$-form $\alpha$
with a vector field $X$ is defined to be the $(k-1)$-form given by
$$
\left( \imath_X\alpha\right)(X_1,\ldots,X_{k-1}):=\alpha (X,X_1,\ldots,X_{k-1}).
$$
For a positive integer $g$, we let
$$
O(g):=\left\{ A\in \BR^{(g,g)}\,|\ A\,{}^t\!A=\,{}^t\!A A=I_g\,\right\}
$$
be the orthogonal group of degree $g$ and let
$$
U(g):=\left\{ A\in \BC^{(g,g)}\,|\ A\,{}^t{\overline A}=\,{}^t{\overline A} A=I_g\,\right\}
$$
be the unitary group of degree $g$. For a $g\times g$ complex matrix $Z$, ${\rm Im}\,Z$
denotes the imaginary part of $Z$.
\end{section}

%%%%%%%%%%%%%%%%%%%%%%%%%%%%%%%%%%%%%%%%%%%%%%%%%%%%%%%%%%%%%%%%%%%%%%%%%%%%%%%%%%%%%%%%%%%%%%%%%%%%%%%%%%%%
%%%%%%%%%%%%%%%%%%%%%%%%%%%%%%%%%%%%%%%%%%%%%%%%%%%%%%%%%%%%%%%%%%%%%%%%%%%%%%%%%%%%%%%%%%%%%%%%%%%%%%%%%%%%
%%%%%%%%%%%%%%%%%%%%%%%%%%%%%%%%%%%%%%%%%%%%%%%%%%%%%%%%%%%%%%%%%%%%%%%%%%%%%%%%%%%%%%%%%%%%%%%%%%%%%%%%%%%%
%%%%%%%%%%%%%%%%%%%%%%%%%%%%%%%%%%%%%%%%%%%%%%%%%%%%%%%%%%%%%%%%%%%%%%%%%%%%%%%%%%%%%%%%%%%%%%%%%%%%%%%%%%%%
%%%%%%%%%%%%%%%%%%%%%%%%%%%%%%%%%%%%%%%%%%%%%%%%%%%%%%%%%%%%%%%%%%%%%%%%%%%%%%%%%%%%%%%%%%%%%%%%%%%%%%%%%%%%
%%%%%%%%%%%%%%%%%%%%%%%%%%%%%%%%%%%%%%%%%%%%%%%%%%%%%%%%%%%%%%%%%%%%%%%%%%%%%%%%%%%%%%%%%%%%%%%%%%%%%%%%%%%%
%%%%%%%%%%%%%%%%%%%%%%%%%%%%%%%%%%%%%%%%%%%%%%%%%%%%%%%%%%%%%%%%%%%%%%%%%%%%%%%%%%%%%%%%%%%%%%%%%%%%%%%%%%%%

\clearpage
\vskip 12mm
\begin{section}{{\bf The brief history of ``symplectic geometry"}}
\setcounter{equation}{0}
\vskip 3mm
First of all we briefly describe the history of the terminology ``Symplectic Geometry".
Hermann Weyl (1885--1955) published his famous book,
{\sf The Classical Group\,: Their Invariants and Representations} \cite{We} in 1939.
In his book, he proposed to change the term ``{\it complex group}" to the term
``{\it symplectic group}". The word ``{\it symplectic}" is the ancient Greek word
for ``{\it complex}". He described the following
remarks\,:
\vskip 2mm
{\footnotesize The name ``complex group" formerly advocated by me in allusion to line complexes, as
these are defined by the vanishing of antisymmetric bilinear forms, has become more
and more embarrassing through collision with the word ``complex" in the connotation
of complex number. I therefore propose to replace it by the corresponding Greek
adjective ``symplectic." Dickson calls the group the ``Abelian linear group" in homage
to Abel who first studied it. \cite[The footnote at p.\,165]{We}}

\vskip 3mm
Thereafter the term ``symplectic" to describe this group was accepted and used by the mathematical community. For a positive integer $g$ and a field $K$, the symplectic group $Sp(2g,K)$ of degree $g$
is defined by
\begin{equation*}
  Sp(2g,K):=\left\{ M\in K^{(2g,2g)}\,|\ {}^tMJ_gM=J_g\,\right\},
\end{equation*}
where
\begin{equation*}
J_g=\begin{pmatrix}
      \ 0 & I_g \\
      -I_g & 0
    \end{pmatrix}\in K^{(2g,2g)}
\end{equation*}
denotes the symplectic matrix of degree $g$.

\vskip 3mm
As far as I know, the terminology ``{\sf Symplectic Geometry}" appeared for the first time
in the article of Carl Ludwig Siegel (1896--1981), {\em Symplectic Geometry} \cite{Si}
which was published in 1943.
In that paper, Siegel studied the geometry of the Siegel upper half plane $\BH_g$ which is biholomorphic to the Hermitian symmetric space $Sp(2g,\BR)/U(g).$ He discovered the explicit fundamental domain with respect to the Siegel modular group $Sp(2g,\BZ).$ In fact, $\BH_g$ is a symplectic manifold of dimension $g(g+1)$ because it is a K{\"a}hler manifold. We note that $\BH_g$ is an Einstein-K{\"a}hler Hermitian symmetric space.
\newcommand\POB{ {{\partial}\over {\partial{\overline \Omega}}} }
\newcommand\PZB{ {{\partial}\over {\partial{\overline Z}}} }
\newcommand\PX{ {{\partial}\over{\partial X}} }
\newcommand\PY{ {{\partial}\over {\partial Y}} }
\newcommand\PU{ {{\partial}\over{\partial U}} }
\newcommand\PV{ {{\partial}\over{\partial V}} }
\newcommand\PO{ {{\partial}\over{\partial \Omega}} }
\newcommand\PZ{ {{\partial}\over{\partial Z}} }
\def\Om{\Omega}
\def\om{\omega}
\newcommand\CCF{{\mathfrak F}_g}
\vskip 2mm
Now we briefly outline Siegel's work on the geometry of $\BH_g$.
For a given fixed positive integer $g$, we let
$$
{\mathbb H}_g=\,\{\,\Omega\in \BC^{(g,g)}\,|\ \Om=\,^t\Om,\ \ \ \rm{Im}\,\Om>0\,\}
$$
be the Siegel upper half plane of degree $g$ and let
$$
Sp(2g,\BR)=\{ \alpha\in \BR^{(2g,2g)}\ \vert \ ^t\!\alpha J_g\alpha= J_g\ \}
$$
be the real symplectic group of degree $g$, where $F^{(k,l)}$ denotes
the set of all $k\times l$ matrices with entries in a commutative
ring $F$ for two positive integers $k$ and $l$, $^t\!\alpha$ denotes
the transpose matrix of a matrix $\alpha$ and
$$J_g=\begin{pmatrix} 0&I_g\\
                   -I_g&0\end{pmatrix}\in \BR^{(2g,2g)}.$$
Then $Sp(2g,\BR)$ acts on $\BH_g$ transitively by
\begin{equation}\label{(2.1)}
\alpha\cdot\Om=(A\Om+B)(C\Om+D)^{-1},
\end{equation} where $\alpha=\begin{pmatrix} A&B\\
C&D\end{pmatrix}\in Sp(2g,\BR)$ and $\Om\in \BH_g.$ Let
$$\G_g=Sp(2g,\BZ)=\left\{ \begin{pmatrix} A&B\\
C&D\end{pmatrix}\in Sp(2g,\BR) \,\big| \ A,B,C,D\
\rm{integral}\ \right\}$$ be the Siegel modular group of
degree $g$. This group acts on $\BH_g$ properly discontinuously.
C. L. Siegel investigated the geometry of $\BH_g$ and automorphic
forms on $\BH_g$ systematically. Siegel\,\cite{Si} found a
fundamental domain ${\mathfrak F}_g$ for $\BH_g$ with respect to $\G_g$ and
described it explicitly. Moreover he calculated the volume of
$\CCF.$ Let ${\mathcal A}_g=\G_g\ba\BH_g$ be
the Siegel modular variety of degree $g$.
In fact, ${\mathcal A}_g$ is one of the important arithmetic varieties
in the sense that it is regarded as the moduli of principally
polarized abelian varieties of dimension $g$. Walter Baily\,(1930--2013) proved
that the Satake compactification of ${\mathcal A}_g$ is a normal projective
variety (cf.\,\cite{B1}).
\vskip 0.21cm
For $\Om=(\om_{ij})\in\BH_g,$ we write $\Om=X+i\,Y$
with $X=(x_{ij}),\ Y=(y_{ij})$ real. We put $d\Om=(d\om_{ij})$ and $d{\overline\Om}=(d{\overline\om}_{ij})$. We also put
$$
\PO=\,\left(\,
{ {1+\delta_{ij}}\over 2}\, { {\partial}\over {\partial \om_{ij} }
} \,\right) \qquad\text{and}\qquad \POB=\,\left(\, {
{1+\delta_{ij}}\over 2}\, { {\partial}\over {\partial {\overline
{\om}}_{ij} } } \,\right).
$$
The Bergman metric $ds_g^2$ on $\BH_g$ which is a $Sp(2g,\BR)$-invariant K{\"a}hler metric is given by
\begin{equation}\label{(2.2)}
  ds_g^2={\rm Tr} (Y^{-1}d\Om\, Y^{-1} d{\overline \Omega})=\sum_{1\leq i\leq j\leq g \atop
  1\leq k\leq l\leq g} h_{[ij][kl]}\,d\om_{ij}\, d\overline{\om}_{kl}
\end{equation}
and its K{\"a}hler form is
\begin{equation*}
  \omega_g=\sum_{1\leq i\leq j\leq g \atop
  1\leq k\leq l\leq g} h_{[ij][kl]}\,d\om_{ij}\wedge d\overline{\om}_{kl}.
\end{equation*}
Put $N=\frac{g(g+1)}{2}$. The volume form is given by
\begin{equation*}
  \frac{\om_g^N}{N!} = \left( \frac{i}{2}\right)^N \det (h_{[ij][kl]})\,dv_g
   =\left( \frac{i}{2}\right)^N \,(\det Y)^{-(g+1)}dv_g,
\end{equation*}
where
\begin{equation*}
  dv_g=2^{\frac{g(g-1)}{2}} \bigwedge_{1\leq i\leq j\leq g}d\om_{ij}\wedge d\overline{\om}_{kl}.
\end{equation*}
The invariant metric $ds_g^2$ is K{\"a}hler-Einstein, that is,
\begin{equation*}
\omega_g =\frac{2\,i}{g+1} \partial {\overline{\partial}} \log \det (h_{[ij][kl]}) =
-\frac{2\,i}{g+1} \partial \overline{\partial} \log (\det Y)^{g+1}.
\end{equation*}
The function $\Phi (\Omega):=\frac{4}{g+1}\partial {\overline{\partial}} \log \det (h_{[ij][kl]})$ is the
K{\"a}hler potential of $\omega_g$.

\vskip2mm
Siegel described the geodesics in $\BH_g$ and calculated the distance between two points
in $\BH_g$ explicitly (cf. \cite[pp.\,289--293]{Si}).

\newcommand\OW{\overline{W}}
\newcommand\OP{\overline{P}}
\newcommand\OQ{\overline{Q}}
\newcommand\Dg{{\mathbb D}_g}
\newcommand\Hg{{\mathbb H}_g}

\newcommand\PW{ {{\partial}\over{\partial W}} }
\newcommand\PWB{ {{\partial}\over {\partial{\overline W}}} }
\newcommand\OVW{\overline W}

\vskip 2mm
It seems that in his paper \cite{Si}, Siegel named ``{\sf Symplectic Geometry}"
for the geometry of the symplectic and hyperbolic manifold $\BH_g$ because
the group of all biholomorphic mappings (or isometries) of $\BH_g$ is the symplectic group
$Sp(2g,\BR)/\{ \pm I_{2g}\}.$ We note that Siegel did not endow $\BH_n$ with
the symplectic structure or the symplectic form.
\vskip 1mm
During the winter semester of 1951--1952
at the University of G{\"o}ttingen, Siegel delivered a series of lectures on celestial
mechanics. At that time J{\"u}rgen Moser~(1928--1999), a doctoral student attended the
Siegel's lectures. In 1955, Moser prepared a careful set of Siegel's lecture notes on
celestial mechanics for the publication. Finally in 1956, Siegel's lecture notes were
published in German in the form of book under the title,
``{\sf Vorlesungen {\"u}ber Himmelsmechanik}" (cf.~\cite{Si1}).
Later this book became the basis
of the joint book, ``{\sf Lectures on celestial mechanics}" which was published in
1971 with an English translation~(cf.~\cite{Si2}). It is known that the joint book
of Siegel and Moser became a standard reference on celestial mechanics.
Siegel wrote twelve papers on celestial mechanics in addition to his above book.
In this book he dealt with the three body problem, periodic solutions to the Euler-Lagrangian
equation and the Hamilton's equation, and the stability problem. He proved that the Lagrangian
system can be transformed into the Hamiltonian system through the Lagrangian transformation.
He characterized the canonical transformation in classical mechanics via the symplectic group.
The relation between the canonical transformation and the symplectic transformation was
clarified. We would like to mention the work of Aurel Wintner~(1903--1958). In the early 1940s
Wintner noticed that the set of completely canonical transformations which map a phase space
onto another phase space forms the symplectic group and developed the theory of canonical
transformations in celestial mechanics~(cf.~\cite{WA}).

\vskip 2mm
Now we briefly outline the history of symplectic geometry. Classically symplectic geometry
is the most {\it canonical\ mathematical\ structure} that appears whenever we describe a physical
system whose configurations can be parametrized by the variables $(p_k)$ and each change of
such a configuration requires the amount of work described by an integral of a one-form
$W=f_k dx_k$. Then the components $\{ f_k\}$ of the work $W$ are called\,: force, momentum,
pressure etc., whereas the space of states $P=\{ (p_k,q_k)\}$ is equipped with a canonical
symplectic structure $\omega=dW=dq_k\wedge dp_k.$ When the work does not depend on the
path but only on its ends, the submanifold of {\it physically\ admissible\ states}
constitutes a Lagrangian submanifold of $P$. However, before thermodynamics has emerged
as an important branch of physics, using the {\it canonical\ structure} of the phase space,
many books concerning ``{\sf Analytical\ Mechanics}" have been published following the
works of Joseph L. Langrange\,{\small (1736--1813)},
Sim{\' e}on D. Poisson\,{\small (1781--1840)}, Sir William R. Hamilton\,{\small (1805--1865)},
Karl G.\,J. Jacobi\,{\small (1804--1851)}, Konstantin Caratheodory\,{\small (1873--1950)}
and others. The word {\it canonical\ structure} means precisely the {\it symplectic geometry}.
\vskip 2mm
It is known that Charles Ehresmann\,{\small (1905--1979)} defined the notion of
``symplectic manifolds" for the first time in 1950 developing the theory of
fibre bundles. In the early 1950s Jean Marie
Souriau\,{\small (1922--2012)} introduced the notions of a symplectic
vector space and Lagrangian submanifolds. He also gave some applications
to classical mechanics, statistical mechanics and quantum mechanics.
For instance, the structure of the Hamiltonian equations
can be described via the symplectic structure. The matrices which map a
Hamiltonian system into a Hamiltonian system are symplectic. In the Hamiltonian
mechanics, if the configuration space is a manifold of dimension $n$, then the
momentum phase space is a symplectic manifold of dimension $2n$.
Siegel used the symplectic group to describe the canonical transformation of
a Hamiltonian system into a Hamiltonian system\,(cf.~\cite{Si1,Si2}). In fact, symplectic geometry
is realized as the geometrical interpretation of classical mechanics. We may say
that Siegel and Souriau introduced the concept of symplectic structures into the
Hamiltonian mechanics for the first time.

\vskip 2mm
In the early 1950s Souriau initiated the ``Geometric Quantization" as a new branch of
mathematical physics. In our times, the most important impulse for the development of
research on symplectic geometry was given by the ground-breaking book,
``{\bf Structure\ des\ syst$\grave{\rm e}$mes\ dynamiques}"~\cite{So,So1} which was
published in 1970. In this book classical mechanics, statistical mechanics and quantum
mechanics are treated by means of the very unusual innovative ideas.
The papers by Guillemin and Sternberg\,\cite{Gu-St1, Gu-St2},
Marsden\,\cite{A-M,M-W}, Weinstein\,\cite{W1,W-1,W-2,W2,W3} and Kirillov\,\cite{K1, K2}
were strongly inspired by Souriau.
In the 1970s Arnold, Marsden, Weinstein et al developed the theory of symplectic geometry
linking to that of classical mechanics and quantum mechanics.
In the mid-1980s, McDuff school
was motivated by Gromov's work on the pseudoholomorphic curve\,(cf.\,\cite{Gr1, Gr2}) and
then began to develop the theory of symplectic topology\,
(cf.\,\cite{LMS, Mc2, Mc3, Mc4, Mc4+, Mc5, Mc6, Mc7}).
The terminology ``symplectic geometry" was listed and classified as a mathematical subject
(53D22 and 53D25) in MSC2010 and MSC2020 databases of the American Mathematical Society.
The first issue of ``The Journal of Symplectic Geometry" was published in 2001. This journal
publishes papers related to symplectic geometry bimonthly.

\end{section}

%%%%%%%%%%%%%%%%%%%%%%%%%%%%%%%%%%%%%%%%%%%%%%%%%%%%%%%%%%%%%%%%%%%%%%%%%%%%%%%%%%%%%%%%%%%%%%%%%%%%%%%%%%%%
%%%%%%%%%%%%%%%%%%%%%%%%%%%%%%%%%%%%%%%%%%%%%%%%%%%%%%%%%%%%%%%%%%%%%%%%%%%%%%%%%%%%%%%%%%%%%%%%%%%%%%%%%%%%
%%%%%%%%%%%%%%%%%%%%%%%%%%%%%%%%%%%%%%%%%%%%%%%%%%%%%%%%%%%%%%%%%%%%%%%%%%%%%%%%%%%%%%%%%%%%%%%%%%%%%%%%%%%%
%%%%%%%%%%%%%%%%%%%%%%%%%%%%%%%%%%%%%%%%%%%%%%%%%%%%%%%%%%%%%%%%%%%%%%%%%%%%%%%%%%%%%%%%%%%%%%%%%%%%%%%%%%%%
%%%%%%%%%%%%%%%%%%%%%%%%%%%%%%%%%%%%%%%%%%%%%%%%%%%%%%%%%%%%%%%%%%%%%%%%%%%%%%%%%%%%%%%%%%%%%%%%%%%%%%%%%%%%
%%%%%%%%%%%%%%              3. Basic notions, definitions and examples                  %%%%%%%%%%%%%%%%%%%
%%%%%%%%%%%%%%%%%%%%%%%%%%%%%%%%%%%%%%%%%%%%%%%%%%%%%%%%%%%%%%%%%%%%%%%%%%%%%%%%%%%%%%%%%%%%%%%%%%%%%%%%%%%%
%%%%%%%%%%%%%%%%%%%%%%%%%%%%%%%%%%%%%%%%%%%%%%%%%%%%%%%%%%%%%%%%%%%%%%%%%%%%%%%%%%%%%%%%%%%%%%%%%%%%%%%%%%%%
%%%%%%%%%%%%%%%%%%%%%%%%%%%%%%%%%%%%%%%%%%%%%%%%%%%%%%%%%%%%%%%%%%%%%%%%%%%%%%%%%%%%%%%%%%%%%%%%%%%%%%%%%%%%
%%%%%%%%%%%%%%%%%%%%%%%%%%%%%%%%%%%%%%%%%%%%%%%%%%%%%%%%%%%%%%%%%%%%%%%%%%%%%%%%%%%%%%%%%%%%%%%%%%%%%%%%%%%%
%%%%%%%%%%%%%%%%%%%%%%%%%%%%%%%%%%%%%%%%%%%%%%%%%%%%%%%%%%%%%%%%%%%%%%%%%%%%%%%%%%%%%%%%%%%%%%%%%%%%%%%%%%%%
%%%%%%%%%%%%%%%%%%%%%%%%%%%%%%%%%%%%%%%%%%%%%%%%%%%%%%%%%%%%%%%%%%%%%%%%%%%%%%%%%%%%%%%%%%%%%%%%%%%%%%%%%%%%

\vskip 12mm
\begin{section}{{\bf Basic notions, definitions and examples}}
\setcounter{equation}{0}
\vskip 3mm
Let $(M,\omega)$ be a symplectic manifold of dimension $2n$, that is, a smooth manifold of dimension $2n$
equipped with a closed $(d\omega=0)$, nondegenerate ($\omega^n\neq 0$) $2-$form $\omega$. The notion of symplectic structures arose in the Hamiltonian formulation of the theory of classical mechanics.
A classical mechanical system can be modelled by the phase space which is a symplectic space. On the other hand, a quantum mechanical system is modelled by a Hilbert space. Each state of the system corresponds to a line in a Hilbert space.

\vskip 2mm
\begin{definition}\label{def:3.1}
A diffeomorphism $\phi:M_1\lrt M_2$ of two symplectic manifolds $(M_1,\omega_1)$ and $(M_2,\omega_2)$
is called a {\sf symplectomorphism} if $\phi^* \omega_2=\omega_1.$ In this case, we say that
$(M_1,\omega_1)$ is symplectomorphic to $(M_2,\omega_2)$. We denote by ${\rm Symp}(M,\omega)$
the group of all symplectomorphisms $\phi:(M,\omega)\lrt (M,\omega).$
\end{definition}

\vskip 2mm
\begin{definition}\label{def:3.2}
A real vector space $(V,\omega)$ is said to be a symplectic vector space if $V$ is equipped
with a non-degenerate alternating bilinear form $\omega:V\times V \lrt \BR.$
A subspace $W$ of a symplectic vector space $(V,\omega)$ is called {\sf Lagrangian} if
$\omega|_{W\times W}\equiv 0$ and $\dim W=\frac{1}{2}\,\dim V.$
A submanifold $N$ of a symplectic manifold $M$ is called {\sf Lagrangian } if,
for each point $p\in N$, $T_p N$ is a Lagrangian subspace of $T_p M$.
\end{definition}

\vskip 2mm
\begin{definition}\label{def:3.3}
Let $(M,\omega)$ be a symplectic manifold. An almost complex structure $J=\{ J_p\}$ on $M$ is said to be
{\sf compatible \ with} $\omega$ (or $\omega$-compatible) if
$$g_J:=\{ g_p:T_p M\times T_p M \lrt \BR,\ p\in M\}$$
defined by
\begin{equation*}
g_p (X,Y):=\omega_p (X,J_p Y),\qquad X,Y\in T_p M, \ p\in M
\end{equation*}
is a Riemannian metric on $M$.
\end{definition}

The following facts (CJ1)--(CJ3) are well known.
\vskip 2mm\noindent
(CJ1) Any symplectic manifold $(M,\omega)$ has a $\omega$-compatible almost complex structure.
\vskip 2mm\noindent
(CJ2) The set of all $\omega$-compatible almost complex structures on a symplectic manifold\\
\indent \ \ \ \ \ $(M,\omega)$
is path-connected and contractible.
\vskip 2mm\noindent
(CJ3) Let $(M,\omega)$ be a symplectic manifold equipped with a $\omega$-compatible almost \\
\indent \ \ \ \ \ complex structure $J$. Then any almost complex submanifold $N$ of $(M,J)$ is a \\
\indent \ \ \ \ \ symplectic submanifold of $(M,\omega)$.

\vskip 5mm
Let $(M,h)$ be a K{\"a}hler manifold of dimension $n$ with hermitian metric $h=(h_{ij})$. Then its K{\"a}hler form
$\omega$ is given by
\begin{equation}\label{(3.1)}
\omega=\frac{i}{2}\,\sum_{i,j=1}^n h_{ij} dz_i\wedge d{\overline z}_j,\qquad [\omega]\in H^{1,1}(M,\BC)\cap H^2(M,\BR).
\end{equation}
Then $\omega$ is a symplectic form on $M$. Thus $(M,\omega)$ is a symplectic manifold. According to the positivity of $h$, we see that the symplectic form $\omega$ satisfies the positivity condition
\begin{equation}\label{(3.2)}
\omega (X,JX) > 0\quad {\rm for\ all}\ X \in {\mathfrak X}(M),
\end{equation}
where ${\mathfrak X}(M)$ denotes the vector space of all smooth vector fields on $M$.
Moreover the complex structure $J$ satisfies the property $J^* \omega=\omega$, where
$J^*\omega (X,Y):=\omega (JX,JY)$ for all $X,Y \in {\mathfrak X}(M).$
The symplectic volume form is
\begin{equation}\label{(3.3)}
\frac{\omega^n}{n!}=\left( \frac{i}{2}\right)^n\, \det (h_{ij})\,dz_1\wedge
d{\overline z}_1\wedge \cdots
\wedge dz_n\wedge d{\overline z}_n.
\end{equation}
We may say that a K{\"a}hler manifold is a symplectic manifold $(M,\omega)$ equipped with an integral
$\omega$-compatible almost complex structure $J$.

\vskip 2mm
It is of interest to understand when the diffeomorphism type of a smooth family of
symplectic forms remains {\sf stable}. We can ask the following natural question.
When are two volume forms on the same manifold diffeomorphic\,?
In 1965 Moser proved the following result.
\begin{theorem}\label{thm:3.4}
  If $\omega$ and $\tau$ are two volume forms on a smooth connected oriented smooth manifold
$M$ with
$$ \int_M\omega = \int_M\tau,$$
then there exists a diffeomorphism $\phi:M\lrt M$ such that $\tau=\phi^*\omega.$
\end{theorem}
\begin{proof}
  The proof can be found in \cite{Mo}.
\end{proof}
\vskip 2mm
Moser also proved the following {\sf stability theorem:}
\begin{theorem}\label{thm:3.5}
  If $\{ \omega_t\,\vert \ t\in [0,1]\}$ is a smooth family of cohomologous symplectic forms on
  a {\sf compact} connected smooth manifold $M$, then there exists a smooth family
  $\{ \varphi_t\,\vert \ t\in [0,1]\}$ of diffeomorphisms of $M$ such that
  $\varphi_t^*\omega_t=\omega_0$, and $\varphi_0$ is the identity map on $M$.
\end{theorem}
\begin{proof}
  The proof can be found in \cite{Mo}.
\end{proof}
\begin{corollary}\label{coro:3.6}
Let $M$ be a {\sf compact} complex manifold. Let $\omega_1$ and $\omega_2$
be K{\"a}hler forms on $M$. Assume that $[\omega_1]=[\omega_2]\in H^{(1,1)}(M,\BC)\cap H^2(M,\BR).$
Then $(M,\omega_1)$ is symplectomorphic to $(M,\omega_2)$.
\end{corollary}
\begin{proof}
  The proof follows immediately from Theorem 3.4.
\end{proof}

\begin{theorem}\label{thm:3.7}
Let $\omega$ be a closed real $(1,1)$-form on a complex manifold $M$ and let $p\in M.$ Then
there exist a neighborhood $U$ of $p$ and a K{\"a}hler potential $\varphi\in C^\infty (U;\BR)$
such that
$$
\omega=\frac{i}{2}\,\partial \overline{\partial}\varphi \qquad {\rm on}\ U.
$$
\end{theorem}

\begin{definition}\label{def:3.8}
Let $G$ be a connected Lie group and let $(M,\omega)$ be a symplectic manifold.
A smooth $G$-action $\phi:G\times M\lrt M$ is said to be {\sf symplectic} if $G$ acts by symplectomorphisms, that is,  for each $g\in G$ the diffeomorphism $\phi_g:M\lrt M$ given by $\phi_g (p):=\phi(g,p),\ p\in M$ satisfies the condition $\phi_g^*\omega=\omega.$ The triple $(M,\omega,\phi)$
is called a {\sf symplectic\ $G$-manifold}.
\end{definition}

\begin{definition}\label{def:3.9}
Let $G$ be a connected Lie group and let $(M,\omega)$ be a symplectic manifold. Let $\frak g$ be the Lie algebra of $G$. Let $\phi:G\times M\lrt M$ be an action of $(M,\omega)$. For any $X\in \frak g$,
we denote by $X^M$ the vector field on $M$ generated by the one-parameter subgroup of global diffeomorphisms
$p\mapsto \phi(\exp (tX),p),\ p\in M$, that is,
\begin{equation}\label{(3.4)}
  X^M (p)=\frac{d\ }{dt}\Big|_{t=0}\phi(\exp (tX),p),\quad p\in M.
\end{equation}
\end{definition}
Given a smooth function $H:(M,\omega)\lrt \BR$, let $X_H$ be the vector field on $M$ defined by the Hamilton's equation
\begin{equation}\label{(3.5)}
\imath_{X_H}\omega=\omega(X_H,\,\cdot\,)=-dH.
\end{equation}

\begin{definition}\label{def:3.10}
A smooth vector field $Y$ on a symplectic manifold $(M,\omega)$ is said to be {\sf symplectic} if its flow preserves the symplectic structure $\omega$, and {\sf Hamiltonian} if there exist a smooth function
$H:M\lrt \BR$ such that $Y=X_H.$
\end{definition}

\begin{definition}\label{def:3.11}
Let $\sigma:G\lrt {\rm Symp}(M,\omega)$ be a symplectic action of a Lie group $G$ on
a symplectic manifold $(M,\omega)$. Let $\mathfrak{g}$ be the Lie algebra of $G$.
The action $\sigma$ is called a {\sf Hamiltonian action} if
there exists a map $\mu:M\lrt \frak g^*$ satisfying the following conditions {\rm (HA1)} and {\rm (HA2)}:
\vskip 2mm\noindent
{\rm (HA1)} For each $X\in \frak g$,
\begin{equation}\label{(3.6)}
d\mu_X=\imath_{X^M}\omega,
\end{equation}
\indent \ \ \ \ \ where $\mu_X:M\lrt \BR$ is a function defined by $\mu_X (p):=\langle \mu (p),X\rangle$ and $X^M$ is \\
\indent \ \ \ \ \ the vector field on $M$ generated by the one-parameter subgroup
$\{ \exp tX\,|\ t\in\BR \}$ \\
\indent \ \ \ \ \ of $G$.
\vskip 2mm\noindent
{\rm (HA2)} $\mu$ is equivariant with respect to the given action $\sigma$ and the coadjoint action \\
\indent \ \ \ \ \ ${\rm Ad}^*:G\lrt GL(\frak g^*)$,
that is,
\begin{equation}\label{(3.7)}
\mu\circ \sigma_g ={\rm Ad}^*(g)\circ \mu \qquad {\rm for\ all}\ g\in G.
\end{equation}
In this case, the quadruple $(M,\omega,G,\mu)$ is called a {\sf Hamiltonian\ $G$-space} and $\mu$ is called the {\sf moment\ map}.
\end{definition}

\begin{remark}\label{rk:3.12}
(1) A $G$-action on a symplectic manifold $(M,\omega)$ is symplectic if and only if all the vector fields
$X^M\ (X\in \frak g)$ are symplectic if and only if all the one forms $\imath_{X^M}\omega$ are closed.
A $G$-action on a symplectic manifold $(M,\omega)$ is Hamiltonian if and only if all the vector fields
$X^M\ (X\in \frak g)$ are Hamiltonian if and only if all the one forms $\imath_{X^M}\omega$ are exact.
\par\noindent
%\vskip 2mm\noindent
(2) Any symplectic $G$-action on a simply connected symplectic manifold $(M,\omega)$ is Hamiltonian.
Indeed, $[\imath_{X^M}\omega]\in H^1(M,\BR)=0$ for all $X^M\ (X\in \frak g)$ and hence all
$\imath_{X^M}\omega$ are exact.
\end{remark}

\begin{remark}\label{rk:3.13}
Let $T$ be a torus. The fixed point set in a Hamiltonian $T$-space $(M,\omega,T,\mu)$ is a finite union of connected symplectic submanifolds of $M$.
\end{remark}

Meyer\,\cite{Me} and Marsden and Weinstein\,\cite{M-W} proved the following.
\begin{theorem}\,\!\!\!{\sf (Meyer\,[1973] and Marsden-Weinstein\,[1974])}\label{thm:3.15}
Let $(M,\omega,G,\mu)$ be a Hamiltonian $G$-space for a {\bf compact} Lie group $G$. Let $i:\mu^{-1} (0)
\hookrightarrow M$ be the inclusion map. Assume that $G$ acts freely on $\mu^{-1} (0)$. Then
\vskip 2mm\noindent
(a) the orbit space $M_{\rm red}:=\mu^{-1} (0)/G$ is a manifold,
\vskip 2mm\noindent
(b) $\pi:\mu^{-1} (0)\lrt M_{\rm red}$ is a principal $G$-bundle, and
\vskip 2mm\noindent
(c) there exists a symplectic form $\omega_{\rm red}$ on $M_{\rm red}$ such that
$i^*\omega=\pi^*\omega_{\rm red}$.
\end{theorem}
We refer to \cite{M-W, Me} for more details.

\begin{definition}\label{def:3.16}
Let $(M,\omega)$ be a symplectic manifold of dimension $2n$. The {\sf symplectic measure} (or
{\sf Liouville\ measure}) of a Borel subset $U$ of $M$ is defined to be
$$m_\omega (U):=\int_U \frac{\omega^n}{n!},$$
where $\omega^n /{n!}$ is the symplectic volume form of $M$.
\end{definition}

\begin{definition}\label{def:3.17}
Let $G$ be a torus of dimension $n$. Let $(M,\omega,G,\mu)$ be a Hamiltonian $G$-space of dimension $2n$ such that the moment map $\mu$ is proper. The {\sf Duistermaat-Heckman measure} (briefly {\sf D-H measure}), $m_{\rm DH}$, on $\frak g^*$ is defined to be the push-forward of $m_\omega$ by $\mu:M\lrt \frak g^*$. More precisely, for any Borel subset $W$ of
$\frak g^*$, we have
$$m_{\rm DH}(W)=(\mu_*m_\omega) (W):=\int_{\mu^{-1}(W)} \frac{\omega^n}{n!}.$$
where ${\omega^n}/{n!}$ is the symplectic volume form of $M$.
\end{definition}

\vskip 3mm
For a function $f\in C_c^{\infty}(\frak g^*),$ we define its integral with respect to the D-H measure to be
\begin{equation*}
  \int_{\frak g^*} f\, dm_{\rm DH}=\int_M (f\circ\mu)\,\frac{\omega^n}{n!}.
\end{equation*}
\vskip 2mm
On $\frak g^*\simeq\BR^n$, there is also the Lebesque measure $m_0$. The relation between $m_{\rm DH}$ and
$m_0$ is governed by the Radon-Nikodym derivative, denote by $\frac{dm_{\rm DH}}{dm_0}$, which is a generalized function satisfying
\begin{equation*}
  \int_{\frak g^*} f\, dm_{\rm DH}=\int_{\frak g^*} f\,\frac{dm_{\rm DH}}{dm_0}\,dm_0.
\end{equation*}

Duistermaat and Heckman\,\cite{D-H} proved the following.
\begin{theorem}\,\!\!\!{\sf (Duistermaat-Heckman\,[1982])}.\label{thm:3.18}
Let $G$ be a torus of dimension $n$. Let $(M,\omega,G,\mu)$ be a a Hamiltonian $G$-space of dimension $2n$ such that the moment map $\mu$ is proper.
The D-H measure $m_{\rm DH}$ on $\frak g^*$ is a piecewise
polynomial multiple of the Lebesque measure $m_0$ on $\frak g^*\simeq\BR^n$, that is, the Radon-Nikodym derivative
$$\vartheta=\frac{dm_{\rm DH}}{dm_0}$$
is a piecewise polynomial. More precisely, for any Borel subset $W$ of $\frak g^*$,
$$ m_{\rm DH}(W)=\int_W \vartheta (x)\,dx,$$
where $dx=dm_0$ is the Lebesque volume form on $W$ and $\vartheta:\frak g^*\lrt \BR$ is a polynomial on any region consisting of regular values of $\mu$. In particular, for the standard Hamiltonian action of $S^1$ on $(S^2,\omega)$, we have $m_{\rm DH}=2\pi\,m_0$ with a constant polynomial $\vartheta=2\pi.$
\end{theorem}
\begin{proof}
The proof can be found in \cite{D-H}.
\end{proof}

\vskip 2mm
Alan Weinstein\,(1943-- ) proved the following remarkable results.
\vskip 1mm\noindent
{\bf Weinstein Lagrangian Neighborhood Theorem:} {\it Let $M$ be a smooth manifold of dimension $2n$, $X$
a compact $n$-dimensional submanifold, $i:X\hookrightarrow M$ the inclusion map, and $\omega_1$ and $\omega_2$ symplectic forms on $M$ such that $i^*\omega_1=i^*\omega_2=0$, i.e., $X$ is a Lagrangian
submanifold of both $(M,\omega_1)$ and $(M,\omega_2)$. Then there exist neighborhoods $U_1$ and $U_2$
of $X$, and a diffeomorphism $\phi:U_1\lrt U_2$ such that $i_2=\phi\circ i_1$ and $\phi^*\omega_2=\omega_1$, where $i_k:X\lrt U_k\ (k=1,2)$ are the inclusion maps.}
\begin{proof} The proof can be found in \cite{W1, W3}.\end{proof}
\vskip 1mm
We recall that a {\rm Lagrangian\ submanifold} of a symplectic manifold is a submanifold of
the greatest possible dimension on which the symplectic structure restricts to zero.
According to the above theorem, a tubular neighborhood of any Lagrangian submanifold $L$ in any
symplectic manifold is isomorphic (as a symplectic manifold) to a tubular neighborhood of the
zero section of the cotangent bundle $T^*L$ of $L$. Weinstein mentions that
``{\it everything\ in\ the\ world\ is a\ Lagrangian\ submanifold}\,".
 By that, he means that everything in symplectic geometry is most naturally expressed in terms of Lagrangian submanifolds.

\vskip 2mm\noindent
{\bf Weinstein Tubular Neighborhood Theorem:} {\it Let $M$ be a smooth manifold of dimension $n$, $X$ is a submanifold of $M$, $NX$ the normal bundle of $X$ in $M$, $i_0:X \hookrightarrow NX$ the zero section, and
$i:X\hookrightarrow M$ inclusion. Then there exist neighborhoods $U_0$ of $X$ in $NX$, $U_2$ of $X$ in $M$
and a diffeomorphism $\phi:U_0\lrt U_1$ such that $i=\phi\circ i_0$.}
\begin{proof} The proof can be found in \cite{W2, W3}.\end{proof}

\vskip 3mm\noindent
{\bf Examples:} {\bf (1)} Let $x_1,\cdots,x_n,x_{n+1},\cdots, x_{2n}$ be linear coordinates on $\BR^{2n}$.
Then
\begin{equation*}
  \omega_0:=\sum_{i=1}^n dx_i \wedge dx_{n+i}
\end{equation*}
is a symplectic form. Thus $(\BR^{2n},\omega_0)$ is a symplectic manifold. Now we consider the case $n=1$.
The one-dimensional torus
$$S^1=SO(2)=\left\{ \begin{pmatrix}
              \cos t & -\sin t \\
              \sin t & \ \cos t
            \end{pmatrix}\in \BR^{(2,2)}\ \Big|\ 0\leq t < 2\pi  \right\}$$
acts on $\BR^2$ by rotations. This is a Hamiltonian action. The Lie algebra of $S^1$ is given by
\begin{equation*}
  \frak s {\frak o}(2)=\left\{ \begin{pmatrix}
              0 & -a \\
              a & \ 0
            \end{pmatrix}\in \BR^{(2,2)}\ \Big|\ a\in \BR  \right\}\cong \BR.
\end{equation*}
$J=\begin{pmatrix}
              0 & -1 \\
              1 & \ 0
            \end{pmatrix}$ is a basis of $\frak s {\frak o}(2)$. Since
\begin{equation*}
  \exp (tJ)=\begin{pmatrix}
              \cos t & -\sin t \\
              \sin t & \ \cos t
            \end{pmatrix},\quad 0\leq t < 2\pi,
\end{equation*}
the vector field $J^{\sharp}$ on $\BR^2$ generated by $J$ is given by
\begin{equation*}
  J^{\sharp}=y\,\frac{\partial\ }{\partial x}-x\,\frac{\partial\ }{\partial y}.
\end{equation*}
The moment map $\mu:\BR^2\lrt \frak s {\frak o}(2)^*$ satisfies
\begin{equation*}
 \imath_{J^{\sharp}}\omega_0=x\,dx +y\,dy=d\langle \mu(\,\cdot\,),J\rangle.
\end{equation*}
Thus we obtain
\begin{equation*}
 \langle \mu(x,y),J\rangle=\frac{1}{2}\,(x^2+y^2),\quad (x,y)\in \BR^2.
\end{equation*}

\vskip 2mm\noindent
{\bf (2)} The standard symplectic form $\omega$ on the two dimensional sphere $S^2$ is induced by
\begin{equation*}
\omega_p (X,Y):=\langle p,X\times Y \rangle, \quad p\in S^2,\ X,Y\in T_p S^2 =\{ p\}^\perp.
\end{equation*}
Thus $(S^2,\omega)$ is a symplectic manifold. Let $\omega_{\rm st}=d\theta\wedge dh$ be the standard area form in cylindrical polar coordinates $\theta,h\ (0\leq \theta \le 2\pi$ and $-1\leq h\leq 1)$ on $S^2$. Endow $(S^2,\omega_{\rm st})$ with the rotational $S^1$-action about the $z$-axis. This is a Hamiltonian action with the moment map $\mu:S^2\lrt \BR$ given by $\mu (\theta,h)=h,$ and the momentum polytope is $\mu (S^2)=[-1,1].$ The Duistermaat-Heckman polynomial is $\vartheta=2\,\pi\,\chi_{[-1,1]}$, where $\chi_{[-1,1]}$ is the characteristic function of $[-1,1]$. Hence $m_{\rm DH}([a,b])=2\pi (b-a)$ for all $[a,b]\subset [-1,1].$
Assume that the $n$-dimensional sphere $S^n\ (n>2)$ admits a symplectic form $\omega$.
Then $d\omega=0$ because $\omega$ is closed. By Stokes' theorem, the cohomology class
$[\omega^k]$ is nontrivial in $H^{2k}(S^n,\BR)$ for $1\leq k\leq \frac{n}{2}$.
Since
\[ H^q(S^n,\BR)=
 \begin{cases}
  \BR, & \mbox{if } q=0,\ n \\
  0, & \mbox{otherwise},
\end{cases}\]
we obtain a contradiction. Therefore $S^2$ is the only sphere $(n\geq 1)$ which admits
a symplectic structure or a symplectic form.

\vskip 2mm\noindent
{\bf (3)} Let $M$ be a smooth manifold of dimension $n$. Let $x_1,\cdots,x_n$ be local coordinates on
an open neighborhood $U$ of $x\in M$. If $\alpha\in T_x^*X$, then
\begin{equation*}
  \alpha=\sum_{i=1}^n \alpha_i (dx_i)_x, \quad \alpha_i \in \BR.
\end{equation*}
Then $(T^*U, x_1,\cdots,x_n,\alpha_1,\cdots,\alpha_n)$ is a local coordinate chart for $T^*M.$
We define
\begin{equation*}
\Omega_{\rm can}:=\sum_{i=1}^{n} dx_i \wedge  d\alpha_i \qquad {\rm and}\qquad
\alpha: = \sum_{i=1}^{n} \alpha_i\, dx_i.
\end{equation*}
Then it is easily seen that $\alpha$ is intrinsically defined and $\Omega_{\rm can}=-d\alpha$.
Thus $(T^*M,\Omega_{\rm can})$ is a symplectic manifold of dimension $2n$. The zero section of
$T^*M$ given by
\begin{equation*}
  (T^*M)_0:=\{ (x,\alpha)\in T^*M\,|\ \alpha=0 \ {\rm in}\ T_x^*M  \}
\end{equation*}
is a Lagrangian submanifold of $(T^*M,\Omega_{\rm can})$. The 1-form $\alpha$ is called the
{\sf tautological\ form} and $\Omega_{\rm can}$ is called the {\sf canonical} symplectic form
on $T^*M$.
\vskip 2mm
Suppose the Lie group $G$ acts on a smooth manifold $M$. Then $G$ acts naturally on $TM$ and $T^*M$.
We can show that the action of $G$ on $(T^*M,\Omega_{\rm can})$ is Hamiltonian. The moment map
$\mu:T^*M\lrt \frak g^*$ is given by
\begin{equation*}
\langle \mu (\,\cdot\,),Y\rangle =-\imath_{Y^{\sharp}}\alpha\qquad {\rm for\ all}\ Y\in \frak g,
\end{equation*}
where $Y^{\sharp}$ is the vector field on $T^*M$ generated by $Y$.

\vskip 2mm\noindent
{\bf (4)} For any submanifold $Y$ of a smooth manifold $M$, the conormal bundle $N^*Y$ is
a Lagrangian submanifold of $(T^*M,\Omega_{\rm can})$.
\vskip 2mm\noindent
{\bf (5)} Let $(M_1,\omega_1)$ and $(M_2,\omega_2)$ be two symplectic manifolds of dimension $2n$.
Let $\pi_i:M_1\times M_2 \lrt M_i\ (i=1,2)$ be the natural projections. For any nonzero real numbers $a$ and $b$, the 2-form $a\,\pi^*_1\omega_1 +b\,\pi^*_2\omega_2$ is a symplectic form on $M_1\times M_2$.
We consider the twisted product form
\begin{equation*}
  {\widetilde \omega}:=\pi^*_1\omega_1 - \pi^*_2\omega_2.
\end{equation*}
Let $\phi:M_1\lrt M_2$ be a diffeomorphism. Then we can show that $\phi$ is a symplectomorphism if and only if the graph of $\phi$ is a Lagrangian submanifold of $(M_1\times M_2,\widetilde \omega).$

\def\pa{\partial}
\def\A{\left[\begin{matrix} A\\ 0\end{matrix}\right]}
\def\Bb{\left[\begin{matrix} A\\ B\end{matrix} \right]}
\def\Dd{\left[\begin{matrix} -A \\ -B \end{matrix} \right]}
\def\x{\left[\begin{matrix} A+\xi\\ B+\eta\end{matrix} \right]}
\def\lam{\left[\begin{matrix} A+\lambda\\ B+\mu\end{matrix}\right]}
\def\ze{\left[\begin{matrix} 0\\ 0\end{matrix}\right]}
\def\J{J\in {\Bbb Z}^{(m,n)}_{\geq 0}}
\def\N{N\in {\Bbb Z}^{(m,n)}}
\def\Dm{\left[\begin{matrix} -A\\ -B\end{matrix}\right]}
\def\dt{{{d}\over {dt}}\bigg|_{t=0}}
\def\lt{\lim_{t\to 0}}
\def\zhg{\BZ^{(m,n)}}
\def\bhg{\BR^{(m,n)}}
\def\ex{\par\smallpagebreak\noindent}
\def\Box{$\square$}
\def\pis{\pi i \sigma}
\def\sd{\,\,{\vartriangleright}\kern -1.0ex{<}\,}
\def\sc{\bf}
\def\wt{\widetilde}
\vskip 2mm\noindent
{\bf (6)} Let $G$ be a Lie group. Then the coadjoint orbit $\mathcal O (F)$ of $F\in \frak g^*$ is a symplectic manifold equipped with the Kostant-Kirillov-Souriau symplectic form $\omega_F$. Let
$$\mathcal O (F)^-:=(\mathcal O (F), -\omega_F).$$
Then the natural product action of $G$ on $M\times \mathcal O (F)^-$ is Hamiltonian with the following
moment map $\mu_F:M\times \mathcal O (F)^- \lrt \frak g^*$ defined by
\begin{equation*}
\mu_F (p,\xi):=\mu (p)-\xi.\qquad p\in M,\ \xi\in \xi\in \mathcal O (F)^-.
\end{equation*}
If the above action of $G$ is free, we obtain a reduced space with respect to the coadjoint orbit
$\mathcal O (F)$.

\vskip 2mm\noindent
{\bf (7)} For a positive integer $n$, let $T_1:=(\BR/\BZ)^{2n-1}$ and $T_2:=(\BR/\BZ)^{n}$ be two tori.
It is easily seen that the $T_1$-action on $(M,\omega)=((\BR/\BZ)^{2n},\sum_{i=1}^{n}dx_i\wedge dy_i)$
with coordinates $(x_1,y_1,\cdots, x_n,y_n)$ in $(\BR/\BZ)^{2n}$ by translation on the first $2n-1$ components $(x_1,y_1,\cdots, x_{n-1},y_{n-1},x_n)$ is free and symplectic. According to Remark 3.12, this action is not Hamiltonian. The $T_1$-orbits are coisotropic submanifolds of $M$ diffeomorphic to $T_1$.
The $T_2$-action on $(M,\omega)$ by translations on $(x_1,\cdots,x_n)$ is free, symplectic and hence not
Hamiltonian. Its $T_2$-orbits are Lagrangian submanifolds of $M$ diffeomorphic to $T_2$. Here a submanifold
$N$ of a symplectic manifold $(M,\omega)$ is said to be {\it coisotropic} if $T_xN$ is a coisotropic subspace of $(T_x M,\omega_x)$ for all $x\in N$.

\vskip 2mm\noindent
{\bf (8)} Let $N:=((\BR/\BZ)^2 \times S^2, dx\wedge dy+d\theta\wedge dh)$ be a $4$-dimensional symplectic manifold. The $2$-torus $T:=(\BR/\BZ)^2$ acts freely by translations on the left factor of $N$. Let $\BZ_2:=\BZ/2\BZ$ act on $S^2$ by rotating each point horizontally by $\pi$ radians, and let $\BZ_2$ act on $(\BR/\BZ)^2$ by the antipodal action on the first circle $\BR/\BZ$. Then the diagonal action of $\BZ_2$ on $N$ is free. So the quotient space $M:=(\BR/\BZ)^2 \times_{\BZ_2} S^2$ is a smooth manifold equipped with the symplectic form $\omega$ and the $T$-action inherited from the ones on $N$. We see easily that the action of $T$ on $M$ is symplectic but not free, and the $T$-orbits are symplectic tori of dimension $2$.
The orbit space $M/T=S^2/\BZ_2$ is an orbifold with two singular points of order $2$, the south and north poles of $S^2$.

\vskip 2mm\noindent
{\bf (9)} The projective curve $M={\mathbb P}^1(\BC)=\BC \cup \{ \infty \}$ is a K{\"a}hler manifold. Its K{\"a}hler metric
$ds^2_{\rm FS}$ on the local chart $U_0=\{ [z_0,z_1]\in M\,|\ z_0\neq 0\,\}$ is given by
$$ds^2_{\rm FS}=\frac{dx^2 +dy^2}{(x^2+y^2+1)^2}=\frac{dz\,d{\overline z}}{(|z|^2+1)^2},$$
where $z_1/z_0=z=x+iy\,(x,y\in\BR)$ is the usual coordinate on $\BC$.
Its K{\"a}hler form $\omega_{\rm FS}$ is given by
$$\omega_{\rm FS}=\frac{dx\wedge dy}{(x^2+y^2+1)^2}=\frac{i}{2}\,\frac{dz\wedge d{\overline z}}{(|z|^2+1)^2}.$$
The function $K(z):=\log (|z|^2 +1)$ is the K{\"a}hler potential, that is,
$$\omega_{\rm FS}=\frac{i}{2}\,\partial_z \partial_{\overline z}\log (|z|^2 +1).$$
Thus $(M,\omega_{\rm FS})$ is a two-dimensional symplectic manifold. $\omega_{\rm FS}$ is called
the Fubini-Study symplectic form.
The total area of
$M=\BC\cup \{ \infty \}$ with respect to $\omega_{\rm FS}$ is given by
\begin{equation*}
\int_M \omega_{\rm FS}=\int_{\BR^2} \frac{dx\wedge dy}{(x^2+y^2+1)^2}=\pi.
\end{equation*}
Since $M$ is diffeomorphic to $S^2$ by stereographic projection, we obtain
\begin{equation*}
\omega_{\rm FS}=\frac{1}{4}\,\omega_{\rm st}. \quad ({\rm see\ Example\ (2)}).
\end{equation*}

\vskip 2mm\noindent
{\bf (10)} Let $\mathbb T^n:=(\BR/\BZ)^n$ be the $n$-torus. For $\lambda>0$,
let $({\mathbb P}^n(\BC),\lambda\,\omega_{\rm FS})$ be the $n$-dimensional symplectic complex projective space with the Fubini-Study form $\lambda\,\omega_{\rm FS}$. We consider the rotational
$\mathbb T^n$-action on ${\mathbb P}^n(\BC)$ induced from the rotational $\mathbb T^n$-action on the $(n+1)$-dimensional complex plane. This is a Hamiltonian action with its moment map
$\mu:{\mathbb P}^n(\BC)\lrt \BR^n$ given by
\begin{equation}
  \mu([z_0:z_1:\cdots:z_n]):= \left( \frac{\lambda\,|z_1|^2}{\sum_{i=0}^{n} |z_i|^2},\cdots\cdot\cdot,\frac{\lambda\,|z_n|^2}{\sum_{i=0}^{n} |z_i|^2}\right).
\end{equation}
If $e_1=(1,0,\cdots,0)\in \BR^n,\cdots,e_n=(0,\cdots,0,1)\in \BR^n$, the momentum polytope
\begin{equation}
  \mu ({\mathbb P}^n(\BC))= {\rm Convex\ Hull}\,\{ 0,\lambda\,e_1,\cdots,\lambda\,e_n \}.
\end{equation}

\end{section}

\vskip 12mm
\begin{section}{{\bf Arnold's conjecture and pseudoholomorphic curves}}
\setcounter{equation}{0}
\vskip 3.5mm
According to the Darboux's theorem, a symplectic form (or symplectic structure) $\omega$ can always be written in the following form
\begin{equation}\label{(4.1)}
  \omega =\sum_{k=1}^n dp_k \wedge dq_k
\end{equation}
in suitable canonical coordinates $p_1,\cdots,p_n,q_1,\cdots,q_n.$ However these canonical coordinates are not uniquely determined. This theorem says that locally all symplectic forms are the same, i.e., all
symplectic structures are indistinguishable. We recall the Moser's stability theorem\,\cite{Mo}
which says that if
$\omega_t\ (t\in [0,1])$ is a smooth path of symplectic forms such that the cohomology class
$[\omega_0]=[\omega_t]\in H^2(M,\BR)$ for all $t\in [0,1]$,   then all these symplectic forms are the same in the sense that one can make them coincide by moving the points of $M$ in a suitable way. In other words,
one cannot change the symplectic form in any way by deforming it, provided that the cohomology class
$[\omega]\in H^2(M,\BR)$ is unchanged.
\vskip 2mm
In the $2$-dimensional case, a symplectomorphism can be characterized as an area preserving diffeomorphism. More precisely, if $S$ is a region in $\BR^2$ that is diffeomorphic to a disc $D$ and has the same area as $D$, then there exists a symplectomorphism $\phi:S\lrt D$.

\vskip 1mm
In 1913, George David Birkhoff\,(1884--1944) proved the geometric theorem of Poincar{\'e}.
\vskip 1mm
\begin{theorem}\label{thm:4.1}
  Suppose $f:S\lrt S$ is an area-preserving diffeomorphism of the closed annulus $S=\BR/\BZ \times [-1,1]$
  which preserves the two components of the boundary, and twists them in opposite directions. Then $f$ has at least two distinct fixed points.
\end{theorem}
\begin{proof} The proof can be found in \cite{B}.\end{proof}
\vskip 2mm
The above classical fixed point theorem says that any area-preserving twist of a closed annulus
has at least two geometrically distinct fixed points. This fixed point result can be traced
to the work of Poincar{\'e} in celestial mechanics \cite{Po1}, where he showed that the study
of the dynamics of certain cases of the restricted Three Body Problem may be reduced to
investigating area-preserving maps. This work led Poincar{\'e} to the theorem, which he stated
in \cite{Po2} in 1912, and proved in several cases. The complete proof was given by Birkhoff
in 1913. Arnold realized that the higher-dimensional version of the result of Poincar{\'e}
and Birkhoff should concern {\it symplectic maps}, that is, maps preserving a symplectic
form, rather than volume-preserving maps.
The above paragraph was excerpted from \cite[pp.\,384--385]{P2}.

\vskip 3mm
Arnold generalized and formulated Poincar{\'e}'s geometric theorem in the following way.
\begin{conjecture}\label{conj:4.2}
A symplectomorphism of a compact symplectic manifold, homologous to the identity transformation,
has at least as many fixed points as a smooth function on the manifold has critical points.
\end{conjecture}

\vskip 2mm
We say that a symplectomorphism is {\it homologous\ to\ the identity} if it belongs to
the commutator subgroup of the group of symplectomorphisms homotopical to the identity.
With coordinates, such a diffeomorphism is given by $x\longmapsto x+f(x)$, where $x$ is
a point of the plane and $f$ is a periodic function. It symplectic if and only if
the Jacobian $\det\left( D(x+f(x))/Dx \right)$ is identically $1$ and it is homologous to
the identity if and only of the mean value of $f$ is 0. We refer to \cite{Sik} for more
details.
\vskip 2mm
An influential precursor in the study of global aspects of symplectic geometry is Arnold's conjecture which is a high-dimensional analogue of the classical fixed point theorem of H. Poincar{\'e} and G. Birkhoff.
Arnold's conjecture provided one of the most influential sources of inspiration in symplectic
topology~(cf.~\cite{Mc6}). Arnold's conjecture is the invention of symplectic topology.
Audin \cite{Au3} says that an official birthdate of symplectic topology is October 27th, 1965.
In his paper \cite{A3} Arnold started with a definition:
\vskip 1mm
{\footnotesize By symplectic topology, I mean the discipline having the same relation to ordinary topology as the theory of Hamiltonian dynamical systems has to the general theory of dynamical systems.}
\vskip 1mm
He doubted and asked whether there is such a thing as symplectic topology.
The existence of symplectic topology was proved by Mikhael Gromov\,(1943-- ).
\vskip 3mm
Arnold's conjecture was proved by Conley-Zehnder \cite{Co-Z} for the standard symplectic
$2\,n$--torus $T^{2n}$, Fortune \cite{Fo, Fo-W} for the projective space $\mathbb{P}^n\BC,$
Floer \cite{F1} for many K{\"a}hler manifolds of negative curvature. Fukaya--Ono \cite{F-O} and Liu--Tian \cite{L-T} proved Arnold's conjecture for a compact symplectic manifold
using Floer homology which is an infinite dimensional analogue of Morse theory.
Nowadays Arnold's conjecture is described in the form of the following theorem:
\begin{theorem}\label{thm:4.3}
  Let $(M,\omega)$ be a closed symplectic manifold of dimension $2n$ and $\phi:M\lrt M$ an exact
  (or Hamiltonian) symplectomorphism of $(M,\omega)$ with only non-degenerate fixed points.
  Then the number of fixed points of $\phi$ is at least the sum of Betti numbers of $M$, that is,
  $\sum_{k=0}^{2n} \dim\,\!H^k(M,\BQ).$
\end{theorem}
The above theorem implies that the fixed points of $\phi$ correspond to the intersection points
of the graph $\{(p,\phi(p))\in M\times M\,\vert\ p\in M\}$ of $\phi$ and the diagonal set
$\Delta_M:=\{ (p,p)\in M\times M\,\vert\ p\in M\}$ in the symplectic manifold
$(M\times M, \omega\oplus -\omega)$. Thus, to study the fixed points of a Hamiltonian
diffeomorphism, it suffices to study how two Lagrangian submanifolds intersect and behave with
each other.

\vskip 2mm
In 1985 using the pseudoholomorphic curve technique and constructing the Gromov radius,
Gromov proved the famous ${\sf Gromov's\ nonsqueezing\ theorem}$ which is also called the principle of the
${\sf symplectic\ camel}$.
\vskip 3mm\noindent
{\bf Gromov's Nonsqueezing Theorem:} {\it Let $(\BR^{2n},\omega_0)$ be the symplectic manifold equipped with the standard symplectic form
\begin{equation}\label{(4.2)}
  \omega_0 =\sum_{k=1}^n dp_k \wedge dq_k ,
\end{equation}
where $(p_1,\cdots,p_n,q_1,\cdots,q_n)$ is a coordinate in $\BR^{2n}$. Let
\begin{equation*}
B(R):= \left\{ (p_1,\cdots, p_n,q_1,\cdots,q_n)\in \BR^{2n}\,\Big|\ \sum_{i=1}^n (p_i^2+q_i^2)
< R^2\,\right\},\quad R> 0
\end{equation*}
be the ball of radius $R$ in $\BR^{2n}$ and let
\begin{equation*}
Z(r):= \{ (p_1,\cdots, p_n,q_1,\cdots,q_n)\in \BR^{2n}\,|\  p_1^2+q_1^2 < r^2\,\},\quad r> 0
\end{equation*}
be the cylinder of radius $r$, each of them equipped with the symplectic form $\omega_0$. If we can find a symplectic embedding $\varphi:B(R)\lrt Z(r)$, then $R\leq r$.}
%\begin{proof} The proof can be found in \cite{Gr1}.\end{proof}
\vskip 2mm
The detailed proof can be found in \cite{Gr1}. The rough idea of the proof is given as follows.
Here, a pseudoholomorphic curve (or J-holomorphic curve) is defined to be a smooth map
$\phi$ from a Riemann surface $(S,j)$ to an almost complex manifold $(M,J)$ which is
almost complex in the sense that
\begin{equation*}
  J\circ d\phi=d\phi \circ j,
\end{equation*}
where $j$ and $J$ are complex structures on $S$ and $M$ respectively.
For instance, we can spell out what a pseudoholomorphic curve amounts to in the case that
$S=\BC$ and $M=\BC^n$. It turns out that the $\BR$-linear map $J:\BC^n\lrt \BC^n$
with $J^2=-1$ can be parametrized by an open set of $n\times n$ complex matrices
$\mu=(\mu_{\alpha\beta}).$ Thus our almost structure is represented by a matrix valued function
$\mu(z_*)$ of $z_*=(z_\alpha)\in \BC^n$. A pseudoholomorphic curve is given by a solution of
the system of partial differential equations
\begin{equation*}
  \frac{\partial z_\alpha}{\partial \bar{\tau}}+\sum_{\beta=1}^{n}\mu_{\alpha\beta}
  \frac{\overline{\partial z_\beta}}{\partial \tau}=0,
\end{equation*}
which can be thought of as a deformation of the ordinary Cauchy--Riemann equations,
for the vector--valued function $z_*(\tau).$
We refer to a good survey article of Donaldson \cite{Do2} for some details on pseudoholomorphic curves.

\vskip 2mm
First we note that the two coordinates $(p_1,q_1)$ that span the disk $B^2(R)$ are symplectic,
i.e.,
$$
\omega_0 \left( \frac{\partial}{\partial p_1},\frac{\partial}{\partial q_1} \right)\neq 0.
$$
For each $\omega_0$--compatible almost structure $J$, the cylinder has a slicing by
pseudoholomorphic disks of area $\pi r^2$. If the ball is embedded in the cylinder, this slicing
will induce a slicing of the ball. If $J$ is suitably compatible with the embedding,
this slicing of the ball has to have some slices of $\omega_0$--area $\geq \pi R^2.$ Hence
we must have $R\leq r.$

\begin{remark}\label{rk:4.4}
Later on, different proofs of Gromov's nonsqueezing theorem were found, some of which are more
elementary than Gromov's. We refer to \cite[\S 3]{H-Z} and \cite{Su-T} for different proofs.
\end{remark}

\begin{definition}\label{def:4.5}
We define a {\sf symplectic ball} (resp.\,{\sf symplectic cylinder}) of radius $R$ to be the image
of the standard ball (resp.\,cylinder) of radius $R$ by a symplectic embedding.
We say that a local diffeomorphism $\varphi$ has the {\sf nonsqueezing property} if there
is no symplectic ball $B$ whose image $\varphi (B)$ is contained in a symplectic cylinder with
radius strictly less than that of $B$.
\end{definition}

\begin{theorem}{\bf [Eliasberg, Ekeland-Hofer]}.\label{thm:4.6}
  If $\varphi$ is a local diffeomorphism of $\BR^{2n}$ such that both $\varphi$ and
  its inverse $\varphi^{-1}$ have the nonsqueezing property, then $\varphi$ is
  either symplectic or antisymplectic, that is,
$$
\varphi^* (\omega_0)=\pm \omega_0.
$$
\end{theorem}

\begin{corollary}{\bf [Symplectic Rigidity]}\label{coro:4.7}
Let $(M,\omega)$ be a symplectic manifold. Let ${\rm Symp}(M,\omega)$
(resp.\,${\rm Diff}(M))$ be the group of all symplectomorphisms (resp.\,diffeomor-
\par\noindent
phisms) of $M$ onto $M$. Then ${\rm Symp}(M,\omega)$ is closed in ${\rm Diff}(M)$ in the
topology of uniform convergence on compact sets.
\end{corollary}

\vskip 3mm
It follows from the above facts that symplectic geometry is intrinsically topological in nature.
Using the existence of a pseudo-holomorphic curve on a compact 4-dimensional symplectic manifold
with $b_2^+>1$, Taubes \cite{Ta3} proved that the symplectic structure on the 2-dimensional projective space ${\mathbb P}^2\BC$ is unique up to scalar multiplication.
Here, $b_2^+$ is the number of $+1$ eigenvalues of the intersection form on
the rational second homology. (Thus, $b_2^+={\frac{1}{2}}
\,\left\{ {\rm rank}(H_2(M))+{\rm sign}(M) \right\}$, where  ${\rm sign}(M)$
denotes the signature of $M$.) (cf. \cite[p.\,222]{Ta2})

\vskip 2mm
Gromov also proved the following result.
\begin{theorem}\label{thm:4.8}
If the limit of a uniformly (${\mathscr C}^0$) converging sequence of symplectomorphisms is a diffeomorphism, then it is symplectic.
\end{theorem}

\vskip 2mm\noindent
{\bf Eliashberg's\ Principle:} {\it An obstruction to symplectic embeddings (beyond the volume condition)
can be described by a J-holomorphic curve.}
\vskip 2mm
We refer to \cite{E1, El} for more details on Eliashberg's Principle.

\vskip 2mm
We recall that if $(M,\omega)$ is a symplectic manifold, a submanifold
$N\subset M$ is called {\sf symplectic} if the restriction of $\omega$
to $N$ is non-degenerate. Using the pseudoholomorphic curve technique,
Donaldson obtained the following results:
\begin{theorem}{\bf [Donaldson \cite{Do1}:1996]}\label{thm:4.9}
  Let $(M,\omega)$ be a compact symplectic manifold of dimension $2n$, and
  suppose that the de Rham cohomology class $[\omega/2\pi]\in H^2(M,\BR)$
  lies in the integral lattice $H^2(M,\BZ)/Torsion.$ Let $\Theta\in H^2(M,\BZ)$
  be a lift of $[\omega/2\pi]$ to an integral class. Then for sufficiently
  large integers $k$ the Poincar{\'e} dual of $k\Theta$, in $H_{2n-2}(M,\BZ)$,
  can be realized by a symplectic submanifold $N\subset M.$
\end{theorem}

\begin{theorem}{\bf [Donaldson \cite{Do1}:1996]}\label{thm:4.10}
  If $(M,\omega)$ is any compact symplectic manifold, then following hold:
  \vskip 2mm\noindent
  (1) $M$ contains symplectic submanifolds of any even codimension.
  \vskip 2mm\noindent
  (2) If $J$ is a compatible almost complex structure on $M$, then there
  are almost complex
  \par\ \
  structures $J'$ arbitrarily close in $C^0$ to $J$
  such that $M$ contains $J'$--holomorphic
  \par\ \
  curves.
\end{theorem}
\vskip 2mm
The coexistence of flexibility and rigidity (called soft and hard in \cite{Gr2})
is a particular and particularly interesting feature of symplectic geometry.
Rigidity has many incarnations: the Arnold conjecture on the number of fixed points
of Hamiltonian diffeomorphisms, $C^0$-rigidity for Hamiltonian diffeomorphisms,
Hofer's metric, the rigidity of the Poisson bracket \cite{Pol} etc. Flexibility is manifest
in $h$--principles, of which new ones have been found in \cite{El}, and in
Donalsdon's theorem on symplectic hypersurfaces \cite{Do1}. Both rigidity and
flexibility are omnipresent in symplectic embedding problems (which will be discussed
in the subsection 6.2). The advantage here is that symplectic embedding problems
give rise to numbers, they can quantify symplectic flexibility and rigidity, and
localize the boundary between them. In Nonsqueezing Theorem there is only rigidity,
for the problem $E(1,a)\stackrel{s}{\hookrightarrow} C^4(A)$ (cf.\,Remark 6.19)
there is a subtle proximity of flexibility and rigidity, and for ball packings of
linear tori there is only flexibility.
\begin{remark}\label{rk:4.11}
The theory of Floer homology is based on pseudoholomophic curves with boundary
lying on a Lagrangian submanifold (cf.\,\cite{F1,F2,F3}). This leads on to the notion
of the Fukaya category (cf.\,\cite{Fu, F-O, FO1, FO2}).
The Fukaya category is related to the phenomenon of mirror symmetry, as formulated by
Kontsevich (cf.\,\cite{Ko}).
\end{remark}

\end{section}
\vskip 12mm
\begin{section}{{\bf Convexity properties of a moment map}}
\setcounter{equation}{0}
\vskip 3.5mm
Atiyah, Guillemin and Sternberg \cite{A, Gu-St1} proved the following {\sf Convexity Theorem} for a Hamiltonian
$\mathbb T^m$-space. Here $\mathbb T^m=\BR^m/\BZ^m$ is an $m$-dimensional torus.

\begin{theorem}\,\!\!\!{\sf (Convexity Theorem: Atiyah, Guillemin-Sternberg\,[1982])}
\label{thm:5.1}
Let $(M,\omega)$ be a compact connected symplectic manifold. Assume that $(M,\omega,\mathbb T^m,\mu)$
 is a Hamiltonian $\mathbb T^m$-space. Then the following properties {\rm (M1)-(M3)} are satisfied
\vskip 2mm\noindent
{\rm (M1)} the levels of $\mu$ are connected;
\vskip 2mm\noindent
{\rm (M2)} the image $\mu (M)$ of $\mu$ is convex;
\vskip 2mm\noindent
{\rm (M3)} the image of $\mu$ is the convex hull of the images of the fixed points of the \\
\indent \ \ \ \ action.
\end{theorem}
\begin{proof} The proof can be found in \cite{A, Gu-St1}.
We briefly sketch the proof of Atiyah\,\cite{A} (cf.\, \cite[pp.\,169--170]{C}).
\vskip 2mm\noindent
{\bf Claim\ 1.} {\it The levels of $\mu$ are connected for any ${\mathbb T}^m$-action, $m=1,2,\cdots$.}
\vskip 2mm
We leave the proof of Claim 1 to the reader.
\vskip 3mm\noindent
{\bf Claim\ 2.} {\it The image of $\mu$ are convex for any ${\mathbb T}^m$-action, $m=1,2,\cdots$.}
\vskip 2mm
Now we prove Claim 2. For a ${\mathbb T}^1$-action, $\mu (M)$ is convex because in $\BR$ connectivity is convexity. For a Hamiltonian ${\mathbb T}^m$-space $(M,\omega,{\mathbb T}^m,\mu)\ (m\geq 2)$, we first take a matrix $A\in \BZ^{(m,m-1)}$ of rank $m-1$. Let $\rho_A:{\mathbb T}^{m-1}\lrt {\rm Symp}(M,\omega)$ be the action of ${\mathbb T}^{m-1}$ on $(M,\omega)$ defined by
\begin{equation}\label{(5.1)}
  \rho_A (t)\!\cdot\! p=:(At)\!\cdot\! p,\quad t\in {\mathbb T}^{m-1},\ p\in M.
\end{equation}
Then $\rho_A$ is the Hamiltonian ${\mathbb T}^{m-1}$-action on $(M,\omega)$ with its moment map
$\mu_A:M\lrt \BR^{m-1}$ given by
\begin{equation*}
\mu_A (p):=\,{}^t\!A\mu (p),\quad p\in M.
\end{equation*}
For $\xi\in \BR^{m-1}$, choose $p_0\in \mu_A^{-1} (\xi).$ Then
\begin{equation*}
p\in \mu_A^{-1} (\xi)\quad {\rm if\ and\ only\ if}\quad {}^t\!A\mu (p)=\xi={}^t\!A\mu (p_0).
\end{equation*}
So we have
\begin{equation*}
\mu_A^{-1} (\xi)=\left\{ p\in M\,|\ \mu(p)-\mu(p_0)\in {\rm ker}\,({}^tA)\,\right\}.
\end{equation*}
According to Claim 1, $\mu_A^{-1} (\xi)$ is connected. For two points $p_0,p_1\in \mu_A^{-1} (\xi)$,
we take a curve $\alpha:[0,1]\lrt \mu_A^{-1} (\xi)$ with $\alpha (0)=p_0$ and $\alpha (1)=p_1$. Then we obtain a curve $\gamma:[0,1]\lrt  {\rm ker}\,({}^t\!A)\subset \BR^m$ defined by
\begin{equation*}
\gamma (t):=\mu (\alpha(t))-\mu (p_0),\quad t\in [0,1].
\end{equation*}
Since ${\rm ker}\,({}^t\!A)$ is a one-dimensional subspace of $\BR^m$,
\begin{equation*}
c \mu (p_0)+(1-c)\mu (p_1)\in \mu (M)\quad {\rm for\ any}\ c\in [0,1].
\end{equation*}
Any $p_0,p_1\in M$ can be approximated arbitrarily by points $q_0$ and $q_1$ in $M$ with
$\mu(q_0)-\mu(q_1)\in {\rm ker}\,({}^t\!A)$ for some matrix $A\in \BZ^{(m,m-1)}$ of rank $m-1$. Taking the limits $q_0\lrt p_0$ and $q_1\lrt p_1$, we see that $\mu (M)$ is convex. This completes the proof of the statement (M2).
\vskip 2mm
Let $N$ be the fixed point set of the ${\mathbb T}^{m}$-action $\rho$ on $(M,\omega)$. Then $N$ is a finite
disjoint union of connected symplectic submanifolds $N_1,\cdots,N_k$. The moment map $\mu$ is constant on each $N_j\,(1\leq j\leq k)$, say, $\mu(N_j)=c_j\in \BR^m\,(1\leq j\leq k)$. By (M2), the convex hull
${\rm Conv} (c_1,\cdots,c_k)$ is contained in $\mu (M)$. Conversely, suppose that $\xi\in \BR^m$ and
$\xi\notin {\rm Conv} (c_1,\cdots,c_k)$. Choose $\zeta\in\BR^m$ with rationally independent components such that
\begin{equation*}
  \langle \xi,\zeta \rangle > \langle c_j,\zeta \rangle,\quad j=1,\cdots,k.
\end{equation*}
Here $\langle\,\,,\,\,\rangle$ is a $\mathbb T^m$-invariant positive definite inner product on $\BR^m$. By the irrationality of $\zeta$, the set $\{ \exp t\zeta\,|\,t\in\BR\,\}$ is dense in $\mathbb T^m$. Hence the zeros of the vector field $X^\zeta$ on $M$ generated by $\zeta$ are the fixed points of the $\mathbb T^m$-action $\rho$. Since the function $\langle \mu(\,\cdot\,),\zeta \rangle $ attains its maximum value on one of the sets $N_j$, we have the relation
\begin{equation*}
  \langle \xi,\zeta \rangle > \sup_{p\in M} \langle \mu(p),\zeta \rangle.
\end{equation*}
So $\xi\notin \mu(M).$ Therefore $\mu(M)$ is contained in ${\rm Conv} (c_1,\cdots,c_k)$. This completes the proof of the statement (M3).
\end{proof}

\begin{theorem}\label{thm:5.2}
Let $M$ be a connected nonsingular complex projective variety with its K{\"a}hler form $\omega$ and $G$ a compact connected subgroup of the group of complex symplectic transformations of $(M,\omega)$. Suppose that for some point $p\in M$, the stabilizer $G_p$ of $p$ in $G$ is finite. Let $\mathfrak t$ be a Cartan subalgebra of the Lie algebra $\mathfrak g$ of $G$, and let ${\mathfrak t}^*$ be the subspace of the dual space ${\mathfrak g}^*$ of $\mathfrak g$ corresponding to $\mathfrak t$. Let ${\mathfrak t}^*_+$ be a Weyl chamber in ${\mathfrak t}^*$, and let $\mu:M\lrt {\mathfrak g}^*$ be the moment map. Then the intersection $\mu (M)\cap {\mathfrak t}^*$ is a convex polytope of dimension equal to the rank of $G$.
\end{theorem}
\begin{proof} The proof can be found in \cite[pp.\,511--513]{Gu-St1}.\end{proof}

\vskip 2mm
Kirwan \cite{Ki} generalized the above convexity theorem to the case of a Hamiltonian $G$-space where $G$ is a compact Lie group.

\begin{theorem}\,\!\!\!{\sf (Kirwan\,[1984])}\label{thm:5.3}
Let $G$ be a compact connected Lie group. Let $(M,\omega)$ be a compact connected symplectic manifold. Assume that $(M,\omega,G,\mu)$ is a Hamiltonian $G$-space. Let $\frak t^*_+$ be a positive Weyl chamber
$\frak t^*_+$ in the dual space $\frak t^*$ of the Lie algebra $\frak t$ of a maximal torus $T$ of $G$. Then the intersection of $\mu (M)\cap \frak t^*_+$ of the image of the moment map with a positive Weyl chamber in $\frak t^*$ is convex.
\end{theorem}
\begin{proof} We give a sketch of Kirwan's proof.
We first fix a $G$-invariant inner product $\langle \,,\,\rangle$ on $\mathfrak g$ and use it to identify $\mathfrak g^*$ with $\mathfrak g$, and $\mathfrak t^*$ with $\mathfrak t$. Here $\mathfrak t$ is the Lie algebra of a maximal torus $T$ of $G$. Let $\mathfrak t_+$ be a positive Weyl chamber in $\mathfrak t$.
Let $\|\,,\,\|$ be the associated norm on $\mathfrak g$.
\vskip 2mm\noindent
{\bf Claim\ 1.} {\it The subset of points of $M$ where $\|\mu\|^2$ takes its minimum value is connected.}
\vskip 2mm
We refer to \cite[(3.1)]{Ki} for more details of Claim 1.
\vskip 2mm
For any point $\alpha\in \mathfrak g$, there exists a natural $G$-invariant symplectic structure $\tau_\alpha$ on the coadjoint orbit $\mathcal O (\alpha):={\rm Ad}(G)\alpha$ of $\alpha$. Then
$\mathcal O (\alpha)$ has the form $G/H$, where $H$ is the centralizer of $\alpha$ in $G$. The inclusion map
\begin{equation*}
  \mu_H:G/H\hookrightarrow \mathfrak g \cong \mathfrak g^*
\end{equation*}
is the moment map on $G/H=\mathcal O (\alpha)$. Put $\omega_\alpha=-\tau_\alpha$. Clearly
$(\mathcal O (\alpha),\omega_\alpha)$ is a symplectic manifold. Then we see that $M\times G/H$ becomes a
symplectic manifold and its moment map
\begin{equation*}
\mu^{(\alpha)}:M\times G/H \lrt \mathfrak g
\end{equation*}
is given by
\begin{equation*}
\mu^{(\alpha)}(p,gH)=\mu (p)-{\rm Ad}(g)\alpha,\qquad p\in M,\ g\in G.
\end{equation*}
\vskip 2mm\noindent
{\bf Claim\ 2.} {\it For any sufficiently small $\varepsilon >0$, there exists an element
$\alpha\in \mathfrak t_+$ such that the ball $B(\alpha,\varepsilon)$ of radius $\varepsilon$ and center
$\alpha$ meets $\mu (M)\cap \mathfrak t_+$ in precisely two points $\alpha_1$ and $\alpha_2$ neither of
which lies in the interior of $B(\alpha,\varepsilon)$.}
\vskip 2mm
The proof of Claim 2 may be found in \cite[pp.\,549--551]{Ki}.
\vskip 2mm\noindent
{\bf Claim\ 3.} {\it The function $\| \mu^{(\alpha)}\|^2$ on $M\times G/H$ takes its minimum value precisely at those points $(x,gH)$ such that
\begin{equation*}
\mu (g^{-1}x)=\alpha_j,\qquad j=1,2.
\end{equation*}
Here $\alpha_1$ and $\alpha_2$ are two points in Claim 2. }
\vskip 2mm
The proof of Claim 3 may be found in \cite[pp.\,551--552]{Ki}.
\vskip 2mm
Using Claim 1, Claim 2 and Claim 3, Kirwan proved the above theorem as follows:
Suppose $\mu (M)\cap \mathfrak t_+$ is {\sf not} convex. Let $\| \mu^{(\alpha)}\|^2$ be the function
$M\times G/H$ where $\alpha$ satisfies the conditions of Claim 2. According to Claim 3, the set
\begin{equation*}
\left\{ (p,gH)\in M\times G/H\,|\ \mu(g^{-1}p)=\alpha_j,\ \ j=1,2\,\right\}
\end{equation*}
is the disjoint union of the following two non-empty closed subsets
\begin{equation*}
G \left( \mu^{-1}(\alpha_1)\times \{ H\}\right) \quad {\rm and}\quad
G \left( \mu^{-1}(\alpha_2)\times \{ H\}\right).
\end{equation*}
This contradicts Claim 1. Thus $\mu (M)\cap \mathfrak t_+$ is convex. Hence we complete the proof of the above theorem.
\end{proof}

\vskip 0.25cm
Though we do not have a classification of symplectic manifolds so far, fortunately we have a classification of symplectic-toric manifolds which are very special Hamiltonian torus-spaces in terms of combinatorial data. We recall that a {\sf symplectic-toric manifold} is a compact connected symplectic manifold $(M,\omega)$ of dimension $2n$ equipped with an effective Hamiltonian action of an $n$-dimensional torus $\mathbb T^n$ and with a corresponding moment map $\mu:M\lrt \BR^n.$ It can be seen that every symplectic-toric manifold is simply connected. For instance,
according to Example (2) in Section 3, $(S^2, d\theta\wedge dh, S^1)$ is a symplectic-toric manifold.
In 1988 T. Delzant classified all symplectic-toric manifolds in terms of a set of very special polytopes. We describe his classification roughly.

\begin{definition}\label{def:5.4}
A {\sf Delzant polytope} $\Delta$ in $\BR^n$ is a convex polytope satisfying the following properties
{\rm (DP1)-(DP3)}:
\vskip 0.25cm \noindent
{\rm (DP1)} it is {\sf simple}, i.e., there are $n$ edges meeting at each vertex;
\vskip 0.25cm \noindent
{\rm (DP2)} it is {\sf rational}, i.e., the edges meeting at the vertex $p$ are rational in the sense
\par \ \ \ \ \
that each edge is of the form $p+t\,\alpha_i,\ t\geq 0,$ where $\alpha_i\in \BZ^n \ (1\leq i\leq n);$
\vskip 0.25cm \noindent
{\rm (DP3)} it is {\sf smooth}, i.e., for each vertex, the corresponding $\alpha_1,\cdots,\alpha_n$
can be chosen
\par \ \ \ \ \
to be a $\BZ$-basis of $\BZ^n$.
\end{definition}
\vskip 3mm
First we give several definitions for the reader's convenience. A {\sf facet} of a polytope $\Delta$ with
$\dim \Delta=n$ in $\BR^n$ is a $(n-1)$-dimensional face. Let $\Delta$ be a Delzant polytope with
$\dim \Delta=n$ and $d=$\,the number of facets of $\Delta$. A lattice vector $v\in\BZ^n$ is said to be
{\sf primitive} if it cannot be written as $v=ku$ with $u\in \BZ^n,\ k\in\BZ$ with $|k|>1.$
\begin{theorem}\,\!\!{\sf (Delzant \cite{De}: 1988)}\label{thm:5.5}
Symplectic-toric manifolds are classified by Delzant polytopes. More precisely, there is the one-to-one correspondence between the set $\frak A$ of all symplectic toric manifolds of dimension $2n$ and
the set $\frak B$ of all Delzant polytopes in $\BR^n$ given by
\begin{equation}
 \frak A \ni (M,\omega,\mathbb T^n,\mu) \mapsto \mu (M)\in \frak B.
\end{equation}
\end{theorem}
\begin{proof} We give a sketchy proof of the ``if" part following \cite{Gu} (cf.\,\cite{C}).
Let $\Delta$ be a Delzant polytope with $d$ facets. Let $v_i\in \BZ^n\ (i=1,2,\cdots,d)$ be the primitive outward-pointing normal vectors to the facets. Then
\begin{equation*}
  \Delta=\{ x\in (\BR^n)^\times\,|\ \langle x,v_i \rangle \leq \lambda_i,\ \, 1\leq i\leq d\ \} \quad
  {\rm for\ some}\ \lambda_i\in \BR\ (1\leq i\leq d).
\end{equation*}
Let $\{ e_1,\cdots,e_d \}$ be the standard basis of $\BR^d$. If $\pi:\BR^d\lrt \BR^n$ be the map defined by
$\pi (e_i)=v_i,\ i=1,2,\cdots,d,$ then it is easily seen that $\pi$ is surjective and
$\pi (\BZ^d)=\BZ^n.$ Thus $\pi$ induces a surjective group homomorphisim
\begin{equation*}
  \theta:\mathbb T^d(:=\BR^d/\BZ^d)\lrt \mathbb T^n(:=\BR^n/\BZ^n)
\end{equation*}
between $\mathbb T^d$ and $\mathbb T^n$. Let $N$ be the kernel of $\theta$ with its Lie algebra
$\mathfrak n$. The exact sequence of tori
\begin{equation*}
0 \lrt N \stackrel{i}{\lrt}\mathbb T^d  \stackrel{\theta}{\lrt} \mathbb T^n \lrt 0
\end{equation*}
induces an exact sequence of Lie algebras
\begin{equation*}
0 \lrt \mathfrak n \stackrel{i}{\lrt} \BR^d \stackrel{\theta}{\lrt} \BR^n \lrt 0
\end{equation*}
with its dual exact sequence
\begin{equation*}
0 \lrt {\rm Hom} (\BR^n,\BR) \stackrel{\theta^*}{\lrt} {\rm Hom} (\BR^d,\BR) \stackrel{i^*}{\lrt}
{\rm Hom}(\mathfrak n,\BR) \lrt 0.
\end{equation*}
Let $(\BC^d,\omega_0,\mathbb T^d,\mu)$ be a symplectic toric manifold equipped with the standard Hamiltonian action of $\mathbb T^d$ on $\BC^n$:
\begin{equation*}
(\xi^{t_1},\cdots,\xi^{t_d})\cdot (z_1,\cdots,z_d):=(\xi^{t_1}z_1,\cdots,\xi^{t_d}z_d),\quad
\xi:=e^{2\pi i}
\end{equation*}
and
\begin{equation*}
\mu:\BC^d \lrt (\BR^d)^*,\qquad \mu (z_1,\cdots,z_d):=-\pi (|z_1|^2,\cdots,|z_d|^2)+(\lambda_1,\cdots,
\lambda_d).
\end{equation*}
We consider the map
\begin{equation*}
i^*\circ \mu: \BC^d \lrt {\rm Hom} (\BR^d,\BR)\lrt \mathfrak n^*:={\rm Hom}(\mathfrak n,\BR)
\end{equation*}
and put
\begin{equation*}
Z:=(i^*\circ \mu)^{-1}(0),\quad {\rm the\ zero\!-\!level\ set\ of}\ i^*\circ \mu.
\end{equation*}
Then we can show that$Z$ is compact and $N$ acts on $Z$ freely. Thus $p:Z\lrt M_\Delta:=Z/N$ is the principal $N$-bundle on $M_\Delta$. According to the Marsden-Weinstein-Meyer theorem, there exists a
symplectic form $\omega_\Delta$ on $M_\Delta$ such that
\begin{equation*}
p^*\omega_\Delta=j^*\omega_0,
\end{equation*}
where $j:Z\hookrightarrow \BC^d$ is the inclusion. Therefore $(M_\Delta,\omega_\Delta)$ is a compact
symplectic manifold of dimension $2n$. Furthermore we can show that the torus $\mathbb T^n$ acts on
$M_\Delta$ in a Hamiltonian fashion and its moment map $\mu_n:M_\Delta \lrt {\rm Hom} (\BR^n,\BR)$ such that $\mu_n (M_\Delta)=\Delta$ (we refer to \cite[pp.\,185-186]{C} for the detailed proof). Finally
the quadruple $(M_\Delta,\omega_\Delta,\mathbb T^n,\mu_n)$ is the required symplectic toric manifold of
dimension $2n$ corresponding to $\Delta\subset {\rm Hom} (\BR^n,\BR)=(\BR^n)^*$.
\end{proof}

\begin{remark}\label{rk:5.6}
Guillemin, Miranda, Pires and Scott proved the analogue of Theorem 5.2 for log symplectic-toric manifolds which are defined to be generically symplectic-toric and degenerate along a normal crossing configuration of smooth hypersurfaces. Log symplectic-toric manifolds belong to a class of Poisson manifolds. Most often degeneracy loci for Poisson structures are singular.
\end{remark}

\begin{theorem}\!\!\!{\sf (Ahara and Hattori, Audin)}\label{thm:5.7}
Suppose $(M,\omega,S^1)$ is a compact connected symplectic $4$-dimensional manifold equipped with an effective Hamiltonian $S^1$-action. Then $(M,\omega,S^1)$ is $S^1$-equivariantly diffeomorphic to a complex surface with a holomorphic $S^1$-action which is obtained from $\mathbb P^2 (\BC)$, a Hirzebruch surface, or a $\mathbb P^1 (\BC)$-bundle over a Riemann surface with appropriate circle actions by a sequence of blow-ups at the fixed points.
\end{theorem}
We refer to \cite{A-H, Au1, Au2} for more details.

\end{section}
%%%%%%%%%%%%%%%%%%%%%%%%%%%%%%%%%%%%%%%%%%%%%%%%%%%%%%%%%%%%%%%%%%%%%%%%%%%%%%%%%%%%%%%%%%%%%%%%%%%%%%%%
%%%%%%%%%%%%%%%%%%%%%%%%%%%%%%%%%%%%%%%%%%%%%%%%%%%%%%%%%%%%%%%%%%%%%%%%%%%%%%%%%%%%%%%%%%%%%%%%%%%%%%%%
%%%%%%%%%%%%%%%%%%%%%%%%%%%%%%%%%%%%%%%%%%%%%%%%%%%%%%%%%%%%%%%%%%%%%%%%%%%%%%%%%%%%%%%%%%%%%%%%%%%%%%%%
%%%%%%%%%%%%%%%%%%%%%%%%%%%%%%%%%%%%%%%%%%%%%%%%%%%%%%%%%%%%%%%%%%%%%%%%%%%%%%%%%%%%%%%%%%%%%%%%%%%%%%%%
%%%%%%%%%%%%%%%%%%%%%%%%%%%%%%%%%%%%%%%%%%%%%%%%%%%%%%%%%%%%%%%%%%%%%%%%%%%%%%%%%%%%%%%%%%%%%%%%%%%%%%%%
%%%%%%%%%%%%%%%%%%%%%%%%%%%%%%%%%%%%%%%%%%%%%%%%%%%%%%%%%%%%%%%%%%%%%%%%%%%%%%%%%%%%%%%%%%%%%%%%%%%%%%%%
%%%%%%%%%%%%%%%%%%%%%%%%%%%%%%%%%%%%%%%%%%%%%%%%%%%%%%%%%%%%%%%%%%%%%%%%%%%%%%%%%%%%%%%%%%%%%%%%%%%%%%%%
%%%%%%%%%%%%%%%%%%%%%%%%%%%%%%%%%%%%%%%%%%%%%%%%%%%%%%%%%%%%%%%%%%%%%%%%%%%%%%%%%%%%%%%%%%%%%%%%%%%%%%%%

\vskip 0.5cm \noindent
\begin{section}{{\bf Modern theory of symplectic geometry}}
\setcounter{equation}{0}
\vskip 2mm
Symplectic geometry and symplectic topology have led the study of the classification of
symplectic and Hamiltonian actions, the study of symplectic embeddings and the theory of
Gromov-Witten invariants from the early 2000s to the present day. In particular,
Gromov's work~\cite{Gr1} stimulated those mathematicians who have researched the above-mentioned
subjects. In the final section, we list more or less recent results obtained by many mathematicians.
%%%%%%%%%%%%%%%%%%%%%%%%%%%%%%%%%%%%%%%%%%%%%%%%%%%%%%%%%%%%%%%%%%%%%%%%%%%%%%%%%%%%%%%%%%%%%%%%%%%%%%%%
%%%%%%%%%%%%%%%%%%%%%%%%%%%%%%%%%%%%%%%%%%%%%%%%%%%%%%%%%%%%%%%%%%%%%%%%%%%%%%%%%%%%%%%%%%%%%%%%%%%%%%%%
%%%%%%%%%%%%%%%%%%%%%%%%%%%%%%%%%%%%%%%%%%%%%%%%%%%%%%%%%%%%%%%%%%%%%%%%%%%%%%%%%%%%%%%%%%%%%%%%%%%%%%%%
%%%%%%%%%%%%%%%%%%%%%%%%%%%%%%%%%%%%%%%%%%%%%%%%%%%%%%%%%%%%%%%%%%%%%%%%%%%%%%%%%%%%%%%%%%%%%%%%%%%%%%%%
%\vskip 5mm \noindent
\subsection{The classification problem of symplectic group actions}
\leavevmode\par\vspace{2mm}
%\noindent
In this subsection, we briefly describe the developments in the study of Hamiltonian symplectic
group actions and non-Hamiltonian group actions the past thirty years.
We note that the phase space of a mechanical system is modelled by a symplectic manifold, and its symmetries are described by symplectic group actions.

\vskip 2mm
Let $G$ be a Lie group acting on a smooth manifold $M$. For $x\in M$ we let $G_x$
be the stabilizer subgroup of the $G$-action at $x$. The $G$-action is said to be {\sf free}
if $G_x=\{ e\}$ for every $x\in M.$ The $G$-action is called {\sf semifree} if for
every $x\in M$ either $G_x=\{ e\}$ or $G_x=G$. The $G$-action is said to be {\sf effective}
if $\bigcap_{x\in M}G_x=\{ e\}$. We say that the $G$-action is {\sf proper} if
the map $f:G\times M\lrt M\times M$ given by $f(g,x):=(gx,x)\ (x\in M)$ is proper.

\begin{definition}\label{def:6.1}
Let $(V,\omega)$ be a symplectic vector space. A subspace $W$ of $V$ is called
${\sf coisotropic}$ if $W^\omega \subseteq W.$ Here,
$$ W^\omega:=\left\{ v\in V\,\vert\ \omega (v,w)=0\ \ {\rm for\ all}\ w\in W\right\}.$$
A submanifold $N$ of a symplectic manifold $(M,\omega)$ is called ${\sf coisotropic}$
if $T_x N$ is a coisotropic subspace of $(T_x M,\omega_x)$ for every $x\in N.$
Here, $T_x N$ and $T_x M$ denote the tangent space of $N$ and $M$ at $x\in N$
respectively. We say that a submanifold $N$ of a symplectic manifold $(M,\omega)$ is
${\sf symplectic}$ if $T_x N$ is a symplectic subspace of $(T_x M,\omega_x)$
for every $x\in N.$ A {\sf symplectic--toric manifold} is a compact connected
symplectic manifold $(M,\omega)$ of dimension $2\,n$ endowed with an effective
Hamiltonian action of a {\it torus} $T$ of dimension $n$.
\end{definition}

\begin{definition}\label{def:6.2}
Let $\circ:G\times M\lrt M$ and $*:G\times N\lrt N$ be smooth $G$-actions, and
$\phi:M\lrt N$ be a smooth map. The map is said to be a
$G$--{\sf equivariant diffeomorphism} if it is a diffeomorphism and satisfies
the following condition
$$
\phi (g\circ x)=g* \phi (x)\qquad {\rm{for \ all}}\ g\in G\ \textrm{and}\ x\in M.
$$
\end{definition}

\begin{definition}\label{def:6.3}
A symplectic manifold with a Hamiltonian action of a compact Lie group is called
{\sf a multiplicity-free space} if the Poisson bracket of any two invariant
smooth functions vanishes.
\end{definition}
\begin{definition}\label{def:6.4}
Suppose that a torus $T$ acts effectively and symplectically on a compact connected symplectic manifold.
The $T$-action is said to be {\sf coisotropic} if there is a coisotropic $T$-orbit.
\end{definition}
\begin{definition}\label{def:6.5}
Suppose that a torus $T$ acts effectively and symplectically on a compact connected symplectic manifold.
If there is a $\dim T$-dimensional symplectic $T$-orbit, we say that the $T$-action is a
{\sf maximal\ symplectic\ action}.
\end{definition}
Over past thirty years the theory of the symplectic actions on symplectic manifolds has been developed by some experts in symplectic geometry.
In the good survey article\,\cite{P2}, {\'A}lvaro Pelayo described classifications on compact connected symplectic manifolds $(M,\omega)$:
\vskip 2mm
(a) ``Maximal Hamitonian case": Hamiltonian $T$-action, $\dim M=2\,\dim T.$
\vskip 2mm
(b) ``$S^1$-Hamiltonian case": Hamiltonian $T$-action, $\dim M=4\,\ \dim T=1.$
\vskip 2mm
(c) ``Four-dimensional case": $\dim M=4\,\ \dim T=2.$
\vskip 2mm
(d) ``Maximal symplectic case": there is a $\dim\ T$-orbit symplectic orbit.
\vskip 2mm
(e) ``Coisotropic case" there is a coisotropic $T$-orbit.
\vskip 2mm
Here $T$ denotes a torus. He outlined connections of these works with algebraic geometry, toric varieties, log-symplectic toric geometry, torus bundles over tori, nilpotent Lie groups, integral systems and the classification of semi-toric systems.
\vskip 2mm
Let $S^1=U(1)$ be a torus of dimension one.
A Hamiltonian $S^1$-action on a compact connected symplectic manifold $(M,\omega)$ of dimension $2n$ has at least $n+1$ fixed points. In fact, the number of fixed points is $\sum_{k=0}^{2n} {\rm rank}~ H^k (M,\BR)$ and $0\neq [\omega^k]\in H^{2k}(M,\BR)$ for $1\leq k\leq n.$  The following natural questions arise:
\vskip 3mm\noindent
\begin{question}\label{que:6.6}
Under which conditions is a symplectic $G$-action Hamiltonian ?  Describe the obstruction to being Hamiltonian. Here $G$ is an $m$-dimensional torus, a compact connected Lie group, a non-compact
abelian group or a semisimple Lie group.
\end{question}
\begin{question}\label{que:6.7}
Are there non-Hamiltonian symplectic $S^1$-actions on compact connected symplectic manifolds
with nonempty discrete fixed point sets ?
\end{question}
\vskip 1mm
Some affirmative answers were given as follows.
\begin{theorem}\!\!{\sf (Frankel \cite{Fr}:~1959)}\label{thm:6.8}
Let $(M,\omega)$ be a compact connected K{\"a}hler manifold admitting an $S^1$-action preserving the
K{\"a}hler structure $\omega$. If the the $S^1$-action has some fixed points, then it is Hamiltonian.
\end{theorem}

\begin{theorem}\!\!{\sf (McDuff \cite{Mc1}:~1988)}\label{thm:6.9}
  A symplectic $S^1$-action on a compact connected symplectic $4$-manifold with some fixed point is
  Hamiltonian.
\end{theorem}

\begin{theorem}\!\!{\sf (Tolman and Weitsman \cite{T-W}:~2000)}\label{thm:6.10}
  Let $(M,\omega)$ be a compact connected symplectic manifold equipped with a semi-free symplectic $S^1$-action with isolated fixed points. If there is at least one fixed point, the $S^1$-action is
  Hamiltonian.
\end{theorem}

\begin{theorem}\!\!{\sf (Giacobbe \cite{Gi}:~2005)}\label{thm:6.11}
An effective symplectic action of an $n$-dimensional torus on a compact connected symplectic $2n$-dimensional manifold with some fixed point must be Hamiltonian.
\end{theorem}

\begin{theorem}\!\!{\sf (Tolman \cite{To}:~2017)}\label{thm:6.12}
There exists a symplectic non-Hamiltonian $S^1$-action on a compact, connected, six-dimensional
symplectic manifold with exactly $32$ fixed points.
\end{theorem}
Since the above result is important, we provide a sketchy proof.
\begin{proof}
  We give a sketch of Tolman's proof. First we need some additional terminology.
Let’s now consider a non-Hamiltonian symplectic circle action on a closed
$2n$-dimensional symplectic manifold $(M,\omega)$. If $[\omega]\in H^2(M,\BR)$
is integral, there exists a {\sf generalized moment map}, that is,
a map $\Phi:M\lrt S^1$ satisfying $\Phi^*(dt)=-\iota_{\xi_M}\omega$.
Given a proper open subset $ U\subseteq S^1$,  the $S^1$--action on the
preimage $\Phi^{-1}(U)$ is Hamiltonian.
Hence, as in the Hamiltonian case, the {\sf symplectic quotient}
$M/\!\!/\!_t S^1 := \Phi^{-1}(t)/S^1$ is a $(2n-2)$--dimensional orbifold endowed with a symplectic form
$\omega_t\in \Omega^2 (M/\!\!/\!_t S^1)$ for all regular values $t$ of $\Phi$.
Since the action is non-Hamiltonian, every fixed point has at least
two non-zero weights. Hence, there exists a unique continuous real-valued function
$\varphi:S^1\lrt \BR$ such that
\begin{equation*}
  \varphi (t)=\int_{M/\!\!/\!_t S^1}\omega_t^{n-1}
\end{equation*}
for all regular values $t$ of function of $M$. We call $\varphi$ the
{\sf Duistermaat}--{\sf Heckman} function of $M$.
\vskip 2mm
A $K3$ surface is a closed, connected complex surface $(\widehat{X},\widehat{I})$ with
$H^1(\widehat{X},\BR)=\{0\}$
and trivial canonical bundle. A closed complex surface $(X,I)$ with at worst simple
singular points (These are alternatively called {\it Du Val singularities} or
{\it rational double points}) and with minimal resolution $q:\widehat{X}\lrt X$
is a {\sf generalized $K3$ surface}
if $\widehat{X}$ is a $K3$ surface. The quotient of the complex torus
$T:=\BC^2/(\BZ^2+\sqrt{-1}\,\BZ^2)$
by the involution $[z]\mapsto [-z]$ is a
generalized $K3$ surface with $16$ isolated $\BZ_2$ singularities, called
the Kummer surface.
A symplectic form
$\s\in \Omega^2(\widehat{X})$ is {\sf tamed} if $\s(v,\widehat{I}(v))>0$
for every nonzero tangent vector $v$. In this case, we say that
the triple $(\widehat{X},\widehat{I},\s)$ is a {\sf tame $K3$ surface}.
A tamed symplectic form $\s\in \Omega^2(\widehat{X})$ is {\sf Kähler}
if, additionally, $\s(\widehat{I}(v),\widehat{I}(w))=\s (v,w)$  for all tangent
vectors $v$ and $w$.

\vskip 2mm\noindent
{\bf Step 1.} Tolman shows that there exists a free $S^1$-action on a symplectic manifold
$(\widehat{M},\widehat{\omega})$ with proper moment map $\widehat{\Phi}:\widehat{M}\lrt (1,3)$
so that for all $t\in (1,3)$, the symplectic quotient $\widetilde{M}/\!\!/\!_t S^1$ is
symplectomorphic
to a tame $K3$ surface $(\widehat{X},\widehat{I}_t,\widehat{\s}_t)$; moreover,
\vskip 2mm
(a) $\int_{\widehat{X}} (\widehat{\s}_t)^2=-4+16t-4t^2$\,; and
\vskip 2mm
(b) $[\widehat{\s}_t]=\widehat{\kappa}-t\widehat{\eta}$, where
$\widehat{\kappa},\ \widehat{\eta}$ induce a primitive embedding
$\BZ^2 \hookrightarrow H^2(X,\BZ).$
%%%%%%%%%%%%%%%%%%%%%%%%%%%%%%%%%%%%%%%%%%%%%%%%%%%%%%%%%%%%%%%%%%%%%%%%%%%%%%%%%%%%%%%%%%%%%%%%%%%%%%%%
%%%%%%%%%%%%%%%%%%%%%%%%%%%%%%%%%%%%%%%%%%%%%%%%%%%%%%%%%%%%%%%%%%%%%%%%%%%%%%%%%%%%%%%%%%%%%%%%%%%%%%%%
%%%%%%%%%%%%%%%%%%%%%%%%%%%%%%%%%%%%%%%%%%%%%%%%%%%%%%%%%%%%%%%%%%%%%%%%%%%%%%%%%%%%%%%%%%%%%%%%%%%%%%%%
%%%%%%%%%%%%%%%%%%%%%%%%%%%%%%%%%%%%%%%%%%%%%%%%%%%%%%%%%%%%%%%%%%%%%%%%%%%%%%%%%%%%%%%%%%%%%%%%%%%%%%%%
%%%%%%%%%%%%%%%%%%%%%%%%%%%%%%%%%%%%%%%%%%%%%%%%%%%%%%%%%%%%%%%%%%%%%%%%%%%%%%%%%%%%%%%%%%%%%%%%%%%%%%%%
%%%%%%%%%%%%%%%%%%%%%%%%%%%%%%%%%%%%%%%%%%%%%%%%%%%%%%%%%%%%%%%%%%%%%%%%%%%%%%%%%%%%%%%%%%%%%%%%%%%%%%%%
\vskip 2mm\noindent
{\bf Step 2.} Tolman constructs the symplectic manifold
$(\widetilde{M},\widetilde{\omega})$ with locally free $S^1$-action and proper moment map $\widetilde{\Phi}:\widetilde{M}\lrt (-1-\widetilde{\varepsilon},1+\widetilde{\varepsilon})$,
where $\widetilde{\varepsilon}>0$ such that the level sets $\widetilde{\Phi}^{-1}(\pm 1)$ each
contains exactly 16 fixed points with weights $\pm \{-2,1,1\}.$ Otherwise, the action is
locally free. The Duistermaat-Heckeman function of $\widetilde{M}$ is given by
\begin{equation*}
 \widetilde{\varphi}=\begin{cases}
                       -4-16t-4t^2, & \mbox{if}\ t\in (-1-\widetilde{\varepsilon},-1]  \\
                       \ 4+4t^2, & \mbox{if }\ t\in [-1,1] \\
                       -4+16t-4t^2, & \mbox{if}\ [1,1+\widetilde{\varepsilon}).
                     \end{cases}
\end{equation*}
The symplectic quotient $\widetilde{M}/\!\!/\!_t S^1$ is diffeomorphic to the Kummer
surface $T/\BZ_2$ for all $t\in (-1,-1),$ and symplectomorphic
to a tame $K3$ surface $(\widetilde{X},\widetilde{I},({\widetilde{\s}_\pm})_t)$
for all $t\in \pm (1,1+\widetilde{\varepsilon})$; moreover,
$[({\widetilde{\s}}_\pm)_t]=\widetilde{\kappa}-t\widetilde{\eta}_\pm$, where
$\widetilde{\kappa},\ {\widetilde{\eta}}_\pm$ induce a primitive embedding
$\BZ^2 \hookrightarrow H^2(\widetilde{X},\BZ).$
%%%%%%%%%%%%%%%%%%%%%%%%%%%%%%%%%%%%%%%%%%%%%%%%%%%%%%%%%%%%%%%%%%%%%%%%%%%%%%%%%%%%%%%%%%%%%%%
%%%%%%%%%%%%%%%%%%%%%%%%%%%%%%%%%%%%%%%%%%%%%%%%%%%%%%%%%%%%%%%%%%%%%%%%%%%%%%%%%%%%%%%%%%%%%%%
\vskip 2mm\noindent
{\bf Step 3.} Combining Step 1 and Step 2, we see that
\begin{eqnarray*}
  \int_{\widetilde{X}} \left( ({\widetilde{\s}_+})_t \right)^2&=& \widetilde{\varphi}(t)
  =\int_{\widehat{X}} (\widehat{\s}_t)^2 \quad {\rm for\ all}\ t\in
  t\in  (1,1+\widetilde{\varepsilon});\quad {\rm and}\\
 \int_{\widetilde{X}} \left( ({\widetilde{\s}_-})_t \right)^2&=& \widetilde{\varphi}(t)
  =\int_{\widehat{X}} (\widehat{\s}_{t+4})^2 \quad {\rm for\ all}\ t\in
  t\in  (-1-\widetilde{\varepsilon},-1).
\end{eqnarray*}
It can be shown that there exists $\epsilon$ and $\widehat{\epsilon}$ with
$0< \epsilon<\widehat{\epsilon}< \widetilde{\varepsilon}$ so that there are equivariant
symplectomorphisms
$$
\widetilde{M}\supset \widetilde{\Phi}^{-1}(1+\epsilon,1+\widehat{\epsilon})\lrt
\widehat{\Phi}^{-1}(1+\epsilon,1+\widehat{\epsilon})\subset \widehat{M}
$$
and
$$
\widetilde{M}\supset \widetilde{\Phi}^{-1}(-1-\widehat{\epsilon},-1+{\epsilon})\lrt
\widehat{\Phi}^{-1}(3-\widehat{\epsilon},3-{\epsilon})\subset \widehat{M}.
$$
Now you can use these two equivariant symplectomorphisms to glue together
$\widetilde{\Phi}^{-1}(-1-\widehat{\epsilon},1+\widehat{\epsilon})\subset \widetilde{M}$
and $\widehat{\Phi}^{-1}(1+\epsilon,3-{\epsilon})\subset \widehat{M}$, thus constructing
$S^1$-action on a symplectic manifold $(M,\omega)$ with a generalized moment map
$\Phi:M\lrt \BR/(4)$. In particular, the action is {\it not Hamiltonian} because there are
no fixed points with all positive weights or negative weights.
\end{proof}

\vskip 1mm
Recently Jang and Tolman improved Theorem \ref{thm:6.12} by reducing the number of fixed points.
\begin{theorem}\!\!{\sf (Jang and Tolman \cite{J-T}:~2023)}\label{thm:6.13}
Given an integer $k\geq 5$, there exists a non-Hamiltonian symplectic $S^1$-action
on a closed connected six-dimensional symplectic manifold with exactly $2k$ fixed points.
\end{theorem}

\vskip 2mm
We present some results on the classification of symplectic actions.
\begin{theorem}\!\!{\sf (Duistermaat and Pelayo \cite{D-P}:~2007)}\label{thm:6.14}
Compact connected symplectic manifolds $(M,\omega)$ with a coisotropic $T$-action are classified up to $T$-equivariant symplectomorphisms by symplectic invariants: the fundamental form $\omega^{\frak t}$, the Hamiltonian torus $T_h$ and its associated polytope $\Delta$, the period lattice $P$ of
$N=(\frak l / \frak t_h)^*$, the Chern class $c:N\times N \lrt \frak l$ of $M_{\rm reg}\lrt M_{\rm reg}/T$, and the holonomy invariant $[\tau:P\lrt T]_B\in {\rm Hom}_c (P,T)/B.$ Moreover, for any such list
$\mathcal L$ of five invariants there exists a compact connected symplectic manifolds $(M,\omega)$ with a coisotropic $T$-action with list of invariants $\mathcal L$.
\end{theorem}
{\'A}lvaro Pelayo \cite{P2} proposed the following natural classification problem:
\vskip 2mm\noindent
\begin{problem}\label{prob:6.15}
Let $G$ be an $m$-dimensional compact connected Lie group. Construct symplectic invariants and classify, up to equivariant symplectomorphisms, effective symplectic $G$-actions on compact connected symplectic $2n$-dimensional manifolds $(M,\omega)$ in terms of these invariants.
\end{problem}
\vskip 2mm
Pelayo\,\cite{P1, P2} proved the following result.
\begin{theorem}\!\!{\sf (Pelayo \cite{P1}:~2010\ {\rm and}\ \cite{P2}:~2017)}\label{thm:6.16}
Let $(M,\omega)$ be a compact connected symplectic $4$-manifold equipped with an effective action
of a $2$-torus $T$. If the symplectic $T$-action is Hamiltonian, then:
\vskip 2mm
(1) $(M,\omega)$ is a symplectic-toric manifold, so classified up to $T$-equivariant \\
\indent \ \ \ \ symplectomorphisms by the image of the moment map $\mu:M\lrt \frak t^*$ of the \\
\indent \ \ \ \ $T$-action.
\vskip 3mm\noindent
If the symplectic $T$-action is not Hamiltonian, then one and only one of the following cases occurs:
\vskip 2mm
(2) $(M,\omega)$ is equivariantly symplectomorphic to $((\BR/\BZ)^2\times S^2,\omega_*).$
The symplectic
\indent \ \ \ \
form $\omega_*$ on $(\BR/\BZ)^2 \times S^2$ is given
by
$$ \omega_*=dx\wedge dy + \omega_{\rm st},$$
\indent \ \ \ \
where $\omega_{\rm st}=d\theta\wedge dh$ is defined in Example (2) in Section 3 (see p.\,14).

\vskip 2mm
(3) $(M,\omega)$ is equivariantly symplectomorphic to $(T\times \frak t^*)/Q$ with the induced form \\
\indent \ \ \ \ and the $T$-action, where $Q\leq T\times \frak t^*$ is a discrete cocompact subgroup for \\ \indent \ \ \ \ the group structure on $T\times \frak t^*$.
\vskip 2mm
(4) $(M,\omega)$ is equivariantly symplectomorphic to a symplectic orbifold bundle
$$P:=\widetilde{\Sigma}\times_{\pi_1^{orb} (\sum,p_0)} T$$
\indent \ \ \ \
over a good orbisurface $\Sigma$, with symplectic form and $T$-action induced by \\
\indent \ \ \ \ the product ones. Here, in order to form the quotient $P$, the orbifold fundamental\\
\indent \ \ \ \  group $\pi_1^{orb} (\Sigma)$ acts on
$\widetilde{\Sigma}\times T$ diagonally, and on $T$ by means of a homomorphism\\
\indent \ \ \ \
$\mu:\pi_1^{orb} (\Sigma)\lrt T.$
\end{theorem}
%%%%%%%%%%%%%%%%%%%%%%%%%%%%%%%%%%%%%%%%%%%%%%%%%%%%%%%%%%%%%%%%%%%%%%%%%%%%%%%%%%%%%%%%%%%%%%%
%%%%%%%%%%%%%%%%%%%%%%%%%%%%%%%%%%%%%%%%%%%%%%%%%%%%%%%%%%%%%%%%%%%%%%%%%%%%%%%%%%%%%%%%%%%%%%%
%%%%%%%%%%%%%%%%%%%%%%%%%%%%%%%%%%%%%%%%%%%%%%%%%%%%%%%%%%%%%%%%%%%%%%%%%%%%%%%%%%%%%%%%%%%%%%%
%%%%%%%%%%%%%%%%%%%%%%%%%%%%%%%%%%%%%%%%%%%%%%%%%%%%%%%%%%%%%%%%%%%%%%%%%%%%%%%%%%%%%%%%%%%%%%%
%%%%%%%%%%%%%%%%%%%%%%%%%%%%%%%%%%%%%%%%%%%%%%%%%%%%%%%%%%%%%%%%%%%%%%%%%%%%%%%%%%%%%%%%%%%%%%%
%%%%%%%%%%%%%%%%%%%%%%%%%%%%%%%%%%%%%%%%%%%%%%%%%%%%%%%%%%%%%%%%%%%%%%%%%%%%%%%%%%%%%%%%%%%%%%%
\vskip 8mm
\subsection{The symplectic embedding problems}
\leavevmode\par\vspace{2mm}
In this subsection, we give a brief description of the developments in the study of
the symplectic embedding problems. J-holomorphic curves and the Taubes-Seiberg-Witten
theory have played an important role in the study of the symplectic embedding.
\vskip 2mm
Symplectic embeddings of simple shapes, such as (collections of)
balls, ellipsoids, and cubes, lie at the heart of symplectic geometry ever since
Gromov’s seminal Nonsqueezing Theorem from 1985. Symplectic embedding results
give a feeling for what symplectic means, and together with the techniques used
in their proofs lead to new connections to other fields, including those mentioned
above. After a very fruitful decade of research starting around 1989, not too much
happened in the subsequent decade. But since 2008 there has been much progress
in old and new questions on symplectic embeddings. This “third revolution” was
instigated by two ingenious constructions by Guth \cite{Gut} and
McDuff \cite{Mc2} (cf.\,\cite[p.\,140]{Sc}). We refer to \cite[\S 4, pp.\,152--158]{Sc}
for the background and motivations of the study of the symplectic embeddings.

\vskip 2mm
Given two open subsets $U$ and $V$ in $\BR^n$, we often write $U\stackrel{s}{\hookrightarrow}V$ instead of ``there exists a symplectic embedding of $U$ into $V$".
That is, if $U,V$ are open subsets of $\BR^n$, a symplectic
embedding ``$f:U\stackrel{s}{\hookrightarrow} V$" is a smooth embedding such that
$f^*\omega_{\BR^n}=\omega_{\BR^n}$. Similarly one defines symplectic embeddings
$f:(M,\omega)\stackrel{s}{\hookrightarrow} (M',\omega')$ between general symplectic
manifolds.
We denote by $D(a)$ the open disk in $\BR^2$ of area $a$, centered at the origin, and $P(a_1,\cdots,a_n)=D(a_1)\times \cdots\times D(a_n)$ the open polydisk in $\BR^{2n}$. We let $C^{2n}(a):=P(a,\cdots,a)$ be the cube in $\BR^{2n}$ and let
$Z^{2n}(a):=D(a)\times \BC^{n-1}$ be the symplectic cylinder. Let
\begin{equation*}
  E(a_1,\cdots,a_n):=\left\{ (z_1,\cdots,z_n)\in\BC^n \big|\ \,\sum_{i=1}^n \frac{\pi |z_i|^2}{a_i} < 1\,\right\}
\end{equation*}
denote the open ellipsoid whose projection to the $j$-th complex coordinate plane is $D(a_j)$ and
let $B^{2n}(a)=E(a,\cdots,a)$  be the ball of radius $\sqrt{a/\pi}$. We put
$\mathbb T^4 (A):=\mathbb T^2 (A)\times \mathbb T^2 (A)$, where $\mathbb T^2 (A)$ is the torus
$\BR^2/(A\BZ\oplus \BZ)$ endowed with the symplectic form $dx\wedge dy$ inherited from $\BR^2$.

\vskip 2mm
We define the symplectic capacities
\begin{equation*}
  c_{\rm EZ}(a):=\inf \left\{ A\,|\ E(1,a)\stackrel{s}{\hookrightarrow} Z^4(A)\,\right\},\qquad a\geq 1,
\end{equation*}
\begin{equation*}
  c_{\rm EC}(a):=\inf \left\{ A\,|\ E(1,a)\stackrel{s}{\hookrightarrow} C^4(A)\,\right\},\qquad a\geq 1
\end{equation*}
and
\begin{equation*}
  c_k(C^4):=\inf \left\{ A\,\Big|\ \bigsqcup_k B^4 (1)\stackrel{s}{\hookrightarrow} C^4(A)\,\right\},
\end{equation*}
where $\bigsqcup_k B^4 (1)$ denotes any collection of $k$ disjoint balls $B^4(1)$ in $\BR^4$.

\begin{theorem}\!\!{\sf (Gromov \cite{Gr1}:~1985)}\label{thm:6.17}
Let $a\geq 1$. Then $E(1,a)\stackrel{s}{\hookrightarrow} Z^4(A)$ if and only if
$A\geq 1$. That is, $c_{\rm EZ}(a)=1.$
\end{theorem}
\begin{proof}
The proof follows from Gromov's Nonsqueezing Theorem (cf.\, p.\,24).
\end{proof}

\begin{theorem}\!\!{\sf (Frenkel and M{\"u}ller \cite{F-M}:~2015)}\label{thm:6.18}
Let $\sigma=1+\sqrt{2}$ be the silver ratio. Then the symplectic capacity
$c_{\rm EC}(a)$ satisfies the following:
\vskip 2mm
(a) On the interval $[1,\sigma^2]$, the function $c_{\rm EC}(a)$ is given by the Pell stairs.
Precisely,
\begin{equation*}
  c_{\rm EC}(a)=\begin{cases}
                  1, & \mbox{if } a\in [1,2], \\
                  \frac{1}{\sqrt{2\alpha_n}}, & \mbox{if } a\in [\alpha_n,\beta_n], \\
               \sqrt{\frac{\alpha_{n+1}}{2}}, & \mbox{if } a\in [\beta_n, \alpha_{n+1}],
                \end{cases}
\end{equation*}
\qquad \ \
for\ all\ $n\geq 0$.
\vskip 2mm
(b) On the interval $[\sigma^2,7\frac{1}{32}]$ we have $c_{\rm EC}(a)=\sqrt{a/2}$ except on seven disjoint\\
\indent \ \ \ \ intervals, where $c_{\rm EC}$ is piecewise linear.
\vskip 2mm
(c) $c_{\rm EC}(a)=\sqrt{a/2}$ for all $a\geq 7\frac{1}{32}$.
\end{theorem}
\noindent
Here, $\alpha_n$ and $\beta_n$ are defined as follows: The {\sf Pell numbers} $P_n$ and the
{\sf half companion Pell numbers} $H_n$ are defined by the following recurrence relations:

\begin{eqnarray*}
  P_0\!\!\! &=&\!\!\! 0,\ \ P_1=1,\ \ P_n=2P_{n-1}+ P_{n-2}, \\
  H_0\!\!\! &=&\!\!\! 1,\ \ H_1=1,\ \ H_n=2H_{n-1}+ H_{n-2}.
\end{eqnarray*}
The two sequences $(\alpha_n)_{n\geq 0}$ and $(\beta_n)_{n\geq 0}$ are then defined
as follows:
$$
\alpha_n:=
\begin{cases}
  \frac{2P_{n+1}^2}{H_n^2} & \mbox{if }\ n\ \rm{is\ even}, \\
  \frac{H_{n+1}^2}{2P_n^2} & \mbox{if}\ n\ \rm{is\ odd;}
\end{cases}
\qquad
\beta_n:=
\begin{cases}
  \frac{H_{n+2}}{H_n} & \mbox{if }\ n\ \rm{is\ even}, \\
  \frac{P_{n+2}}{P_n} & \mbox{if}\ n\ \rm{is\ odd}.
\end{cases}
$$
\noindent
Then, for all $n\geq 0,$
\[
\alpha_0=2 < \beta_0=3 < \alpha_1=\frac{9}{2}< \beta_1=5 < \cdots< \alpha_n
< \beta_n < \alpha_{n+1} < \beta_{n+1} < \cdots,
\]
\noindent
and both sequences  converge to $\sigma^2=3+2\sqrt{2}=5.82842712\cdots, $ which is
the square of the silver ratio $\sigma=1+\sqrt{2}$.

\vskip 2mm
In the proof of Theorem 6.18, $J$-holomorphic curves (or pseudoholomorphic curves) play
an important role to show the existence of certain symplectic embeddings.
It is a long journey from Gromov's Nonsqueezing Theorem to Theorem 6.18 (it took 30 years):
the first step is the solution of the ball packing problem
$$
\coprod_{i=1}^{k}B^4(a_i)\stackrel{s}{\hookrightarrow} C^4(A)
$$
that started with \cite{Mc4+} in 1994, and the second step is the translation of the problem
$$
E(1,a)\stackrel{s}{\hookrightarrow} C^4(A)
$$
to a ball packing problem \cite[\S 6]{Mc2} in 2009.

\begin{remark}\label{rk:6.19}
$\bigsqcup_k B^4(1) \stackrel{s}{\hookrightarrow} C^4(A)$ if and only if
$E(1,k)\stackrel{s}{\hookrightarrow} C^4(A)$, that is, $c_k (C^4)=c_{\rm EC}(k)$
for all $k\in\BZ^+.$ That is, the ball packing problem
$\coprod_k B^4(1)\stackrel{s}{\hookrightarrow} C^4(A)$
is included in the $1$--parametric problem
$E(1,a)\stackrel{s}{\hookrightarrow} C^4(A).$
\end{remark}

\vskip 3mm
For two positive real numbers $a$ and $b$, let
${\mathcal N} (a,b)=\{ N_k(a,b)\}_{k\geq 0}$ be the sequence of numbers formed by arranging
all the linear combinations $ma+nb$ with $m,n\geq 0$ in non-decreasing order with repetitions.
We say that $\mathcal N (a,b)\leq \mathcal N (c,d)$ if, for every $k$, the $k$-th entry of
$\mathcal N (a,b)$ is smaller than or equal to the $k$-th entry of $\mathcal N (c,d)$.
McDuff proved the Hofer conjecture:
\begin{theorem}\!\!{\sf (McDuff \cite{Mc3}:~2011)}\label{thm:6.20}
 $E(a,b)$ embeds symplectically into $E(c,d)$ if and only if
$\mathcal N (a,b)\leq \mathcal N (c,d)$.
\end{theorem}
She proved the Hofer conjecture using the Taubes-Seiberg-Witten theory and pseudo-holomorphic curves technique. We can show that
\begin{equation}\label{(6.1)}
E(1,a)\stackrel{s}{\hookrightarrow} E(A,2A)\ {\rm if\ and\ only\ if}\
E(1,a) \stackrel{s}{\hookrightarrow} C^4(A).
\end{equation}
Combining Theorem \ref{thm:6.20} and \eqref{(6.1)}, we obtain
\begin{equation}\label{(6.2)}
c_{\rm EC}(a)=\sup_{k\in\BZ^+} \frac{N_k(1,a)}{N_k(1,2)}
\end{equation}
\vskip 2mm
Given a real number $a\geq 1$ and a positive real number $\mu>0$, we let
$$
\widetilde{E}(1,a):=\{ (x_i)\in \BR^4\,\vert \ x_1^2+x_2^2 +a^{-1}(x_3^2+x_4^2)\leq 1 \}
\subset \BR^4
$$
and
$$
B^4(\mu):=\{ (x_i)\in \BR^4\,\vert \ x_1^2+x_2^2 +x_3^2+x_4^2\leq \mu\}\subset \BR^4
$$
be the closed ellipsoid and the closed ball in $\BR^4$ respectively.
We define the function $c:[1,\infty)\lrt \BR$ by
\begin{equation}\label{(6.3)}
  c(a):=\inf \{\mu\,\vert\ \widetilde{E}(1,a)\stackrel{s}{\hookrightarrow}B^4(\mu)\},
  \quad a\in [1,\infty).
\end{equation}
Then it is seen that the function $c(a)$ is continuous and non-decreasing. Furthermore,
it has the following scaling property:
\begin{equation}\label{(6.4)}
  \frac{c(\lambda a)}{\lambda a}\leq \frac{c(a)}{a}\quad {\rm when}\ \lambda>1.
\end{equation}
Let $\{ f_n\,\vert\ n\geq 0\}=\{ 0,1,1,2,3,\ldots\}$ be the Fibonacci sequence.
We put $g_0=1$ and $g_n:=f_{2n-1}$ for for each $n\geq 1$.
And for each $n\geq 0$, we define
$$
a_n=\left( \frac{g_{n+1}}{g_n}\right)^2 \quad {\rm and}\quad
b_n=\frac{g_{n+2}}{g_n}.
$$
More generally, we get
$$ \cdots < a_n < b_n < a_{n+1} < b_{n+1} < a_{n+2} < b_{n+2} < \cdots $$
and
$$\lim_{n\rightarrow \infty} a_n =\lim_{n\rightarrow \infty} b_n
=\tau^4={\small{\frac{7+3\sqrt{5}}{2}}}
 =6.854101966\cdots,$$
where $\tau=\frac{1+\sqrt{5}}{2}=1.61803398\cdots$ is the golden ratio.
For $a\geq 1$, we define
\begin{equation}\label{(6.5)}
  c_{\rm{ECH}}(a):=\inf \{ \mu>0\,\vert \ \mathcal{N}(1,a)\leq
  \mathcal{N} (\mu,\mu) \}.
\end{equation}

\vskip 2mm
McDuff and Schlenk proved the following:
\begin{theorem}\!\!{\sf (McDuff and Schlenk \cite{Mc7}: 2012)}\label{thm:6.21}\ \ \
  \rm{(i)} For each $n\geq 0,$
$$ c(a)=\frac{a}{\sqrt{a_n}}\qquad {\rm{for}}\ a\in [a_n,b_n], $$
and $c$ is constant with value $\sqrt{a_{n+1}}$ on the interval $[b_n,a_{n+1}]$.
\vskip 2mm\noindent
\rm{(ii)} $c(a)=\frac{a+1}{3}$ on $[\tau^4,7].$
\vskip 2mm\noindent
\rm{(iii)} There are a finite number of closed disjoint intervals
$I_j\subset [7,8\frac{1}{36}]$ such that $c(a)=\sqrt{a}$ for all $a>7,\ a\in I_j$.
Moreover, $c$ is piecewise linear in each $I_j$, with one non-smooth point in the
interior of $I_j.$
\vskip 2mm\noindent
\rm{(iv)} $c(a)=\sqrt{a}$ {\rm for}\ $a\geq 8\frac{1}{36}.$

\vskip 2mm\noindent
\rm{(v)} $c_{\rm{ECH}}(a)\geq c(a)$ \ \, for\ all\ $a\geq 1.$
\end{theorem}
\vskip 2mm
I would like to make comments on Theorem 6.21. Due to Gromov's nonsqueezing theorem,
we see that the preservation of volume is the only obstruction for the symplectic
embedding problem. Later on, it was known that other constraints come from the
characteristic flow on the boundary of a domain, for example, the Ekeland--Hofer
capacities. Also there are other obstructions related to the theory of ECH (embedded
contact homology) developed by M.\,L. Hutchings (cf.\,\cite{Hu2,Hu3}).
\vskip 2mm
In their paper, McDuff and Schlenk compute $c(a)$ for all $a\in [1,\infty)$, and
show that the structure of the graph of the function $c(a)$ is very rich.
It is shown by them that the function $c(a)$ is related to Fibonacci numbers,
weight expansions, continued fractions, exceptional curves in the blowups of
$\mathbb{P}^2\BC$, and to the problem of counting lattice points in the
right-angled triangles. Pseudoholomorphic curves, the ECH theory and
the Taubes-Seiberg-Witten theory are used as tools in their proof.

\begin{theorem}\!\!{\sf (Latschev, McDuff and Schlenk\cite{LMS}:~2013;
Entov and Verbitsky \cite{E-V}:~2016)}
\label{thm:6.22}
Let $a\geq 1$. Then
  $E(1,a)\stackrel{s}{\hookrightarrow} \mathbb T^4 (A)$ whenever
  ${\rm Vol}(E(1,a))< {\rm Vol}(\mathbb T^4 (A))$.
\end{theorem}

\begin{theorem}\!\!{\sf (Schlenk \cite[p.\,154]{Sc}:~2018)}\label{thm:6.23}
  If $B^{2n}(a) \bigsqcup B^{2n}(a)\stackrel{s}{\hookrightarrow} B^{2n}(A),$ then $2a\leq A.$
\end{theorem}

Let $(M,\omega)$ be a symplectic manifold of dimension $2n$. The {\sf Gromov width} is defined to be
\begin{equation*}
  c_B^n(M,\omega):=\sup \left\{ a\,|\ B^{2n}(a)\stackrel{s}{\hookrightarrow} (M,\omega)\ \right\}.
\end{equation*}

\vskip 2mm\noindent
\begin{problem}\label{prob:6.24}
Compute $ c_B^n(M,\omega)$ for $n\geq 2.$ Is it finite ?
\end{problem}
\vskip 2mm
We recall the following principle in Section 4.
\vskip 2mm\noindent
{\bf Eliashberg's\ Principle}\,{\sf \cite{El}:} {\it An obstruction to symplectic embeddings
(beyond the volume condition) can be described by a J-holomorphic curve.}

\vskip 2mm
Surprisingly J-holomorphic curves can be used to construct symplectic embeddings.

\begin{theorem}\!\!{\sf (Schlenk \cite[p.\,177]{Sc}:~2018)}\label{thm:6.25}
(1) There exists a symplectic embedding
$$
P(1,\infty,\infty) \stackrel{s}{\hookrightarrow} P(2,2,\infty).
$$
\vskip 2mm
(2) There exists a symplectic embedding $P(1,a,a) \stackrel{s}{\hookrightarrow} P(2,2,\infty)$ for all
$a\geq 1.$
\vskip 2mm
(3) There exists a symplectic embedding $P(1,\infty,\infty) \stackrel{s}{\hookrightarrow} P(2+\varepsilon,2+\varepsilon,\infty)$ \\
\indent \ \ \ \ \ for all $\varepsilon>0.$
\end{theorem}

\begin{theorem}\!\!{\sf (Schlenk \cite[p.\,178]{Sc}:~2018)}\label{thm:6.26}
For any $n\geq 3$ and for every $\varepsilon>0,$ there exists a symplectic embedding
$F:Z^{2n}(1) \stackrel{s}{\hookrightarrow} \BC^n$ such that
\begin{equation*}
  {\rm vol}_{2k} \left( \pi_k (F(Z^{2n}(1)))\right) < \varepsilon
\end{equation*}
for $k=2,3,\cdots,n-1,$ where ${\rm vol}_{2k}(U):=(k!)^{-1}\int_U \omega_0^k$ denotes the Euclidean volume of a domain $U$ in $\BC^k$ and $\pi_k:\BC^n\lrt \BC^k$ is the projection given by
$(z_1,\cdots,z_n)\mapsto (z_1,\cdots,z_k)$.
\end{theorem}
\vskip 2mm
For a positive integer $k$, we consider the function defined by
\begin{equation}\label{(6.6)}
  c_k(x):=\inf \{ A\,|\ E(1,x)\times \BR^{2k}\stackrel{s}{\hookrightarrow} B^4 (A)\times \BR^{2k}\,\}.
\end{equation}
\vskip 2mm
R.\,K. Hind \cite{H} proved that $c_k(x)\leq \frac{3x}{x+1}$ if $x> \tau^4$, where
$\tau=\frac{1+\sqrt{5}}{2}$ is the golden ratio. This imples that $c_k(x)< c_0 (x)$ for $k\geq 1.$
In 2018, D. McDuff showed that Hind's bound is sharp for certain values of $x$.
%%%%%%%%%%%%%%%%%%%%%%%%%%%%%%%%%%%%%%%%%%%%%%%%%%%%%%%%%%%%%%%%%%%%%%%%%%%%%%%%%%%%%%%%%%%%%%%%%%%
\begin{theorem}\!\!{\sf (McDuff \cite{Mc4}:~2018)}\label{thm:6.27}
  If $x=3m-1$ with $m\in \BZ^+$ and $x> \tau^4$, then $c_k(x)= \frac{3x}{x+1}.$
\end{theorem}
\vskip 2mm
McDuff proposed the following conjecture.
\begin{conjecture}\label{conj:6.28}
$c_k(x)= \frac{3x}{x+1}$ \ if\ $x> \tau^4$.
\end{conjecture}
\vskip 2mm
Let $K$ be a compact connected Lie group and let $\mathfrak k^*$ be the dual of the Lie algebra
$\mathfrak k$ of $K$. Each coadjoint orbit ${\mathcal O}$ in $\mathfrak k$ is equipped with the Kostant-Kirillov-Souriau symplectic form $\omega$ canonically defined by
\begin{equation*}
  \omega_\eta (X^{\sharp}, Y^{\sharp})=\langle \eta,[X,Y]\rangle,\quad \eta\in {\mathcal O}_\lambda, \
  X,Y\in \mathfrak k,
\end{equation*}
where $X^{\sharp}, Y^{\sharp}$ are the vector fields on $\mathfrak k^*$ corresponding to
$X,Y\in \mathfrak k$ induced by the coadjoint action. Each coadjoint orbit intersects a positive Weyl chamber in a single point. So there is a bijection between the coadjoint orbits and points in the positive Weyl chamber. Points in the interior of the positive Weyl chamber are called {\it regular} points. The orbits corresponding to regular points are called {\sf generic orbits} that are diffeomorphic to $K/T$ for $T$ a maximal torus of $K$. Coadjoint orbits intersecting the positive Weyl chamber at its boundary are called {\sf degenerate orbits}.

\vskip 2mm
Caviedes Castro \cite{CC} proved the following:
\begin{theorem}\!\!{\sf (Castro \cite{CC}:~2016)}\label{thm:6.29}
  Let $K$ be a compact connected Lie group. The Gromov width of a coadjoint orbit ${\mathbb O}_\lambda$ through a point $\lambda$ lying on some rational line in $\mathfrak t^*,$ equipped with the Kostant-Kirillov-Souriau symplectic form $\omega_\lambda$, can not be greater than the following quantity
\begin{equation}\label{(6.7)}
\min \left\{ \,|\langle \lambda,\alpha^{\vee} \rangle|\,|\ \alpha^{\vee}\ {\rm is\ a \ coroot\ and}\
\langle \lambda,\alpha^{\vee} \rangle \neq 0\ \right\}.
\end{equation}
Here $\mathfrak t^*$ is the dual of the Lie algebra of a maximal torus $T$ of $K$.
\end{theorem}
Fang, Littelmann and Pabiniak \cite{FLP} gave a uniform proof for the conjectured Gromov width of rational coadjoint orbits of all compact connected simple Lie groups by analyzing simplices in Newton-Okounkov bodies.
\begin{theorem}\!\!{\sf (Fang, Littelmann and Pabiniak \cite{FLP}:~2018)}\label{thm:6.30}
  Let $K$ be a compact connected Lie group. Then the Gromov width of a coadjoint orbit ${\mathbb O}_\lambda$ through a point $\lambda$ lying on some rational line in $\mathfrak t^*,$ equipped with the Kostant-Kirillov-Souriau symplectic form $\omega_\lambda$, is at least
\begin{equation*}
\min \left\{ \,|\langle \lambda,\alpha^{\vee} \rangle|\,|\ \alpha^{\vee}\ {\rm is\ a \ coroot\ and}\
\langle \lambda,\alpha^{\vee} \rangle \neq 0\ \right\}.
\end{equation*}
Here $\mathfrak t^*$ is the dual of the Lie algebra of a maximal torus $T$ of $K$.
\end{theorem}
Furthermore they proved the following fact in \cite{FLP}.
\begin{proposition}\label{prop:6.31}
Let $K$ be a compact connected Lie group, not of type $G_2,\,F_4$ of $E_8$ and let
$({\mathbb O}_\lambda,\omega_\lambda)$ be its generic coadjoint orbit through a point $\lambda$ lying on some rational line in $\mathfrak t^*,$ equipped with the Kostant-Kirillov-Souriau symplectic form $\omega_\lambda$. Then there exists a symplectic embedding of a ball of capacity (6.7) into
$({\mathbb O}_\lambda,\omega_\lambda)$.
\end{proposition}

\vskip 2mm
We would like to mention the recent work on infinite staircases obtained by
Cristofaro-Gardiner, Holm, Mandini and Pires. Let $(X,\omega)$ be a symplectic
four-manifold. We define the {\sf ellipsoid embedding function} of $X$ by
\begin{equation}\label{(6.8)}
  c_X (a):=\inf \left\{ \lambda\,\big|\ E(1,a)\stackrel{s}{\hookrightarrow} \lambda X \right\},
  \quad {\rm for}\ a\geq 1,
\end{equation}
where $\lambda X$ represents the symplectic scaling $(X,\lambda\cdot\omega)$ of $(X,\omega)$.
It is seen that $c_X$ is a continuous nondecreasing function on $[1,\infty)$
(cf.\,\cite[Proposition 2.1]{Cr-2}).

\begin{definition}\label{def:6.32}
We say that the ellipsoid embedding function $c_X (a)$ has an {\sf infinite staircase} if
its graph has infinitely many non-smooth points.
\end{definition}

\vskip 2mm
%\noindent
A 4-dimensional {\sf toric domain} $X$ is defined to be the preimage of a domain
$\Delta\subset \BR_{\geq 0}^2$ under the moment map
\begin{equation*}
  \mu:\BC^2\lrt \BR^2,\qquad (z_1,z_2)\mapsto (\pi \vert z_1\vert^2,\pi \vert z_2\vert^2).
\end{equation*}
\noindent
Here, $\BR_{\geq 0}^2=\{ (x_1,x_2)\in \BR^2\,\vert \ x_1 \geq 0,\ x_2 \geq 0 \}.$
A {\sf convex toric domain} is the preimage under a closed region $\Delta\subset \BR_{\geq 0}^2$
that is convex, connected, and contains the origin in its interior (in the subspace topology
on $\BR_{\geq 0}^2 \subset \BR^2).$ We denote this $X_\Delta=\mu^{-1}(\Delta)$ and call
the region $\Delta$ the {\sf moment polygon} of $X_\Delta$. We note that $X_\Delta$ is
a manifold with boundary.
\vskip 2mm
Cristofaro-Gardiner et al. proved the following.
\begin{theorem}\!\!{\sf (Cristofaro-Gardiner, Holm, Mandini and Pires \cite{Cr-2}:~2025)}
\label{thm:6.33}
\par\noindent
 Let $\Delta \subset \BR_{\geq 0}^2$ be a convex region that is also a Delzant polygon for a closed
 toric symplectic four-manifold $M$. Then there exists a symplectic embedding
 \begin{equation}\label{(6.9)}
  E(a,b)\stackrel{s}{\hookrightarrow}  M
 \end{equation}
 if and only if there exists a symplectic embedding
 \begin{equation}\label{(6.10)}
  E(a,b)\stackrel{s}{\hookrightarrow}  X_\Delta.
 \end{equation}
\end{theorem}
\vskip 2mm
Following \cite[\S 2.2]{Cr-1}, to a convex toric domain $X_\Delta$ we can associate
a {\sf negative weight expansion} $(b;b_1,b_2,\ldots).$ We note that two different
convex toric domains can have the same negative weight expansion.

\begin{definition}\label{(6.34)}
We say that a negative weight expansion
$$(b;b_1,b_2,\ldots)$$
\noindent
is {\sf finite} if there are finitely many non-zero $b_i'$s, and we say that such
a convex tori domain $X_\Delta$ has {\sf finite type}. Given a finite type
convex toric domain $X$ with negative weight expansion $(b;b_1,\ldots,b_n),$
we define:
\begin{equation*}
  {\rm per}=3b-\sum_{i=1}^n b_i\qquad {\rm and}\qquad
  {\rm vol}=b^2-\sum_{i=1}^n b_i^2.
\end{equation*}
We say that $X_\Delta$ is a {\sf rational convex toric domain} if
$b,b_1,\ldots,$ and $b_n$ are all rational numbers.
\end{definition}

They also proved the following:
\begin{theorem}\!\!{\sf (Cristofaro-Gardiner, Holm, Mandini and Pires \cite{Cr-2}:~2025)}
\label{thm:6.35}
\par\noindent
 Let $X$ be finite type convex toric domain. If the ellipsoid embedding function $c_X(a)$ has
 an infinite staircase then it accumulate at $a_0$, a real solution of the quadratic equation
 \begin{equation}\label{(6.11)}
 x^2-\left( \frac{{\rm per}^2}{{\rm vol}}-2 \right)x+1=0.
 \end{equation}
 Furthermore, $a_0$, the ellipsoid embedding function touches the volume curve:
 $$ c_X(a_0)=\sqrt{\frac{a_0}{\rm vol}}.$$
 Here, $c_X(a)=\sqrt{\frac{a}{\rm vol}}\ (a\geq 1)$  is called the volume curve.
\end{theorem}
\vskip 2mm
They proposed the following conjecture:
\begin{conjecture}\label{conj:6.36}
If the ellipsoid embedding function of a {\sf rational} convex toric domain has
an infinite staircase, then its moment polygon is a scaling of a {\it reflexive} polygon.
\end{conjecture}
\vskip 2mm\noindent
Here, a convex lattice polygon is called {\sf reflexive} if it has exactly one
interior lattice point.
%%%%%%%%%%%%%%%%%%%%%%%%%%%%%%%%%%%%%%%%%%%%%%%%%%%%%%%%%%%%%%%%%%%%%%%%%%%%%%%%%%%%%%%%%%%%%%%%%%%%%%%%
%%%%%%%%%%%%%%%%%%%%%%%%%%%%%%%%%%%%%%%%%%%%%%%%%%%%%%%%%%%%%%%%%%%%%%%%%%%%%%%%%%%%%%%%%%%%%%%%%%%%%%%%
%%%%%%%%%%%%%%%%%%%%%%%%%%%%%%%%%%%%%%%%%%%%%%%%%%%%%%%%%%%%%%%%%%%%%%%%%%%%%%%%%%%%%%%%%%%%%%%%%%%%%%%%
%%%%%%%%%%%%%%%%%%%%%%%%%%%%%%%%%%%%%%%%%%%%%%%%%%%%%%%%%%%%%%%%%%%%%%%%%%%%%%%%%%%%%%%%%%%%%%%%%%%%%%%%
%%%%%%%%%%%%%%%%%%%%%%%%%%%%%%%%%%%%%%%%%%%%%%%%%%%%%%%%%%%%%%%%%%%%%%%%%%%%%%%%%%%%%%%%%%%%%%%%%%%%%%%%
%%%%%%%%%%%%%%%%%%%%%%%%%%%%%%%%%%%%%%%%%%%%%%%%%%%%%%%%%%%%%%%%%%%%%%%%%%%%%%%%%%%%%%%%%%%%%%%%%%%%%%%%

\vskip 5mm
\subsection{The theory of the Gromov-Witten invariants}
\leavevmode\par\vspace{2mm}
%\noindent
Let $(M,\omega,J)$ be a symplectic manifold of dimension $2n$ with a $\omega$-compatible almost complex structure $J$. For two nonnegative integers $g,k\geq 0$, let $\overline{\mathcal M}_{g,k}$ be the Deligne-Mumford moduli space of stable curves of genus $g$ with $k$ marked points, and
$\overline{\mathcal M}_{g,k,A}$ be the moduli space of stable maps into $M$ of homology class
$A\in H_2(M,\BZ)/{\rm torsion}$. The elements of $\overline{\mathcal M}_{g,k,A}$ are of the form
\begin{equation*}
  (\Sigma,p_1,\cdots,p_k,f),
\end{equation*}
where $\Sigma$ is a (not necessarily stable) curve of genus $g$ with $k$ marked points and
$f:\Sigma\lrt M$ is a pseudoholomorphic curve. Let
\begin{equation*}
Y:=\overline{\mathcal M}_{g,k} \times M^k.
\end{equation*}
Then we have the evaluation map ${\sf ev}:\overline{\mathcal M}_{g,k,A}\lrt Y$ defined by
\begin{equation*}
  {\sf ev}(\Sigma,p_1,\cdots,p_k,f)=({\rm st}(C,p_1,\cdots,p_k),f(p_1),\cdots,f(p_k)),
\end{equation*}
where ${\rm st}(C,p_1,\cdots,p_k)$ is the stabilization of $C$. By the Atiyah-Singer index theorem, we obtain
\begin{equation*}
d:=\dim_{\BR} \overline{\mathcal M}_{g,k,A}=2(n-3)(1-g)+2k+2\,c_1(M)\cdot A.
\end{equation*}
The evaluation map send the fundamental class of $\overline{\mathcal M}_{g,k,A}$ to a $d$-dimensional rational homology class $\Phi_{g,k}^A \in H_d(Y,\BQ).$ The homology class $\Phi_{g,k}^A \in H_d(Y,\BQ)$
is called the {\sf Gromov-Witten\ invariant} of $M$ for the data $(g,k,A).$ It is an invariant of the symplectic isotopy class of $(M,\omega)$.
\vskip 2mm
Let us interpret $\Phi_{g,k}^A$ geometrically. If $\beta\in H_* (\overline{\mathcal M}_{g,k})$ and
$\alpha_1,\cdots,\alpha_k\in H_* (M)$ such that the sum of the codimensions of $\beta,\alpha_1,\cdots,\alpha_k$ is equal to $d$, we define
\begin{equation}\label{(6.12)}
\Phi_{g,k}^A (\beta,\alpha_1,\cdots,\alpha_k):=\Phi_{g,k}^A \cdot \beta\cdot\alpha_1\cdots\alpha_k \in H_0(Y,\BQ),
\end{equation}
where $\cdot$ denotes the intersection product in $H_*(Y,\BQ).$ This is a rational number.
Let $e_i:\overline{\mathcal M}_{g,k,A}\lrt M\ (1\leq i\leq k)$ be the evaluation map at the $i$-th marked point. If $\delta_1,\cdots,\delta_k\in H^* (M,\BZ)/{\rm torsion}$ such that
$\sum_{i=1}^k \deg (\delta_i)=d$, we put
\begin{equation}\label{(6.13)}
\Phi_{g,k}^A (\delta_1,\cdots,\delta_k):=\int_{\overline{\mathcal M}_{g,k,A}} e_1^*\delta_1\wedge\cdots\wedge e_k^*\delta_k.
\end{equation}
The number $\Phi_{g,k}^A (\delta_1,\cdots,\delta_k)$ can roughly be understood as the number of pseudoholomorphic curves of genus $g$ representing the homology class $A$ and intersecting $k$ given cycles ${\rm PD}(\delta_i)$ Poincar{\'e} dual to the cohomolgy classes $\delta_i\ (1\leq i\leq k).$
\vskip 2mm
E. Witten \cite{Wi} defined the so-called ${\sf Gromov\!-\!Witten\ potential}\ \Phi_{\omega}^M:H^*(M,\BC)\lrt \BC$ by
\begin{equation}\label{(6.14)}
 \Phi_{\omega}^M (\delta_0,\delta_1,\cdots,\delta_{2n}):=\sum_k\sum_A\sum_{i_1,\cdots,i_k}
 \frac{\exp (-\int_A \omega)}{k!} \Phi_{g,k}^A (\delta_{i_1},\cdots,\delta_{i_k}),
\end{equation}
where $\delta_i\in H^i (M,\BC)\ (0\leq i\leq 2n),\ A$ runs over $H_2 (M,\BZ)$ and unordered
$i_1,\cdots,i_k \in \{ 0,1,\cdots, 2n \}$ with $\sum_{\nu=1}^{k}\deg (\delta_{i_{\nu}})=d.$ The convergence problem arises. The fact that the Gromov-Witten potential $\Phi_{\omega}^M$ satisfies the WDVV equations
was proved by Y. Ruan and G. Tian \cite{R-T}. Maxim Kontsevich and Yuri Manin \cite{K-M} described how the WDVV equations yields a potential Dubrovin structure on $H^*(M,\BC)$, understood as a supermanifold. Therefore each tangent space of $H^*(M,\BC)$ is endowed with a metric given by Poincar{\'e} duality, and a multiplication
\begin{equation}\label{(6.15)}
x * y=\sum_{i,j,k}A_{ij}^k x_i y_j \phi_k,\qquad A_{ij}^k=\sum_{\ell}\partial_i\partial_j\partial_{\ell}\Phi_\omega^M g^{k\ell},
\end{equation}
where the set $\{ \phi_i \}$ is a homogeneous basis of $H^*(X,\BC)$ with $\phi_i\cdot\phi_j=g_{ij}$, the
$(g^{ij})$ denotes the inverse matrix $(g_{ij})$, $x=\sum_i x_i\phi_i$ and $y=\sum_i y_i\phi_i$. This is the quantum deformation of the cup product and the WDVV equations are equivalent to the associativity of the multiplication. The ordinary cup product is the limit of $\Phi_{t\omega}^M$ with $t\rightarrow \infty.$

\vskip 2mm
Y. Ruan and G. Tian\,\cite{R-T} nicely developed a theory of quantum cohomology which is related to symplectic topology, algebraic geometry, quantum field theory, mirror symmetry, and differential topology in 4-dimensional manifolds. They first constructed the Gromov-Witten theory for semi-Fano symplectic manifolds. Thereafter the Gromov-Witten theory was generalized by other people.
The quantum cohomology ring endows the affine space $H^*(X,\BC)$ with the structure of a Frobenius manifold, namely, a Riemannian manifold with an associative product on the tangent spaces and various compatibilities. C. Taubes\,\cite{Ta1,Ta2}  proved that in $4$-dimensional symplectic manifolds, certain Gromov-Witten invariants coincide with the gauge-theoretic Seiberg-Witten invariants.
This relates symplectic topology and differential topology via the Gromov-Witten theory.
W. Chen and Y. Ruan\,\cite{CR1, CR2} defined the Gromov-Witten invariants for compact symplectic {\sf orbifolds} extending the Gromov-Witten invariants for compact symplectic manifolds. Based on a proposal by E. Witten,
H. Fan, T. Jarvis and Y. Ruan\,\cite{FJR} introduced and developed the new Gromov-Witten type theory of geometric invariants, known as the {\sf FJRW\ theory} that is the mathematically rigorous development of topological gravity coupled with $A$-type topological Landau-Ginzburg matter, as an intersection theory on the moduli space of solutions of the Witten equation. The FJRW theory is believed to be the counterpart of the Gromov-Witten theory in the Landau-Ginzburg model \cite{C-R1, C-R2, K-S, Sh1}. The relationship between these two theories is referred to as Landau-Ginzburg/Calabi-Yau (briefly LG/CY) correspondence.
\vskip 2mm
Recently Y. Shen and J. Zhou\,\cite{Sh3} proved the LG/CY correspondence between the Gromov-Witten theories of
elliptic orbifold curves $\mathbb P_{3,3,3},\,\mathbb P_{4,4,2},\,\mathbb P_{6,3,2}$ and
$\mathbb P_{2,2,2,2}$ and their FJRW theory counterparts via the theory of quasi-modular forms. We briefly describe this correspondence. Let $W:\BC^3\lrt \BC$ be the weighted homogeneous polynomial with weights $q_1,q_2,q_3$, the so-called {\sf superpotential} of the LG-model that satisfies the {\sf Calabi-Yau condition}
$q_1+q_2+q_3=1$. Let
\begin{equation}\label{(6.16)}
  G_W:=\{ (\lambda_1,\lambda_2,\lambda_3)\in (\BC^*)^3\,|\ W(\lambda_1x_1,\lambda_2x_2,\lambda_3x_3)=W(x_1,x_2,x_3)\,\}
\end{equation}
be the {\sf group of diagonal symmetries} and let $G$ be a subgroup of $G_W$ containing the exponential grading element
\begin{equation}\label{(6.17)}
 J_W:=(\exp(2\pi i\,q_1),\exp(2\pi i\,q_2),\exp(2\pi i\,q_3) ),\qquad i=\sqrt{-1}.
\end{equation}
The hypersurface $X_W$ defined by $\{ W=0\}$ is a one-dimensional Calabi-Yau variety in a weighted projective space. Then $G_W$ acts on $X_W$, and $J_W$ acts trivially. Thus we obtain the CY orbifold curve which is a global quotient
\begin{equation}\label{(6.18)}
\frak X_W:=X_W/(G/\langle J_W\rangle).
\end{equation}
The elliptic orbifold curves $\frak X_W$ are as follows:
\begin{eqnarray}\label{(6.19)--(6.22)}
   W &=& x_1^3 +x_2^3+x_3^3,\quad G=G_W,\qquad \frak X_W=\mathbb P_{3,3,3};\\
   W &=& x_1^4 +x_2^4+x_3^2,\quad G=G_W,\qquad\frak X_W=\mathbb P_{4,4,2}; \\
   W &=& x_1^6 +x_2^3+x_3^2,\quad G=G_W,\qquad\frak X_W=\mathbb P_{6,3,2}; \\
   W &=& x_1^4 +x_2^4+x_3^2,\quad G=G_1\times G_{x_3^2},\qquad\frak X_W=\mathbb P_{2,2,2,2},
\end{eqnarray}
where $G_1:=\langle (i,i),(1,-1)\rangle.$
\vskip 2mm
For the pair $(W,G)$, both GW theory and FJRW theory come with a graded vector space equipped with a non-degenerate pairing, which we denote by
\begin{equation*}
   \left( \mathscr H^{\rm GW},\eta^{\rm GW}\right),\qquad
\left( \mathscr H^{\rm FJRW},\eta^{\rm FJRW}\right).
\end{equation*}
Here $\mathscr H^{\rm GW}$ is the {\sf Chen-Ruan cohomology} \cite{CR1, CR2} of $\frak X_W$, and
$\mathscr H^{\rm FJRW}$ is the {\sf FJRW state space} \cite{FJR} of $(W,G)$.
Let $\overline{\mathcal M}_{g,k}$ be the Deligne-Mumford moduli space of $k$-pointed stable curves of genus $g$ and $\psi_j\in H^2 (\overline{\mathcal M}_{g,k},\BQ)$ be the $j$-th $\psi$-class. Let $\beta$ be an effective curve class in the underlying coarse moduli of $\frak X_W, \{ \alpha_j\}$ be elements in
$\mathscr H^{\rm GW}$ and $\{ \gamma_j \}$ be elements in $\mathscr H^{\rm FJRW}$. Then one can define the
{\sf ancester GW invariant}
$\langle \alpha_1 \psi_1^{\ell_1},\cdots,\alpha_k \psi_1^{\ell_k}\rangle_{g,k,\beta}^{\rm GW}$
(cf.\,\cite[p.\,6\,(2.2)]{Sh3}) and the
{\sf FJRW invariant}
$\langle \gamma_1 \psi_1^{\ell_1},\cdots,\gamma_k \psi_1^{\ell_k}\rangle_{g,k}^{\rm FJRW}$
(cf.\,\cite[p.\,9\ (2.13)]{Sh3}).
\vskip 2mm
We parametrize a K{\"a}hler class ${\mathcal P}\in \mathscr H^{\rm GW}$ by $t$ and set $q=e^t$. The Divisor Axiom in GW theory allows us define a {\sf GW correlation function} as a formal $q$-series
\begin{equation}\label{(6.23)}
  \langle\!\langle \alpha_1 \psi_1^{\ell_1},\cdots,\alpha_k \psi_1^{\ell_k}
  \rangle\!\rangle_{g,k}^{\rm GW} (q):=\sum_{\beta}\langle \alpha_1 \psi_1^{\ell_1},\cdots,\alpha_k \psi_1^{\ell_k}\rangle_{g,k,\beta}^{\rm GW}\,\, q^{\int_\beta \mathcal P}.
\end{equation}
The GW invariants give rise to various structures on $\mathscr H^{\rm GW}$. Among them the
{\sf quantum multiplication} $\bigstar_q$ is defined by
\begin{equation}\label{(6.24)}
  \alpha_1  \bigstar_q \alpha_2 := \sum_{\mu,\nu} \langle\!\langle \alpha_1,\alpha_2,\mu \rangle\!\rangle_{0,3}^{\rm GW}\,\eta_{\rm GW}^{(\mu,\nu)}\nu.
\end{equation}
Here both $\mu,\nu$ run over a basis of $\mathscr H^{\rm GW}$ and $\eta_{\rm GW}^{(\cdot,\cdot)}$ is the inverse of the pairing $\eta^{\rm GW}(\cdot,\cdot)$. At the large volume limit $t=-\infty$, the quantum multiplication
$\bigstar_q$ becomes the Chen-Ruan product.
\vskip 2mm
Similarly we parametrize a degree $2$ element $\phi\in \mathscr H^{\rm FJRW}$ by $u$ and define an
{\sf FJRW correlation function}
\begin{equation}\label{(6.25)}
  \langle\!\langle \gamma_1 \psi_1^{\ell_1},\cdots,\gamma_k \psi_1^{\ell_k}
  \rangle\!\rangle_{g,k}^{\rm FJRW} (u):=\sum_{n\geq 0} \frac{u^n}{n!}\,
  \langle \gamma_1 \psi_1^{\ell_1},\cdots,\gamma_k \psi_1^{\ell_k},\phi,\cdots,
  \phi\rangle_{g,k+n}^{\rm FJRW}.
\end{equation}
We have a Frobenius algebra $(\mathscr H^{\rm FJRW},\bullet)$, where the multiplication $\bullet$ is defined from the pairing $\eta^{\rm FJRW}$ on $\mathscr H^{\rm FJRW}$ and the genus zero $3$-point invariants through the following formula
\begin{equation}\label{(6.26)}
 \eta^{\rm FJRW}(\gamma_1, \gamma_2 \bullet \gamma_3)=\langle \gamma_1,\gamma_2,\gamma_3
 \rangle_{0,3}^{\rm FJRW}.
\end{equation}
We define the {\sf quantum multiplication} $\bullet_u$ by
\begin{equation}\label{(6.27)}
 \gamma_1\bullet_u \gamma_2:=\sum_{\gamma,\zeta}
 \langle\!\langle \gamma_1,\gamma_2,\gamma
  \rangle\!\rangle_{0,3}^{\rm FJRW}\,\eta_{\rm FJRW}^{(\gamma,\zeta)}\,\zeta.
\end{equation}
Here both $\gamma,\zeta$ run over a basis of $\mathscr H^{\rm FJRW}$ and $\eta_{\rm FJRW}^{(\cdot,\cdot)}$ is the inverse of the pairing $\eta^{\rm FJRW}(\cdot,\cdot)$. The quantum multiplication $\bullet_u$ is a
deformation of the multiplication $\bullet$ in Formula \eqref{(6.26)} as $\bullet_{u=0}=\bullet$.
\vskip 2mm
We refer to \cite{Sh1, Sh2, Sh3} for more precise details on
$\left( \mathscr H^{\rm GW},\eta^{\rm GW}\right)$ and
$\left( \mathscr H^{\rm FJRW},\eta^{\rm FJRW}\right)$.
The LG/CY correspondence \cite{FJR, Wi} says that the two enumerative theories should be equivalent under an appropriate transformation.
\vskip 2mm
Let $\BH$ (resp. $\mathbb D$) be the Poincar{\'e} upper half plane (resp. the Poincar{\'e} unit disk).
We denote by $\widehat{M}(\G)$ (resp. $\widetilde{M}(\G)$) the ring of almost-holomorphic modular forms (resp. quasi-modular forms) on $\BH$ for an arithmetic subgroup $\G$ of $SL(2,\BZ)$. Let $C_{\BH}^\omega$
(resp. $C_{\BD}^\omega$) be the ring of real analytic functions on $\BH$ (resp. $\BD$). We denote by
$\mathcal O_{\BD}$ the ring of holomorphic functions on $\BD$. Shen and Zhou\,\cite{Sh3} introduced the
{\sf Cayley transformation} $\mathfrak{C}:\widehat{M}(\G)\subset C_{\BH}^\omega \lrt C_{\BD}^\omega$ and its variant $\mathfrak{C}_{\rm hol}:\widetilde{M}(\G)\lrt \mathcal O_{\BD}$. They are induced by the Cayley transform $T:\BH\lrt \BD$ based at a point $\tau_0\in\BH$ given by
\begin{equation*}
T(\tau):=\frac{\tau-\tau_0}{ \frac{\tau}{\tau_0-\overline{\tau}_0}-
\frac{\overline{\tau}_0}{\tau_0-\overline{\tau}_0} },\quad \tau\in\BH.
\end{equation*}
Indeed, the Cayley transform $T:\mathbb{H} \lrt \mathbb{D}$ is a geometric variant of
the classical Cayley transform.
Shen and Zhou showed that the GW correlation functions
of elliptic orbifold curves are elliptic $q$--expansions of some quasi-modular forms
near the cusp of the Poincar{\'e} upper half plane $\mathbb{H}$, and the FJRW correlation functions
are represented as a Taylor series ($s$-expansions: $s\in \mathbb{D}$) in the unit disk $\mathbb{D}.$
Using the Cayley transform $T:\BH\lrt \BD$ and the Cayley transformation
$\mathfrak{C}_{\rm hol}:\widehat{M}(\G)\lrt C_{\BD}^\omega$
the GW correlation functions can be identified with
the FJRW correlation functions. Therefore the LG/CY correspondence is established.
Indeed, the Cayley transform is a powerful mathematical bridge linking A-model Gromov-Witten (GW) invariants on the Calabi-Yau side to Fan-Jarvis-Ruan-Witten (FJRW) invariants on the Landau-Ginzburg side.

\begin{theorem}{\sf (Shen and Zhou \cite{Sh3}:\,2018)}\label{thm:6.37}
Let $(W,G)$ be a pair in Formulas (6.19)--(6.22). Then there exists a degree and pairing
isomorphism between the graded vector spaces
\begin{equation*}
{\mathscr G}:\left( \mathscr H^{\rm GW},\eta^{\rm GW}\right) \lrt
\left( \mathscr H^{\rm FJRW},\eta^{\rm FJRW}\right)
\end{equation*}
and the Cayley transform $\mathfrak{C}_{\rm hol}$, based at an elliptic point $\tau_0\in\BH$, such that for any $\{ \alpha_j \}\subseteq \mathscr H^{\rm GW}$,
\begin{equation*}
\mathfrak{C}_{\rm hol}\left( \langle\langle  \alpha_1\psi_1^{\ell_1},\cdots, \alpha_k\psi_k^{\ell_k}             \rangle\rangle^{\rm GW}_{g,k}(q) \right)=
\langle\langle  {\mathscr G}(\alpha_1)\psi_1^{\ell_1},\cdots, {\mathscr G}(\alpha_k)\psi_k^{\ell_k}             \rangle\rangle^{\rm FJRW}_{g,k}(u).
\end{equation*}
Here $\psi_j\in H^2(\overline{\mathcal M}_{g,k},\BQ) \,(1\leq j\leq k)$ is the $j$-th $\psi$-class and $q=e^t$ where $t$ was parametrized by a K{\"a}hler class $\mathcal P\in \mathscr H^{\rm GW}.$
\end{theorem}
\vskip 1mm
The above theorem provides an insight into a number-theoretic understanding of the LG/CY correspondence
by means of the theory of quasi-modular forms. It may be said that it sheds light on understanding
physical phenomena from the viewpoint of number theory.

\vskip 10mm
\subsection{Floer homology, the Fukaya category and mirror symmetry}
\leavevmode\par\vspace{2mm}
%\noindent
In this subsection, we give a very sketchy description of Lagrangian Floer homology, the Fukaya category and homological mirror symmetry (briefly, {\bf HMS}). The Fukaya category is constructed on the basis of Lagrangian Floer homology, and plays an important role in the study of symplectic geometry and the HMS
conjecture. This subsection was written on the base of the references
\cite{Aur,Ci, ER, F1, F1+, F2, F3, Fu, F-O, FO1, FO2, Ko, K-M, Mc3+, Oh1, Oh2, P-Z, Se1, Se2, Yao}

%%%%%%%%%%%%%%%%%%%%%%%%%%%%%%%%%%%%%%%%%%%%%%%%%%%%%%%%%%%%%%%%%%%%%%%%%%%%%%%%%%%%%%%%%%%%%%%%%%%%%%%
%%%%%%%%%%%%%%%%%%%%%%%%%%%%%%%%%%%%%%%%%%%%%%%%%%%%%%%%%%%%%%%%%%%%%%%%%%%%%%%%%%%%%%%%%%%%%%%%%%%%%%%
%%%%%%%%%%%%%%%%%%%%%%%%%%%%%%%%%%%%%%%%%%%%%%%%%%%%%%%%%%%%%%%%%%%%%%%%%%%%%%%%%%%%%%%%%%%%%%%%%%%%%%%
\vskip 3mm\noindent
{\bf (I) Lagrangian Floer homology}
\vskip 2mm
\begin{definition}\label{def:6.38}
A symplectic manifold $(M,\omega)$ is said to be {\sf symplectically aspherical} if
for any smooth map $f:S^2\lrt M$, we have $\int_{S^2}f^*\omega=0.$
\end{definition}

\vskip 2mm
Now we let $(M,\omega)$ be a symplectic manifold. We consider a pair
($\mathfrak{L}_0,\mathfrak{L}_1)$ of compact Lagrangian submanifolds inside $(M,\omega)$
such that $\mathfrak{L}_0$ and $\mathfrak{L}_1$ intersect transversally.
Let
\begin{equation*}
  \wp(\mathfrak{L}_0,\mathfrak{L}_1)=\wp:=\left\{ \gamma:[0,1]\lrt M\,\vert\
  \gamma\ {\rm is\ smooth},\ \gamma(0)\in \mathfrak{L}_0\ {\rm and}\
  \g (1)\in\mathfrak{L}_1 \right\}
\end{equation*}
be the space of smooth paths in $(M,\omega)$. Then $\wp$ has a natural topological structure
given by the $C^\infty$ topology and is not connected in general. To get a nicer space,
we fix any path $\g_0\in \wp$ and consider the universal covering $\widetilde{\wp}$
of its connected component. The universal covering $\widetilde{\wp}$ is given by
$(\g,[u])$, where $\g$ is a path in $\wp$ and $u$ is the homotopy class of maps:
\begin{equation*}
  u:[0,1]\times [0,1] \lrt M, \qquad (s,t)\mapsto u_s(t),
\end{equation*}
for which $u_0=\g_0,\ u_1=\g_1,\ u_s(0)\in \mathfrak{L}_0$ and $u_s (1)\in\mathfrak{L}_1.$
We define the action functional $A:\widetilde{\wp}\lrt \BR$ by
\begin{equation*}
  A((\g,[u])):=-\int_{I^2}u^*\omega,\qquad (\g,[u])\in \widetilde{\wp},\ I:=[0,1].
\end{equation*}
If we assume that $M$ is symplectically aspherical, and $\mathfrak{L}_0,\mathfrak{L}_1$
are simply connected, then the functional $A$ is well defined. The differential map of
$A$ at $(\g,[u])$ is given by
\begin{equation*}
  dA_{(\g,[u])}(X)=-\int_0^1 \omega (\g'(t),X(t))dt
  =-\int_0^1 g_{\g (t)}(J_{\g (t)}(\g'(t)),X(t))dt
\end{equation*}
where $X$ is a vector field along the path $\g$ and $J$ is an almost complex structure
which is compatible with $\omega$, denoting by $g$ the Riemannian metric associated to
$\omega$. Thus, we have the gradient of $A$:
\begin{equation*}
  {\rm grad}\,A([(\g,[u])])=-J\g'.
\end{equation*}
Here, we see $\g'$ as an infinitesimal variation of the point $[(\g,[u])]$ and, therefore,
as an element of the tangent space of $\widetilde{\wp}$.
\vskip 2mm We note that a flow line connecting two critical points $p$ and $q$ is a
curve on the space: in our case, a map $f:\BR\lrt \widetilde{\wp}$ satisfying the
following equations:
\[
\frac{df}{ds}={\rm grad}\,A(f(s)),\quad \lim_{s\rightarrow\infty}f(s)=q,
\quad \lim_{s\rightarrow -\infty}f(s)=p.
\]
We observe that a map $f$ as above is actually a map $f:\BR\times [0,1]\lrt M$ and
that the first equation is the generalized Cauchy-Riemann equation given by
\begin{equation*}
  \frac{\partial f}{\partial s}+J \frac{\partial f}{\partial t}=0.
\end{equation*}
Thus, $f$ is a $J$-holomorphic map called a {\sf $J$-holomophic strip}.

\vskip 2mm
Let $\mathfrak{L}(n)$ be the space of all Lagrangian subspaces of the symplectic vector
space $(\BR^{2n},\omega_0)$ with the standard symplectic form $\omega_0$. It is known that
$\mathfrak{L}(n)$ can be identified with the quotient space $U(n)/O(n).$ We define the map
\begin{equation*}
  \rho: \mathfrak{L}(n)\lrt S^1, \quad \rho(\alpha):=\det (\alpha^2) \ \ (\alpha\in U(n)).
\end{equation*}

\begin{definition}\label{def:6.39}
Let $\gamma$ be a loop in $\mathfrak{L}(n)$. The {\sf Maslov index} of $\gamma$ is defined
to be the degree of the map $\rho\circ \g$. We note that it is independent of the homotopy class
of the loop.
\end{definition}

\begin{proposition}\label{prop:6.40}
By the definition above, the Maslov index may be seen as a cocycle $\mu\in H^1(\mathfrak{L}(n)).$
Then, its Poincar{\'e} dual is given by
$$
\Sigma=\left\{ \mathfrak{L}\in \mathfrak{L}(n)\,\vert \
\dim \mathfrak{L}\cap \mathfrak{L}_0 >0 \right\},
$$
where $\mathfrak{L}_0$ is a fixed Lagrangian submanifold.
\end{proposition}
\begin{proof}
The proof can be found in \cite{FO1, FO2}.
\end{proof}

\begin{definition}\label{def:6.41}
For a $J$-holomorphic strip $u:\BR\times [0,1]\lrt M$, the {\sf Maslov index} of the strip
is defined to be the number of non-transverse intersections (counted with sign and
multiplicity) of $u^*T{\mathfrak{L}_0}\vert_{\BR\times \{0\}}$ and
$u^*T{\mathfrak{L}_{0}}\vert_{\BR\times \{1\}}$ in the Lagrangian Grassmanian after a trivialization
of $u^*TM$.
\end{definition}

\begin{definition}\label{def:6.42}
The {\sf Novikov field} $\Lambda_{\rm nov}$ over a ground field $k$ is given by:
\begin{equation*}
  \Lambda_{\rm nov}:=\left\{ \sum_{i=1}^\infty a_i T^{\beta_i}\ \Big|\ a_i\in k,\ \beta_i\in \BR
  \ {\rm and}\ \lim_{i\rightarrow \infty}\beta_i=\infty \right\}.
\end{equation*}
For given two Lagrangian submanifolds $\mathfrak{L}_0$ and $\mathfrak{L}_1$ in a
symplectic manifold $(M,\omega)$ which intersect transversely, we define the {\sf Floer complex}
\begin{equation*}
  {\rm CF}(\mathfrak{L}_0,\mathfrak{L}_1)=\bigoplus_{q\in \mathfrak{L}_0\cap\mathfrak{L}_1}
  \Lambda_{\rm nov}\,q
\end{equation*}
which is the free $\Lambda_{\rm nov}$--module generated by the intersection points.
For an intersection point $p\in \mathfrak{L}_0\cap\mathfrak{L}_1$, the {\sf Floer boundary operator}
$\partial:{\rm CF}(\mathfrak{L}_0,\mathfrak{L}_1)\lrt
{\rm CF}(\mathfrak{L}_0,\mathfrak{L}_1)$ is defined by
\begin{equation*}
  \partial (p)=\sum_{q\in \mathfrak{L}_0\cap \mathfrak{L}_1}
  \sum_{A\in \pi_2 (M,\mathfrak{L}_0\cup\mathfrak{L}_1)\atop \mu (A)=1}
  \sharp M(p,q,J,A)\cdot T^{\int_A \omega}q.
\end{equation*}
Here,
$M(p,q,J,A)$ is the space consisting of all J-holomorphic maps $u$ joining $p$ and $q$
such that $[u]=A\in \pi_2(M,\mathfrak{L}_0\cup\mathfrak{L}_1).$
\end{definition}

 Floer proved the following theorem:
\begin{theorem}\label{thm:6.43}
  If $M,\ \mathfrak{L}_0$ and $\mathfrak{L}_1$ are symplectically aspherical
  with respect to $\omega,\ \omega\vert_{\mathfrak{L}_0}$ and $\omega\vert_{\mathfrak{L}_1}$
  respectively and $J$ is a generic $t$--dependent family of almost complex structures,
  then $\partial^2=0.$  Here, the {\sf Lagrangian Floer homology} is defined to be
$$
{\rm HF}(\mathfrak{L}_0,\mathfrak{L}_1,J)=
  \frac{{\rm ker}\,\partial}{{\rm im}\,\partial}.
$$
\end{theorem}
\begin{proof}
The proof can be found in \cite{F1+}.
\end{proof}

Floer proved the following theorem:
\begin{theorem}\label{thm:6.44}
  Let $\mathfrak{L}_0$ and $\mathfrak{L}_1$ be Lagrangian submanifolds in a symplectic
  manifold $(M,\omega)$. Assume the same hypothesis of Theorem 6.43. Assume $\varphi:M\lrt M$
  is a {\sf Hamiltonian diffeomorphism}.
  Then we have the following
  properties (a),\,(b) and (c):
\vskip 2mm\noindent
\ {\rm (a)} ${\rm HF}(\mathfrak{L}_0,\mathfrak{L}_1,J)$ does not depend on $J$, and hence
will be denoted by ${\rm HF}(\mathfrak{L}_0,\mathfrak{L}_1)$.
\vskip 2mm\noindent
\ {\rm (b)} ${\rm HF}(\mathfrak{L}_0,\mathfrak{L}_1)=
{\rm HF}(\mathfrak{L}_0,\varphi(\mathfrak{L}_1)).$
\vskip 2mm\noindent
\ {\rm (c)} ${\rm HF}(\mathfrak{L},\mathfrak{L}):=
{\rm HF}(\mathfrak{L},\varphi(\mathfrak{L}))=H_* (\mathfrak{L},\Lambda_{\rm nov})$
for any Lagrangian submanifold $\mathfrak{L}$ in $M$.
\vskip 2mm\noindent
Here, $H_* (\mathfrak{L},\Lambda_{\rm nov})$ is the singular homology with coefficients in
$\Lambda_{\rm nov}$.
\end{theorem}
\begin{proof}
The proof can be found in \cite{F1+}.
\end{proof}

\begin{remark}\label{rk:6.45}
Y.-G. Oh \cite{Oh1, Oh2} generalized Floer's work on the Lagrangian intersection problem.
Indeed, he constructed Lagrangian intersection cohomology for a pair of monotone Lagrangian
submanifolds of a symplectic manifold.
\end{remark}
%%%%%%%%%%%%%%%%%%%%%%%%%%%%%%%%%%%%%%%%%%%%%%%%%%%%%%%%%%%%%%%%%%%%%%%%%%%%%%%%%%%%%%%%%%%%%%%%%%%%%%%
%%%%%%%%%%%%%%%%%%%%%%%%%%%%%%%%%%%%%%%%%%%%%%%%%%%%%%%%%%%%%%%%%%%%%%%%%%%%%%%%%%%%%%%%%%%%%%%%%%%%%%%
%%%%%%%%%%%%%%%%%%%%%%%%%%%%%%%%%%%%%%%%%%%%%%%%%%%%%%%%%%%%%%%%%%%%%%%%%%%%%%%%%%%%%%%%%%%%%%%%%%%%%%%
%%%%%%%%%%%%%%%%%%%%%%%%%%%%%%%%%%%%%%%%%%%%%%%%%%%%%%%%%%%%%%%%%%%%%%%%%%%%%%%%%%%%%%%%%%%%%%%%%%%%%%%
%%%%%%%%%%%%%%%%%%%%%%%%%%%%%%%%%%%%%%%%%%%%%%%%%%%%%%%%%%%%%%%%%%%%%%%%%%%%%%%%%%%%%%%%%%%%%%%%%%%%%%%
%%%%%%%%%%%%%%%%%%%%%%%%%%%%%%%%%%%%%%%%%%%%%%%%%%%%%%%%%%%%%%%%%%%%%%%%%%%%%%%%%%%%%%%%%%%%%%%%%%%%%%%
\vskip 3mm\noindent
{\bf (II) The Fukaya category}
\vskip 2mm
In 1993 Kenji Fukaya\,(1959-- ) introduced the concept of the Fukaya category which is an important
object in the study of symplectic geometry. The Fukaya category encodes many important properties
of a symplectic manifold and Lagrangian submanifolds.
Loosely speaking, the Fukaya category is a $\mathcal{A}_\infty$--category with objects given by
Lagrangian submanifolds and with {\rm hom}--sets given by the Floer complexes such that
the composition maps are given by counting holomorphic strips.
%%%%%%%%%%%%%%%%%%%%%%%%%%%%%%%%%%%%%%%%%%%%%%%%%%%%%%%%%%%%%%%%%%%%%%%%%%%%%%%%%%%%%%%%%%%%%%%%%%%%%%%
%%%%%%%%%%%%%%%%%%%%%%%%%%%%%%%%%%%%%%%%%%%%%%%%%%%%%%%%%%%%%%%%%%%%%%%%%%%%%%%%%%%%%%%%%%%%%%%%%%%%%%%
\vskip 2mm
\begin{definition}\label{def:6.46}
We fix a base field $k$. An $\mathcal{A}_\infty$--{\sf category} $\mathscr{C}$ consists of objects,
graded vector spaces over $k$ denoted by ${\rm hom}_{\mathscr{C}} (X,Y)$ for any pair of
objects $X,Y$ in $\mathscr{C}$, and composition maps
\begin{equation*}
  \mu_{\mathscr{C}}^d:{\rm hom}_{\mathscr{C}} (X_0,X_1)\otimes\cdots\otimes
 {\rm hom}_{\mathscr{C}} (X_{d-1},X_d) \lrt {\rm hom}_{\mathscr{C}} (X_0,X_d)[2-d]
\end{equation*}
for every $d\geq 1,$ where $V[n]$ means that the grading is shifted down by $n$.
If there is no fusion, we write ${\rm hom} (X,Y)$ and $\mu^d$
instead of ${\rm hom}_{\mathscr{C}} (X,Y)$ and $\mu_{\mathscr{C}}^d$ respectively.
\end{definition}
%%%%%%%%%%%%%%%%%%%%%%%%%%%%%%%%%%%%%%%%%%%%%%%%%%%%%%%%%%%%%%%%%%%%%%%%%%%%%%%%%%%%%%%%%%%%%%%%%%%%%%%
%%%%%%%%%%%%%%%%%%%%%%%%%%%%%%%%%%%%%%%%%%%%%%%%%%%%%%%%%%%%%%%%%%%%%%%%%%%%%%%%%%%%%%%%%%%%%%%%%%%%%%%
\vskip 2mm
\begin{definition}\label{def:6.47}
Let $\mathscr{C}$ be an $\mathcal{A}_\infty$--category. The {\sf cohomological category}
$\mathcal{H}(\mathscr{C})$ of $\mathscr{C}$ is defined to be a category with objects given by
the same objects of $\mathscr{C}$, and {\rm hom}-{\rm sets} given by the cohomology groups
$H^*({\rm hom}(X,Y),\mu^1)$ of the complex ${\rm hom}(X,Y)$ for any pair $(X,Y)$ of objects.
The composition on $\mathcal{H}(\mathscr{C})$ is defined by $\mu^2$ as follows:
\begin{equation*}
  [a_1]\circ [a_2]=(-1)^{\deg (a_1)}[\mu^2(a_1,a_2)].
\end{equation*}
The composition is associative.
\end{definition}

%%%%%%%%%%%%%%%%%%%%%%%%%%%%%%%%%%%%%%%%%%%%%%%%%%%%%%%%%%%%%%%%%%%%%%%%%%%%%%%%%%%%%%%%%%%%%%%%%%%%%%%
%%%%%%%%%%%%%%%%%%%%%%%%%%%%%%%%%%%%%%%%%%%%%%%%%%%%%%%%%%%%%%%%%%%%%%%%%%%%%%%%%%%%%%%%%%%%%%%%%%%%%%%
\begin{definition}\label{def:6.48}
For $d+1$ Lagrangian submanifolds $\mathfrak{L}_0,\ldots,\mathfrak{L}_d$ inside a symplectic
manifold, we define
$$
\mu^d:{\rm CF}(\mathfrak{L}_0,\mathfrak{L}_1)\otimes \cdots \otimes
{\rm CF}(\mathfrak{L}_{d-1},\mathfrak{L}_d)\lrt {\rm CF}(\mathfrak{L}_0,\mathfrak{L}_d)
$$
by
\begin{equation*}
  \mu^d(p_1,p_2,\ldots,p_d)=\sum_{q\in \mathfrak{L}_0\cap \mathfrak{L}_d}
  \sum_{A\in \pi_2(M,\mathfrak{L}_1\cup \mathfrak{L}_2 \cup \cdots\cup \mathfrak{L}_d)\atop
  {\rm ind}(A)=2-d}
  \sharp M(p_1,p_2,\ldots,p_d,q,J,A)\cdot
  T^{\int_A\omega}q,
\end{equation*}
where the index is otained by concatenating paths and $M(p_1,p_2,\ldots,p_d,q,J,A)$ is
the manifold of $J$--holomorphic maps with homotopy class $A$ joining
$p_1,p_2,\ldots,p_d$ and $q$.
\end{definition}

\begin{theorem}\label{thm:6.49}
  The maps $\mu^k$ satisfy the $\mathcal{A}_\infty$--relations.
\end{theorem}
\begin{proof}
  The proof can be found in \cite{Se2}.
\end{proof}

\vskip 2mm
Now we define the Fukaya category as follows:
\begin{definition}\label{def:6.50}
Let $(M,\omega)$ be a symplectic manifold with the first class $c_1(TM)=0.$
\vskip 2mm
{\bf $\bullet$}
The objects of the {\sf Fukaya category} $\mathfrak{F}(M,\omega)$ are
compact closed oriented Lagrangian submanifolds $\mathfrak{L}\subset M$
such that $[\omega]\cdot \pi_2 (M,\mathfrak{L})=0$ and with vanishing Maslov class.
\vskip 2mm
{\bf $\circledcirc$}
For every pair of objects $(\mathfrak{L}_0,\mathfrak{L}_1)$, we consider
a family of almost complex structures and Hamiltonian (to achieve transversality)
in order to define ${\rm CF}(\mathfrak{L}_0,\mathfrak{L}_1)$. We also need to
get some perturbation data to have transversality of more than two Lagrangian submanifolds
in order to define the maps $\mu^d.$

\vskip 2mm
{\bf $\bullet$} Given this, we set
$$
{\rm hom}(\mathfrak{L}_0,\mathfrak{L}_1):={\rm CF}(\mathfrak{L}_0,\mathfrak{L}_1)\quad{\rm and}
\quad \mu^1:={\rm the\ Floer\ differential}.
$$
The higher composition maps are given by counts of perturbed $J$--holomorphic disks
as in the definition above. By Theorem 6.49, $\mathfrak{F}(M,\omega)$ becomes
a $\Lambda_{\rm nov}$--linear, $\BZ$--graded, non-unital (but cohomologically unital)
$\mathcal{A}_\infty$--category.
\end{definition}
%%%%%%%%%%%%%%%%%%%%%%%%%%%%%%%%%%%%%%%%%%%%%%%%%%%%%%%%%%%%%%%%%%%%%%%%%%%%%%%%%%%%%%%%%%%%%%%%%%%%%%%
%%%%%%%%%%%%%%%%%%%%%%%%%%%%%%%%%%%%%%%%%%%%%%%%%%%%%%%%%%%%%%%%%%%%%%%%%%%%%%%%%%%%%%%%%%%%%%%%%%%%%%%
%%%%%%%%%%%%%%%%%%%%%%%%%%%%%%%%%%%%%%%%%%%%%%%%%%%%%%%%%%%%%%%%%%%%%%%%%%%%%%%%%%%%%%%%%%%%%%%%%%%%%%%
%%%%%%%%%%%%%%%%%%%%%%%%%%%%%%%%%%%%%%%%%%%%%%%%%%%%%%%%%%%%%%%%%%%%%%%%%%%%%%%%%%%%%%%%%%%%%%%%%%%%%%%
%%%%%%%%%%%%%%%%%%%%%%%%%%%%%%%%%%%%%%%%%%%%%%%%%%%%%%%%%%%%%%%%%%%%%%%%%%%%%%%%%%%%%%%%%%%%%%%%%%%%%%%
%%%%%%%%%%%%%%%%%%%%%%%%%%%%%%%%%%%%%%%%%%%%%%%%%%%%%%%%%%%%%%%%%%%%%%%%%%%%%%%%%%%%%%%%%%%%%%%%%%%%%%%
\vskip 3mm\noindent
{\bf (III) Homological mirror symmetry}
\vskip 2mm

\begin{definition}\label{def:6.51}
A $n$--dimensional compact K{\"a}hler manifold is called a {\sf Calabi--Yau manifold} if its
canonical line bundle is trivial, that is, if there exists a nowhere vanishing holomorphic
$n$--form.
\end{definition}
\vskip 2mm
\begin{definition}\label{def:6.52}
Let $\mathcal{A}$ be an abelian category. The {\sf derived category} {\rm D}($\mathcal{A}$) is defined by a
universal property with respect to the category {\sf Kom}($\mathcal{A}$) of cochain complexes
with terms in $\mathcal{A}$. The objects of {\sf Kom}($\mathcal{A}$) are of the form
\begin{equation*}
  \cdots\lrt X^{-1}\stackrel{d^{-1}}{\lrt} X^{0}\stackrel{d^{0}}{\lrt}
  X^{1}\stackrel{d^{1}}{\lrt}X^{2}\stackrel{d^{2}}{\lrt}\cdots,
\end{equation*}
where each $X^i$ is an object of $\mathcal{A}$ and each of the composites $d^{i+1}\circ d^i$
is zero. The $i$--th cohomology group of the complex is
$H^i(X^{\bullet})=\ker d^i/ {\rm im}\,d^{i-1}.$ If ($X^{\bullet},d_X^{\bullet}$) and
($Y^{\bullet},d_Y^{\bullet}$) are two objects in this category, then a morphism
$f^{\bullet}:(X^{\bullet},d_X^{\bullet})\lrt (Y^{\bullet},d_Y^{\bullet})$ is defined
to be a family of morphisms $f_i:X^i\lrt Y^i$ such that $f_{i+1}\circ d_X^i=d_Y^i\circ f_i$.
Such a morphism induces morphisms on cohomology groups
$H^i(f^{\bullet}):H^i(X^{\bullet})\lrt H^i(Y^{\bullet})$, and $f^{\bullet}$ is called
a {\sf quasi-isomorphism} if each of these morphisms is an isomorphism in $\mathcal{A}$.
\end{definition}
\vskip 2mm
We recall a triangulated category.

\begin{definition}\label{def:6.53}
A {\sf shift} or {\sf translation functor} on a category $D$ is an additive automorphism
$\Sigma$ from $D$ to $D$. It is common to write $X[n]=\Sigma^n X$ for integers $n$.
A {\sf triangle} $(X, Y, Z, u, v, w)$ consists of three objects $X, Y$, and $Z$, together with morphisms
$u:X\lrt Y,\ v:Y\lrt Z$ and $w:Z\lrt X[1].$ Triangles are generally written in the unravelled form:
$$
X\stackrel{u}{\lrt} Y \stackrel{v}{\lrt} Z \stackrel{w}{\lrt} X[1].
$$
A {\sf triangulated category} is an additive category D with a translation functor and a class of triangles, called {\sf exact triangles} (or {\sf distinguished triangles}), satisfying the following properties {\rm (TR\,1)}, {\rm (TR\,2)}, {\rm (TR\,3)} and {\rm (TR\,4)}. (These axioms are not entirely independent, since {\rm (TR\,3)} can be derived from the others.
\vskip 2mm\noindent
{\sf (TR\,1)}
\vskip 2mm\noindent
\ $\blacklozenge$ For every object X, the following triangle is exact:
\[
X\stackrel{{\rm id}}{\lrt} X \lrt 0\lrt X[1]
\]
\vskip 2mm\noindent
\ $\blacklozenge$ For every morphism $u:X\lrt Y$, there is an object $Z$
(called a {\sf cone} or {\sf cofiber} of
\par \
the morphism $u$) fitting into an exact triangle
\[
X\stackrel{{\rm u}}{\lrt} Y \lrt Z\lrt X[1].
\]
\vskip 2mm\noindent
\ $\blacklozenge$ Every triangle isomorphic to an exact triangle is exact.
%%%%%%%%%%%%%%%%%%%%%%%%%%%%%%%%%%%%%%%%%%%%%%%%%%%%%%%%%%%%%%%%%%%%%%%%%%%%%%%%%%%%%%%%%%%%%%%%%%%
%%%%%%%%%%%%%%%%%%%%%%%%%%%%%%%%%%%%%%%%%%%%%%%%%%%%%%%%%%%%%%%%%%%%%%%%%%%%%%%%%%%%%%%%%%%%%%%%%%%
\vskip 3mm\noindent
{\sf (TR\,2)} If
\[
X\stackrel{u}{\lrt} Y\stackrel{v}{\lrt} Z\stackrel{w}{\lrt} X[1]
\]
is an exact triangle, then so are the following rotated triangles
\[
Y\stackrel{v}{\lrt} Z\stackrel{w}{\lrt} X[1]\stackrel{ {-u[1]}}{\lrt} Y[1]
\qquad {\rm and}\qquad
Z[-1]\stackrel{ {-w[-1]}}{\lrt} X\stackrel{u}{\lrt} Y\stackrel{v}{\lrt} Z.
\]
%%%%%%%%%%%%%%%%%%%%%%%%%%%%%%%%%%%%%%%%%%%%%%%%%%%%%%%%%%%%%%%%%%%%%%%%%%%%%%%%%%%%%%%%%%%%%%%%%%%
%%%%%%%%%%%%%%%%%%%%%%%%%%%%%%%%%%%%%%%%%%%%%%%%%%%%%%%%%%%%%%%%%%%%%%%%%%%%%%%%%%%%%%%%%%%%%%%%%%%

\vskip 3mm\noindent
{\sf (TR\,3)} Given two exact triangles and a map between the first morphisms in each triangle, there exists a morphism between the third objects in each of the two triangles that makes everything commute.
%%%%%%%%%%%%%%%%%%%%%%%%%%%%%%%%%%%%%%%%%%%%%%%%%%%%%%%%%%%%%%%%%%%%%%%%%%%%%%%%%%%%%%%%%%%%%%%%%%%
%%%%%%%%%%%%%%%%%%%%%%%%%%%%%%%%%%%%%%%%%%%%%%%%%%%%%%%%%%%%%%%%%%%%%%%%%%%%%%%%%%%%%%%%%%%%%%%%%%%
\vskip 3mm\noindent
{\sf (TR\,4) The octahedral axiom}
\vskip 1mm\noindent
Let $u:X\lrt Y$ and $v:Y\lrt Z$ be morphisms, and consider the composed morphism
$vu:X\lrt Z$. Form exact triangles for each of these three morphisms according to {\sf (TR\,1)}.
The octahedral axiom states (roughly) that the three mapping cones can be made into the vertices
of an exact triangle so that ``everything commutes."
\end{definition}
%%%%%%%%%%%%%%%%%%%%%%%%%%%%%%%%%%%%%%%%%%%%%%%%%%%%%%%%%%%%%%%%%%%%%%%%%%%%%%%%%%%%%%%%%%%%%%%%%%%
%%%%%%%%%%%%%%%%%%%%%%%%%%%%%%%%%%%%%%%%%%%%%%%%%%%%%%%%%%%%%%%%%%%%%%%%%%%%%%%%%%%%%%%%%%%%%%%%%%%
%%%%%%%%%%%%%%%%%%%%%%%%%%%%%%%%%%%%%%%%%%%%%%%%%%%%%%%%%%%%%%%%%%%%%%%%%%%%%%%%%%%%%%%%%%%%%%%%%%%
%%%%%%%%%%%%%%%%%%%%%%%%%%%%%%%%%%%%%%%%%%%%%%%%%%%%%%%%%%%%%%%%%%%%%%%%%%%%%%%%%%%%%%%%%%%%%%%%%%%
%%%%%%%%%%%%%%%%%%%%%%%%%%%%%%%%%%%%%%%%%%%%%%%%%%%%%%%%%%%%%%%%%%%%%%%%%%%%%%%%%%%%%%%%%%%%%%%%%%%
%%%%%%%%%%%%%%%%%%%%%%%%%%%%%%%%%%%%%%%%%%%%%%%%%%%%%%%%%%%%%%%%%%%%%%%%%%%%%%%%%%%%%%%%%%%%%%%%%%%
\vskip 2mm
We recall the notion of a coherent sheaf on a ringed space.
\begin{definition}\label{def:6.54}
A sheaf of $\mathscr{O}_X$--modules $\mathfrak{S}$ over a ringed space $(X,\mathscr{O}_X)$
is called {\sf coherent} if satisfies the following properties {\sf (C1)} and {\sf (C2)}:
\vskip 2mm\noindent
\ {\sf (C1)} $\mathfrak{S}$ is of finite type over $\mathscr{O}_X$;
\vskip 2mm\noindent
\ {\sf (C2)} For any open subset $U\subset X$, any $m\in\BZ^+$ and any morphism
$\varphi:\mathscr{O}_X^m\vert_U\lrt \mathfrak{S}\vert_U$
\par \ \ \ \
of $\mathscr{O}_X$--modules, the kernel of $\varphi$ is of finite type.
\end{definition}

\begin{definition}\label{def:6.55}
Let $X$ be a compact complex manifold with symplectic structure $\omega$. We define $D^\pi (X)$
to be the triangulated derived category of the Fukaya category $\mathfrak{F}(X,\omega)$.
Let {\sf Coh}$(X)$ be the category of coherent sheaves over $X$. It is known that
{\sf Coh}$(X)$ is an abelian category. We define
$D^b(X)$ to be the bounded derived category $D^b$({\sf Coh}$(X))$.
\end{definition}
%%%%%%%%%%%%%%%%%%%%%%%%%%%%%%%%%%%%%%%%%%%%%%%%%%%%%%%%%%%%%%%%%%%%%%%%%%%%%%%%%%%%%%%%%%%%%%%%%%%%%%%
%%%%%%%%%%%%%%%%%%%%%%%%%%%%%%%%%%%%%%%%%%%%%%%%%%%%%%%%%%%%%%%%%%%%%%%%%%%%%%%%%%%%%%%%%%%%%%%%%%%%%%%
The Fukaya category is determined by the symplectic structure of the Calabi-Yau manifold $X$,
while the coherent sheaves base on the complex structure of its mirror $\widehat{X}$.
The correspondence between symplectic geometry and complex geometry might originate from the
perspective of homological algebra and category theory which goes beyond the differential
and algebraic geometry.
%%%%%%%%%%%%%%%%%%%%%%%%%%%%%%%%%%%%%%%%%%%%%%%%%%%%%%%%%%%%%%%%%%%%%%%%%%%%%%%%%%%%%%%%%%%%%%%%%%%%%%%
%%%%%%%%%%%%%%%%%%%%%%%%%%%%%%%%%%%%%%%%%%%%%%%%%%%%%%%%%%%%%%%%%%%%%%%%%%%%%%%%%%%%%%%%%%%%%%%%%%%%%%%
\vskip 2mm
There are two triangulated derived categories $D^\pi (X)$ and $D^b (\widehat{X})$.
We know that $D^\pi (X)$ encodes the symplectic geometry of a symplectic manifold $M$
and $D^b (\widehat{X})$ encodes the complex and algebraic geometry of $X$.
%%%%%%%%%%%%%%%%%%%%%%%%%%%%%%%%%%%%%%%%%%%%%%%%%%%%%%%%%%%%%%%%%%%%%%%%%%%%%%%%%%%%%%%%%%%%%%%%%%%%%%%
%%%%%%%%%%%%%%%%%%%%%%%%%%%%%%%%%%%%%%%%%%%%%%%%%%%%%%%%%%%%%%%%%%%%%%%%%%%%%%%%%%%%%%%%%%%%%%%%%%%%%%%
%%%%%%%%%%%%%%%%%%%%%%%%%%%%%%%%%%%%%%%%%%%%%%%%%%%%%%%%%%%%%%%%%%%%%%%%%%%%%%%%%%%%%%%%%%%%%%%%%%%%%%%
In 1994 Kontsevich \cite{Ko} proposed the following conjecture in his ICM plenary address:
\begin{conjecture}\!\!{\bf [HMS Conjecture]}.\label{conj:6.56}
 Let $X$ be a Calabi-Yau manifold with its mirror
$\widehat{X}$. Then the derived Fukaya category $D^\pi (X)$ of $X$
is equivalent to the bounded derived
category $D^b (\widehat{X})$ of the coherent sheaves of its mirror $\widehat{X}$ as triangulated categories.
\end{conjecture}
This rather abstract conjecture took some time to be appreciated by either mathematicians or
physicists. Nowadays, the HMS conjecture has been proved for some simple cases, for example
the elliptic curve case \cite{P-Z} and the quartic surface \cite{Se2}. Furthermore there are many results concerning other spaces beyond
Calabi-Yau manifolds, as Fano varieties, abelian varieties and others. Today mirror symmetry
has provided a profound insight into symplectic geometry, complex-algebraic geometry, homological
algebra and noncommutative geometry.

\vskip 5mm\noindent
{\bf Acknowledgements}
\vskip 2mm\noindent
I would like to give my deep and hearty thanks to the anonymous referee
for his (or her) careful reading of the manuscript, and insightful-critical comments and
useful suggestions.

\end{section}

%%%%%%%%%%%%%%%%%%%%%%%%%%%%%%%%%%%%%%%%%%%%%%%%%%%%%%%%%%%%%%%%%%%%%%%%%%%%%%%%%%%%%%%%%%%%%%%%%%%%%%%%
%%%%%%%%%%%%%%%%%%%%%%%%%%%%%%%%%%%%%%%%%%%%%%%%%%%%%%%%%%%%%%%%%%%%%%%%%%%%%%%%%%%%%%%%%%%%%%%%%%%%%%%%
%%%%%%%%%%%%%%%%%%%%%%%%%%%%%%%%%%%%%%%%%%%%%%%%%%%%%%%%%%%%%%%%%%%%%%%%%%%%%%%%%%%%%%%%%%%%%%%%%%%%%%%%
%%%%%%%%%%%%%%%%%%%%%%%%%%%%%%%%%%%%%%%%%%%%%%%%%%%%%%%%%%%%%%%%%%%%%%%%%%%%%%%%%%%%%%%%%%%%%%%%%%%%%%%%
%%%%%%%%%%%%%%%%%%%%%%%%%%%%%%%%%%%%%%%%%%%%%%%%%%%%%%%%%%%%%%%%%%%%%%%%%%%%%%%%%%%%%%%%%%%%%%%%%%%%%%%%
%%%%%%%%%%%%%%%%%%%%%%%%%%%%%%%%%%%%%%%%%%%%%%%%%%%%%%%%%%%%%%%%%%%%%%%%%%%%%%%%%%%%%%%%%%%%%%%%%%%%%%%%

\vskip 1cm
\makeatletter
\def\sectionname{}
\renewcommand{\@bibtitlestyle}{%
  \section*{\refname}%
  \markboth{\refname}{\refname}%
}
\makeatother

\end{document}